\documentclass[11pt]{amsart}

\usepackage{epigamath}

\usepackage[english]{babel}

\usepackage[utf8]{inputenc}
\usepackage[T1]{fontenc}
\usepackage{amsmath,amsthm, amssymb,amsfonts,nicefrac}
\usepackage{mathabx} 
\usepackage{graphicx}
\usepackage{textcomp}
\usepackage{graphpap}
\usepackage{caption}
\usepackage{setspace}
\usepackage{xspace}

\usepackage[cmtip]{xy}
\xyoption{pdf} \xyoption{color} \xyoption{all}

\usepackage{tikz}
\usetikzlibrary{cd,arrows,positioning, decorations.pathreplacing}
\tikzset{>=stealth} \tikzcdset{arrow style=tikz}
\tikzset{link/.style={column sep=1.8cm,row sep=0.16cm}}
\tikzset{link2/.style={column sep=0.4cm,row sep=0.1cm}} 
\tikzset{map/.style={row sep=0em, column sep=0em}}
\tikzset{c/.style={every coordinate/.try}}

\newtheorem{thm}{Theorem}[section]
\newtheorem*{thm*}{Theorem}
\newtheorem{cor}[thm]{Corollary}
\newtheorem{lem}[thm]{Lemma}
\newtheorem{prop}[thm]{Proposition}
\theoremstyle{definition}
\newtheorem{defi}[thm]{Definition}

\newtheorem{rmk}[thm]{Remark}
\newtheorem{ex}[thm]{Example}

\tikzset{
hexagon/.pic = {code={ \tikzset{scale = 0.6}
\coordinate (n) at (-1.75,1){};
\coordinate (E1) at (1,1.3){}; \coordinate (E2) at (1.75,0){};
\coordinate (E3) at (1,-1.3){}; \coordinate (E4) at (-1,-1.3){};
\coordinate (E5) at (-1.75,0){}; \coordinate (E6) at (-1,1.3){};
\draw (E1) to (E2)
    (E2) to (E3)
    (E3) to (E4)
    (E4) to (E5)
    (E5) to (E6)
    (E6) to (E1);
\coordinate (D1) at (0,1){}; \coordinate (D2) at (1.2,0.6){};
\coordinate (D3) at (1.2,-0.6){}; \coordinate (D4) at (0,-1){};
\coordinate (D5) at (-1.2,-0.6){}; \coordinate (D6) at (-1.2,0.6){};
\coordinate (D11) at (0.1,1){}; \coordinate (D22) at (1.2,0.5){};
\coordinate (D222) at (1.55,0.3){};
\coordinate (D2222) at (1.2,0.8){};
\coordinate (D33) at (1.1,-0.7){}; \coordinate (D44) at (-0.1,-1){};
\coordinate (D444) at (0.3,-1.2){};
\coordinate (D55) at (-1.2,-0.4){}; \coordinate (D66) at
(-1.2,0.7){};
\coordinate (D666) at (-1.2,0.8){};
 }},
quadrat/.pic = {code={ \tikzset{scale = 0.8}
\coordinate (A1) at (1,1){}; \coordinate (A2) at (1,-1){};
\coordinate (A3) at (-1,-1){}; \coordinate (A4) at (-1,1){}; \draw
(A1) to (A2)
    (A2) to (A3)
    (A3) to (A4)
    (A4) to (A1);
\coordinate (B1) at (0,0.9){}; \coordinate (B2) at (0.9,0){};
\coordinate (B3) at (0,-0.9){}; \coordinate (B4) at (-0.9,0){}; }},
hexagon-leer/.pic = {code={ \tikzset{scale = 0.8}
\coordinate (n) at (-1.75,1){};
\coordinate (E1) at (1,1.3){}; \coordinate (E2) at (1.75,0){};
\coordinate (E3) at (1,-1.3){}; \coordinate (E4) at (-1,-1.3){};
\coordinate (E5) at (-1.75,0){}; \coordinate (E6) at (-1,1.3){};
\coordinate (D1) at (0,1){}; \coordinate (D2) at (1.2,0.6){};
\coordinate (D3) at (1.2,-0.6){}; \coordinate (D4) at (0,-1){};
\coordinate (D5) at (-1.2,-0.6){}; \coordinate (D6) at (-1.2,0.6){};
\coordinate (D11) at (0.1,1){}; \coordinate (D22) at (1.2,0.5){};
\coordinate (D33) at (1.1,-0.7){}; \coordinate (D44) at (-0.1,-1){};
\coordinate (D55) at (-1.2,-0.4){}; \coordinate (D66) at
(-1.2,0.7){}; }},
quadrat-leer/.pic = {code={ \tikzset{scale = 0.8} 
\coordinate (A1) at (1,1){}; \coordinate (A2) at (1,-1){};
\coordinate (A3) at (-1,-1){}; \coordinate (A4) at (-1,1){};
\coordinate (B1) at (0,0.9){}; \coordinate (B2) at (0.9,0){};
\coordinate (B3) at (0,-0.9){}; \coordinate (B4) at (-0.9,0){};
}} }

\renewcommand{\to}{\longrightarrow}
\newcommand{\rat}{\dashrightarrow}

\newcommand{\simeqv}{\mathrel{\rotatebox[origin=c]{-90}{$\simeq$}}}

\newcommand{\Z}{\ensuremath{\mathbb{Z}}}
\newcommand{\Q}{\ensuremath{\mathbb{Q}}}
\newcommand{\F}{\ensuremath{\mathbb{F}}}

\newcommand{\R}{\ensuremath{\mathbb{R}}}
\newcommand{\C}{\ensuremath{\mathbb{C}}}
\newcommand{\p}{\ensuremath{\mathbb{P}}}
\newcommand{\A}{\ensuremath{\mathbb{A}}}

\newcommand{\Ql}{\mathcal{Q}}
\newcommand{\Rl}{\mathcal{R}}
\newcommand{\Sl}{\ensuremath{\mathcal{S}}}

\DeclareMathOperator{\PGL}{PGL}
\newcommand{\SlO}{{\mathrm{SO}}}
\DeclareMathOperator{\Aut}{Aut} \DeclareMathOperator{\Bir}{Bir}
\newcommand{\Pic}{{\mathrm{Pic}}}
\newcommand{\NS}{{\mathrm{NS}}}
\newcommand{\rk}{{\mathrm{rk}}}
\DeclareMathOperator{\GL}{GL} 
\DeclareMathOperator{\Gal}{Gal} 

\DeclareMathOperator{\sym}{Sym} 
 \DeclareMathOperator{\pt}{pt}

\renewcommand{\k}{\mathbf{k}}

\newcommand{\bk}{\overline{\mathbf{k}}}

\def\dashmapsto{\mapstochar\dashrightarrow}

 \DeclareMathOperator{\id}{id}

\newcommand{\supth}[1]{\ensuremath{#1^{\mathrm{th}}}}

\EpigaVolumeYear{5}{2021} \EpigaArticleNr{14} \ReceivedOn{August 14,
2020}
\InFinalFormOn{July 20, 2021} \AcceptedOn{August 22, 2021}

\title{Algebraic subgroups of the plane Cremona group \\ over a perfect field}
\titlemark{Algebraic subgroups of the plane Cremona group}

\author{Julia Schneider}
\address{\'Ecole Polytechnique F\'ed\'erale de Lausanne, chair of Algebraic Geometry (B\^atiment MA), Station 8, CH-1015 Lausanne, Switzerland}
\email{julia.schneider@epfl.ch}

\author{Susanna Zimmermann}
\address{Institut de math\'ematiques d'Orsay, Universit\'e Paris-Saclay, Orsay Ville, France}
\email{susanna.zimmermann@universite-paris-saclay.fr}
  
\authormark{J. Schneider and S. Zimmermann}

\AbstractInEnglish{We show that any infinite algebraic subgroup of the
  plane Cremona group over a perfect field is contained in a maximal
  algebraic subgroup of the plane Cremona group. We classify the
  maximal groups, and their subgroups of rational points, up to
  conjugacy by a birational map. In the last two sections, added in
  2023, we correct the missing case in the classification of infinite
  algebraic subgroups of the plane Cremona group over a perfect
  field.
  
  \medskip\noindent \textbf{[\emph{Corrigendum added on September 2023}]}}



\TitleInFrench{Sous-groupes alg\'ebriques du groupe de Cremona du
plan sur un corps parfait}

\AbstractInFrench{Nous montrons que tout sous-groupe alg\'ebrique
infini du groupe de Cremona du plan sur un corps parfait est contenu
dans un sous-groupe maximal. Nous classifions les
sous-groupes maximaux ainsi que leurs sous-groupes des points
rationnels, \`a conjugaison par une transformation birationnelle
pr\`es.
Dans les deux derni\`eres sections, ajout\'ees en 2023, nous corrigeons le cas manquant dans la classification des sous-groupes alg\'ebriques infinis du groupe de Cremona du plan sur un corps parfait.}

\acknowledgement{J.~S. is supported by the Swiss National Science
Foundation project P2BSP2\_200209 and hosted by the Institut de
Math\'ematiques de Toulouse. S.~Z. is supported by the ANR Project FIBALGA ANR-18-CE40-0003-01,
the Projet PEPS 2019 ``JC/JC'' and the Project ``\'Etoiles montantes
of the R\'egion Pays de la Loire''.}

\begin{document}


\removeabove{0.2cm} \removebetween{0.2cm} \removebelow{0.2cm}

\maketitle

\begin{prelims}

\DisplayAbstractInEnglish

\bigskip

\DisplayKeyWords

\medskip

\DisplayMSCclass

\bigskip

\languagesection{Fran\c{c}ais}

\bigskip

\DisplayTitleInFrench

\medskip

\DisplayAbstractInFrench

\end{prelims}


\newpage

\setcounter{tocdepth}{1}

\tableofcontents


\section{Introduction}

We study algebraic groups acting birationally and faithfully on a
rational smooth projective surface over a perfect field $\k$. Any
choice of birational map from that surface to the projective plane
$\p^2$ induces an action of the algebraic group on $\p^2$ by
birational transformations. Its subgroup of rational points can thus
be viewed as a subgroup of the plane Cremona group $\Bir_\k(\p^2)$,
which motivates the name {\em algebraic subgroup} of
$\Bir_\k(\p^2)$. The full classification - up to conjugacy - of
algebraic subgroups of the plane Cremona group is open over many
fields, because classifying the finite algebraic groups is very
hard. Here is a selection of classification results over various
perfect fields: \cite{BeauvilleBayle, Blanc-prime, Blanc-conj-cl,
MR2641179, BlancLinearisation, DIprime, Robayo, Yasinsky2,Yasinsky}.
The full classification of maximal algebraic subgroups of
$\Bir_\C(\p^2)$ (finite and infinite) can be found in
\cite{Blanc_alg_subgroups} and the classification of the real locus
of infinite algebraic subgroups of $\Bir_\R(\p^2)$ can be found in
\cite{RZ}. In this article, we restrict ourselves to consider
infinite algebraic subgroups of $\Bir_\k(\p^2)$ over a perfect field
$\k$ and we classify these groups up to conjugacy by elements of
$\Bir_\k(\p^2)$ and up to inclusion. We also classify their
subgroups of $\k$-rational points up to conjugation by elements of
$\Bir_\k(\p^2)$ and up to inclusion. The two classifications are
different as soon as $\k$ has a quadratic extension, see
Corollary~\ref{cor:1}(\ref{cor1:3})--(\ref{cor1:4}).

Let us explain why we work over a perfect field. Given an algebraic
subgroup $G$ of $\Bir_\k(\p^2)$, the strategy is to find a rational,
regular and projective surface on which $G$ acts by automorphisms
and then use a $G$-equivariant Minimal Model Program to arrive on a
conic fibration or a del Pezzo surface. It then remains to describe
the automorphism group of that surface. Over a perfect field $\k$,
regular implies smooth, and a smooth projective surface over $\k$ is
a smooth projective surface over the algebraic closure $\bk$ of $\k$
equipped with an action of the Galois group $\Gal(\bk/\k)$ of $\bk$
over $\k$. In particular, the classification of rational smooth del
Pezzo surfaces is simply the classification of
$\Gal(\bk/\k)$-actions on smooth del Pezzo surfaces over $\bk$ with
$\Gal(\bk/\k)$-fixed points. This is straightforward if they have
degree $\geq6$, as we will see in \S\ref{sec:DP89} and
\S\ref{sec:DP6}. Over an imperfect field, regular does not imply
smooth and a finite field extension may make appear singularities.
The classification of regular del Pezzo surfaces is still open. In
characteristic $2$, there are regular, geometrically non-normal del
Pezzo surfaces of degree $6$ \cite[Proposition 14.3, Proposition
14.5]{FanelliSchroer20} and there are regular del Pezzo surfaces of
degree $2$ that are geometrically non-reduced \cite[Proposition
3.4.1]{Maddock}. In particular, we cannot use directly the
classification of regular del Pezzo surfaces over a separably closed
field to describe the automorphism group of regular del Pezzo
surfaces over an imperfect field, nor directly the classification of
non-normal del Pezzo surfaces given in \cite{Reid94}. \smallskip

Now, assume again that $\k$ is a perfect field. Theorem~\ref{thm:1},
Theorem~\ref{thm:2}, Theorem~\ref{thm:3} and Corollary~\ref{cor:1}
recover the classification results of \cite{Blanc_alg_subgroups} and
\cite{RZ} over $\C$ and $\R$ for infinite algebraic subgroups, and
we will see that these results extend without any surprises over a
perfect field with at least three elements. We leave it up to the
reader to decide how surprising they find the results over the field
with two elements.

By a theorem of Rosenlicht and Weil, for any algebraic subgroup $G$
of $\Bir_\k(\p^2)$ there is a birational map $\p^2\rat X$ to a
smooth projective surface $X$ on which $G$ acts by automorphisms,
see Proposition~\ref{thm:projective model}. It conjugates $G$ to a
subgroup of $\Aut(X)$, the group scheme of automorphisms of $X$, and
$G(\k)$ is conjugate to a subgroup of $\Aut_\k(X)$. For a conic
fibration $\pi\colon X\to\p^1$ we denote by
$\Aut(X,\pi)\subset\Aut(X)$ the subgroup preserving the conic
fibration, by $\Aut(X/\pi)\subset\Aut(X,\pi)$ its subgroup inducing
the identity on $\p^1$, and by $\Aut_\k(X,\pi)$ and $\Aut_\k(X/\pi)$
their $\k$-points. For a $\Gal(\bk/\k)$-invariant collection  $p_1,\dots,p_r\in X(\bk)$ of
points, we denote by
$\Aut_\k(X,p_1,\dots,p_r)$, resp. $\Aut_\k(X,\{p_1,\dots,p_r\})$,
the subgroup of $\Aut_\k(X)$ fixing each $p_i$, resp.  preserving
the set $\{p_1,\dots,p_r\}$. A {\em splitting field} of
$\{p_1,\dots,p_r\}$ is a finite normal extension $L/\k$ of smallest degree
such that $p_1,\dots,p_r\in X(L)$ and such that $\{p_1,\dots,p_r\}$
is a union of $\Gal(L/\k)$-orbits.

Suppose that $\k$ has a quadratic extension $L/\k$ and let $g$ be
the generator of $\Gal(L/\k)\simeq\Z/2$. By $\Ql^L$ we denote the
$\k$-structure on $\p^1_L\times\p^1_L$ given by $(x,y)^g=(y^g,x^g)$. By
$\Sl^{L,L'}$ we denote a surface obtained by blowing up $\Ql^L$ in a
point $p$ of degree $2$, where $L'/\k$ is the splitting field of
$p$, whose geometric components are not on the same ruling of
$\p^1_L\times\p^1_L$. We will show in Lemma~\ref{prop:DP9} that its isomorphism class depends only on the isomorphism classes of $L,L'$. In Theorem~\ref{thm:1}(\ref{1:62}), we denote
by $E\subset\Sl^{L,L'}$ its exceptional divisor.

\begin{thm}\label{thm:1}
Let $\k$ be a perfect field and $G$ an infinite algebraic subgroup
of $\Bir_\k(\p^2)$. Then there is a $\k$-birational map $\p^2\rat X$
that conjugates $G$ to a subgroup of $\Aut(X)$, with $X$ one of the
following surfaces, where no indication of the $\Gal(\bk/\k)$-action
means the canonical action.
\begin{enumerate}
\item\label{1:1} $X=\p^2$ and $\Aut(\p^2)\simeq\PGL_3$
\item\label{1:2} $X=\F_0$ and $\Aut(\F_0)\simeq\Aut(\p^1)^2\rtimes\Z/2\simeq\PGL_2^2\rtimes\Z/2$
\item\label{1:3} $X=\Ql^L$ and  $\Aut(\Ql^L)$ is the $\k$-structure on $\Aut(\p^1_L)^2\rtimes\Z/2$ given by the $\Gal(L/\k)$-action $(A,B,\tau)^g=(B^g,A^g,\tau)$, where $L/\k$ is a quadratic extension.
\item\label{1:4} $X=\F_n$, $n\geq2$, and the action of $\Aut(\F_n)$ on $\p^1$ induces a split exact sequence
\[
 1\to V_{n+1}\to \Aut(\F_n)\to \GL_2/\mu_n\to 1
\]
where $\mu_n=\{a\id\mid a^n=1\}$ and $V_{n+1}$ is a vector space of dimension $n+1$.

\item\label{1:5} $X$ is a del Pezzo surface of degree $6$ with $\NS(X_{\bk})^{\Aut_{\bk}(X)}=1$. The action of $\Aut_{\bk}(X)$ on $\NS(X_{\bk})$ induces the split exact sequence
\[1\rightarrow (\bk^*)^2\to \Aut_{\bk}(X)\to \sym_3\times\Z/2\rightarrow 1.\]
Moreover, we are in one of the following cases.
    \begin{enumerate}
    \item\label{1:5c}
        $\rk\,\NS(X)=1$ and there is a quadratic extension $L/\k$ and a birational morphism $\pi\colon X_L\to\p^2_L$ blowing up a point $p=\{p_1,p_2,p_3\}$ of degree $3$ with splitting field $F$ over $\k$, and one of the following cases holds:
        \begin{enumerate}
        \item\label{1:5c1} $\Gal(F/\k)\simeq\Z/3$ and the action of $\Aut_\k(X)$ on $\NS(X)$ induces the split exact sequence
            \[1\rightarrow \Aut_L(\p^2,p_1,p_2,p_3)^{\pi\Gal(L/\k)\pi^{-1}}\to \Aut_\k(X) \to \Z/6\rightarrow 1,\]
        \item\label{1:5c2}
        $\Gal(F/\k)\simeq\sym_3$ and the action of $\Aut_\k(X)$ on $\NS(X)$ induces the split exact sequence
            \[1\rightarrow \Aut_L(\p^2,p_1,p_2,p_3)^{\pi\Gal(L/\k)\pi^{-1}}\to \Aut_\k(X) \to \Z/2\rightarrow 1,\]
        \end{enumerate}
    \item\label{1:5b} $\rk\,\NS(X)\geq2$, $\rk\,\NS(X)^{\Aut_\k(X)}=1$ and $X$ is one of the following:
        \begin{enumerate}
        \item\label{1:51} $X$ is the blow-up of $\p^2$ in the coordinate points, and the action of $\Aut_\k(X)$ on $\NS(X)$ induces the split exact sequence
    \[ 1\rightarrow (\k^*)^2\to \Aut_\k(X)\to\sym_3\times\Z/2\rightarrow 1.\]
        \item\label{1:52} $X$ is the blow-up of $\F_0$ in a point $p=\{(p_1,p_1),(p_2,p_2)\}$ of degree $2$. The action of $\Aut_\k(X)$ on $\NS(X)$ induces the exact sequence,
    \[1\rightarrow \Aut_\k(\p^1,p_1,p_2)^2\to \Aut_\k(X)\to \sym_3\times\Z/2\rightarrow 1\]
    which is split if $\mathrm{char}(\k)\neq2$.
        \item\label{1:53} $X$ is the blow-up of $\p^2$ in a point $p=\{p_1,p_2,p_3\}$ of degree $3$ with splitting field $L$ such that $\Gal(L/\k)\simeq\Z/3$. The action of $\Aut_\k(X)$ on $\NS(X)$ induces the split exact sequence
    \[1\rightarrow \Aut_\k(\p^2,p_1,p_2,p_3)\to \Aut_\k(X) \to \Z/6\rightarrow 1\]
        \item\label{1:54} $X$ is the blow-up of $\p^2$ in a point $p=\{p_1,p_2,p_3\}$ of degree $3$ with splitting field $L$ such that $\Gal(L/\k)\simeq\sym_3$. The action of $\Aut_\k(X)$ on $\NS(X)$ induces the split exact sequence
    \[1\rightarrow \Aut_\k(\p^2,p_1,p_2,p_3)\to \Aut_\k(X)\to \Z/2\rightarrow 1\]
    where $\Z/2$ is generated by a rotation.
        \end{enumerate}
    \item\label{1:5a} $\rk\,\NS(X)^{\Aut_\k(X)}=2$ and there is a quadratic extension $L/\k$ and a birational morphism $\nu\colon X\to \Ql^L$ contracting two curves onto rational points $p_1,p_2$ or one curve onto a point $\{p_1,p_2\}$ of degree $2$ with splitting field $L'/\k$. The action of $\Aut_\k(X)$ on $\NS(X)$ induces the split exact sequence
    \[1\rightarrow T^{L,L'}(\k)\to\Aut_\k(X)\to \Z/2\times\Z/2\rightarrow1\]
    where $\nu\Aut_\k(X)\nu^{-1}=\Aut_\k(\Ql^L,\{p_1,p_2\})$ and $T^{L,L'}$ is the subgroup of $\Aut_\k(\Ql^L,p_1,p_2)$ preserving the rulings of $\Ql^L_L$.
    \end{enumerate}
\item\label{1:6} $\pi\colon X\to\p^1$ is one of the following conic fibrations with
\[\rk\,\NS(X_{\bk}/\p^1)^{\Aut_{\bk}(X,\pi)}=\rk\,\NS(X/\p^1)^{\Aut_\k(X,\pi)}=1:\]
    \begin{enumerate}
    \item\label{1:61} $X/\p^1$ is the blow-up of points $p_1,\dots,p_r\in\F_n$, $n\geq2$, contained in a section
    $S_n\subset \F_n$ with $S_n^2=n$.
    The geometric components of the $p_i$ are on pairwise distinct geometric fibres and $\sum_{i=1}^r\deg(p_i)=2n$. There are split exact sequences
            \[
\begin{tikzcd}[link2]
& (T_1/\mu_n)\rtimes\Z/2&\Aut(X)&&\\[-5pt]
& \simeqv&\vert\vert&&\\[-5pt]
1\ar[r] & \Aut(X/\pi_X)\ar[r]&\Aut(X,\pi_X)\ar[r] & \Aut(\p^1,\Delta)\ar[r] &1\\
1\ar[r]& \Aut_\k(X/\pi_X)\ar[r]&\Aut_\k(X,\pi_X)\ar[r] & \Aut_\k(\p^1,\Delta)\ar[r] &1\\[-5pt]
& \simeqv&\vert\vert&&\\[-5pt]
& (\k^*/\mu_n(\k))\rtimes\Z/2&\Aut_\k(X)&&
\end{tikzcd}
\]
    where $\Delta=\pi(\{p_1,\dots,p_r\})\subset\p^1$, $T_1$ is the split one-dimensional torus and $\mu_n$ its subgroup of $n^\mathrm{th}$ roots of unity.
    \item\label{1:62} There exist quadratic extensions $L$ and $L'$ of $\k$ such that $X/\p^1$ is the blow-up of $\Sl^{L,L'}$ in points $p_1,\dots,p_r\in E$, $r\geq1$. The $p_i$ are all of even degree, their geometric components are on pairwise distinct geometric components of smooth fibres and each geometric component of $E$ contains half of the geometric components of each $p_i$. There are exact sequences
        \[
\begin{tikzcd}[link2]
& \SlO^{L,L'}\rtimes\Z/2&\Aut(X)&&\\[-5pt]
& \simeqv&\vert\vert&&\\[-5pt]
1\ar[r] & \Aut(X/\pi_X)\ar[r]&\Aut(X,\pi_X)\ar[r] & \Aut(\p^1,\Delta)\ar[r] &1\\
1\ar[r]& \Aut_\k(X/\pi_X)\ar[r]&\Aut_\k(X,\pi_X)\ar[r] & (D^{L,L'}_\k\rtimes\Z/2)\cap\Aut_\k(\p^1,\Delta)\ar[r] &1\\[-5pt]
& \simeqv&\vert\vert&&\\[-5pt]
& \SlO^{L,L'}(\k)\rtimes\Z/2&\Aut_\k(X)&&
\end{tikzcd}
\]
with $\Delta=\pi(\{p_1,\dots,p_r\})\subset\p^1$ and
$\SlO^{L,L'}=\{(a,b)\in T^L\mid ab=1\}$, and
\begin{itemize}
\item if $L,L'$ are $\k$-isomorphic, then $\SlO^{L,L'}(\k)\simeq\{a\in L^*\mid aa^g=1\}$\\
 and $D^{L,L'}_\k\simeq\{\alpha\in k^*\mid \alpha=\lambda\lambda^g,\lambda\in L\}$, where $g$ is the generator of $\Gal(L/\k)$,
\item if $L,L'$ are not $\k$-isomorphic, then
$\SlO^{L,L'}(\k)\simeq\k^*$ and \\
$D^{L,L'}_{\k}\simeq\{\lambda\lambda^{gg'}\in F\mid \lambda\in
K,\lambda\lambda^{g'}=1\}$, where $\k\subset F\subset LL'$ is the
intermediate extension such that $\Gal(F/\k)\simeq\langle
gg'\rangle\subset\Gal(L/\k)\times\Gal(L'/\k)$, where $g,g'$ are the
generators of $\Gal(L/\k),\Gal(L'/\k)$, respectively.
\end{itemize}
    \end{enumerate}
\end{enumerate}
\end{thm}

We consider a family among
(\ref{1:3}), (\ref{1:5a}), (\ref{1:5c}), (\ref{1:52}), (\ref{1:53}), (\ref{1:54}), and (\ref{1:62})
empty if the point of requested degree or the requested field
extension does not exist.

Theorem~\ref{thm:1}(\ref{1:5}) is in fact the classification of
rational del Pezzo surfaces of degree $6$ over $\k$ up to
isomorphism, and for any of the eight classes there is a field over
which a surface in the class exists, see \S\ref{sec:DP6}.

The next theorem lists the conjugacy classes in $\Bir_\k(\p^2)$ of
the groups in Theorem~\ref{1:1}. Let $G$ be an affine algebraic group and
$X/B$ a $G$-Mori fibre space (see Definition~\ref{def:GMfs}). We
call it {\em $G$-birationally rigid} if for any $G$-equivariant
birational map $\varphi\colon X\rat X'$ to another $G$-Mori fibre
space $X'/B'$ we have $X'\simeq X$. In particular,
$\varphi\Aut(X)\varphi^{-1}=\Aut(X')$. We call it {\em
$G$-birationally superrigid} if any $G$-equivariant birational map
$X\rat X'$ to another $G$-Mori fibre space $X'/B'$ is an
isomorphism. If we replace $G$ by $G(\k)$ everywhere, we get the
notion of $G(\k)$-Mori fibre space, $G(\k)$-birationally rigid and
$G(\k)$-birationally superrigid. The following theorem also shows
that $G$-birationally (super)rigid does not imply
$G(\k)$-birationally (super)rigid.

The del Pezzo surfaces $X$ and the conic fibrations $X/\p^1$ in
Theorem~\ref{thm:1} are $\Aut(X)$-Mori fibre spaces, and, except for
the del Pezzo surfaces from (\ref{1:5a}), they are also
$\Aut_\k(X)$-Mori fibre spaces.

\begin{thm}\label{thm:2}
Let $\k$ be a perfect field.
\begin{enumerate}
\item\label{2:1} Any del Pezzo surface $X$ and any conic fibration $X/\p^1$ from Theorem~\ref{thm:1} is $\Aut(X)$-birationally superrigid.
\item\label{2:2} Any del Pezzo surface $X$ in Theorem~\ref{thm:1}(\ref{1:1})--(\ref{1:4}), (\ref{1:52})--(\ref{1:54}) and any conic fibration $X/\p^1$ from (\ref{1:62}) is $\Aut_\k(X)$-birationally superrigid.
\item\label{2:5} Let $X$ be a del Pezzo surface from Theorem~\ref{thm:1}(\ref{1:5c}). \\
If $|\k|\geq3$, then $X$ is $\Aut_\k(X)$-birationally superrigid. \\
If $|\k|=2$, then there is an $\Aut_\k(X)$-equivariant birational
map $X\rat X'$, where $X'$ is the del Pezzo surface from
Theorem~\ref{thm:1}(\ref{1:52}).
\item\label{2:3} Let $X$ be the del Pezzo surface from Theorem~\ref{thm:1}(\ref{1:51}). \\
If $|\k|\geq3$, then $X$ is $\Aut_\k(X)$-birationally superrigid.
\\
If $|\k|=2$, there are $\Aut_\k(X)$-equivariant birational maps
$X\rat\F_0$ and $X\rat X'$, where $X'$ is the del Pezzo surface of
degree $6$ from Theorem~\ref{thm:1}(\ref{1:52}).
\item\label{2:4} Any conic fibration $X/\p^1$ from Theorem~\ref{thm:1}(\ref{1:61}) is $\Aut_\k(X)$-birationally superrigid if $\k^*/\mu_n(\k)$ is non-trivial. If $\k^*/\mu_n(\k)$ is trivial and $X\rat Y$ is an $\Aut_\k(X)$-equivariant birational map to a surface $Y$ from Theorem~\ref{thm:1}, then $Y\simeq X$.
\end{enumerate}
\end{thm}

We say that an algebraic subgroup $G$ of $\Bir_\k(\p^2)$ is {\em
maximal} if it is maximal with respect to inclusion among the
algebraic subgroups of $\Bir_\k(\p^2)$. We say that $G(\k)$ is {\em
maximal} if for any algebraic subgroup $G'$ of $\Bir_\k(\p^2)$
containing $G(\k)$, we have $G(\k)=G'(\k)$.

By Theorem~\ref{thm:2}(\ref{2:3}), if $|\k|=2$ and $X$ is a del
Pezzo surface from (\ref{1:51}), then $\Aut_\k(X)$ is not maximal:
It is conjugate to a subgroup of $\Aut_\k(\F_0)$ and this inclusion
is strict, because $\Aut_\k(X)\simeq\sym_3\times\Z/2$ has $12$
elements, whereas $\Aut_\k(\F_0)$ has $72$ elements. Similarly,
$\Aut_\k(X)$ is not maximal if $X$ is a del Pezzo surface from
(\ref{1:5c}) and $|\k|=2$.

\begin{cor}\label{cor:1}
Let $\k$ be a perfect field and $H$ an infinite algebraic subgroup
of $\Bir_\k(\p^2)$.
    \begin{enumerate}
    \item\label{cor1:1} Then $H$ is contained in a maximal algebraic subgroup $G$ of $\Bir_\k(\p^2)$.
    \item\label{cor1:3} Up to conjugation by a birational map, the maximal infinite algebraic subgroups of $\Bir_\k(\p^2)$ are precisely the groups $\Aut(X)$ in Theorem~\ref{thm:1}.
    Two maximal infinite subgroups $\Aut(X)$ and $\Aut(X')$ are conjugate by a biratonal map if and only if $X\simeq X'$.
    \item\label{cor1:4} $H(\k)$ is maximal if and only if it is conjugate to one of the $\Aut_\k(X)$ from
        \begin{itemize}
        \item (\ref{1:1})--(\ref{1:4}), (\ref{1:52})--(\ref{1:54}), (\ref{1:6}),
        \item (\ref{1:5c}), (\ref{1:51}) if $|\k|\geq3$.
        \end{itemize}
    Two such groups $\Aut_\k(X)$ and $\Aut_\k(X')$ are conjugate by a birational map if and only if $X\simeq X'$.
    \end{enumerate}
\end{cor}

\begin{thm}\label{thm:3}
Let $\k$ be a perfect field. The conjugacy classes of the maximal
subgroups $\Aut_\k(X)$ of $\Bir_\k(\p^2)$  from Theorem~\ref{thm:1}
are parametrised by
    \begin{itemize}
    \item (\ref{1:1}), (\ref{1:2}): one point
    \item (\ref{1:3}): one point for each $\k$-isomorphism class of quadratic extensions of $\k$
    \item (\ref{1:4}): one point for each $n\geq2$
    \item (\ref{1:5c1}) one point for any pair $(L,F)$ of $\k$-isomorphism classes of quadratic extensions $L$ and Galois extensions $F/\k$ with $\Gal(F/\k)\simeq\Z/3$ if $|\k|\geq3$
    \item (\ref{1:5c2}): one point for any pair $(L,F)$ of $\k$-isomorphism classes of quadratic extensions $L$ and Galois extensions $F/\k$ with $\Gal(F/\k)\simeq\sym_3$
    \item (\ref{1:51}): one point if $|\k|\geq3$
    \item (\ref{1:52}): one point for each $\k$-isomorphism class of quadratic extensions of $\k$
    \item (\ref{1:53}): one point for each $\k$-isomorphism class of Galois extensions $F/\k$ with $\Gal(F/\k)\simeq\Z/3$.
    \item (\ref{1:54}): one point for any $\k$-isomorphism class of Galois extensions $F/\k$ with $\Gal(F/\k)\simeq\sym_3$.
    \item (\ref{1:61}): for each $n\geq2$ the set of points $\{p_1,\dots,p_r\}\subset\p^1$ with $\sum_{i=1}^r\deg(p_i)=2n$ up to the action of $\Aut_\k(\p^1)$
    \item (\ref{1:62}): for each $n\geq1$ and for each pair of $\k$-isomorphism classes of quadratic extensions $(L,L')$, the set of points $\{p_1,\dots,p_r\}\subset\p^1$ of even degree with $\sum_{i=1}^r\deg(p_i)=2n$ up to the action of $D^{L,L'}_{\k}(\k)\rtimes\Z/2$
    \end{itemize}
\end{thm}

We show the following consequence of \cite{Schneider} and
\cite{Zimmermannb, Zimmermann}.

\begin{prop}\label{prop:juliaz}
For any perfect field $\k$ there is a surjective homomorphism
\[\Phi\colon\Bir_\k(\p^2)\to\bigast_J\bigoplus_I\Z/2,\]
where $J$ is the set of points of degree $2$ in $\p^2$ up to
$\Aut_\k(\p^2)$ and $I$ is at least countable. If $[\bk:\k]=2$, then
$|I|=|\k|$.
\end{prop}

If $\k=\R$ (or more generally $[\bk:\k]=2$) then the abelianisation
map of $\Bir_\R(\p^2)$ is a homomorphism as in
Proposition~\ref{prop:juliaz}. By \cite[Theorem 1.3]{RZ} any
infinite algebraic group acting on $\Bir_\R(\p^2)$ that has
non-trivial image in the abelianisation is a subgroup of the group
in (\ref{1:62}), and this holds also if $[\bk:\k]=2$. We will show a
slightly more general statement over perfect fields with
$[\bk:\k]>2$, for which we need to introduce equivalence classes of
Mori fibre spaces and links of type II.

We call two Mori fibre spaces $X_1/\p^1$ and $X_2/\p^1$ equivalent
if there is a birational map $X_1\rat X_2$ that preserves the
fibration. In particular, if $\varphi\colon X_1\rat X_2$ is a link
of type II between Mori fibre spaces $X_1/\p^1$ and $X_2/\p^1$, then
these two are equivalent. There is only one class of Mori fibre
spaces birational to the Hirzebruch surface $\F_1$
\cite[Lemma]{Schneider}, because all rational points in $\p^2$ are
equivalent up to $\Aut(\p^2)$. We denote by $J_6$ the set of classes
of Mori fibre spaces birational to some $\Sl^{L,L'}$, and by $J_5$
the set of classes birational to a blow-up of $\p^2$ in a point of
degree $4$ whose geometric components are in general position. We
call two Sarkisov links $\varphi$ and $\varphi'$ of type II between
conic fibrations equivalent if the conic fibrations are equivalent
and if the base-points of $\varphi$ and $\varphi'$ have the same
degree. For a class $C$ of equivalent rational Mori fibre spaces, we
denote by $M(C)$ the set of equivalence classes of links of type II
between conic fibrations in the class $C$ whose base-points have
degree $\geq16$.

\begin{thm}[{\cite[Theorem 3, Theorem 4]{Schneider}}]\label{thm:julia}
For any perfect field with $[\bk:\k]>2$ there is a non-trivial
homomorphism
\begin{equation}\tag{$\ast$}\label{quotient}
\Psi\colon\Bir_\k(\p^2)\to\bigoplus_{\chi\in
M(\F_1)}\Z/2\ast(\bigast_{C\in J_6}\bigoplus_{\chi\in
M(C)}\Z/2)\ast(\bigast_{C\in J_5}\bigoplus_{\chi\in M(C)}\Z/2).
\end{equation}
\end{thm}
In fact, the homomorphism from Proposition~\ref{prop:juliaz} for
$[\bk:\k]>2$ is induced by the one in Theorem~\ref{thm:julia}.

We show that an infinite algebraic group acting birationally on
$\p^2$ is killed by the homomorphism $\Psi$ unless it is
conjugate to a group of automorphisms acting on $\Sl^{L,L'}$ or a
Hirzebruch surface.

\begin{prop}\label{thm:4}
Let $\k$ be a perfect field with $[\bk:\k]>2$ and let $\Psi$ be the
homomorphism (\ref{quotient}). Let $G$ be an infinite algebraic
subgroup of $\Bir_\k(\p^2)$. Then $\Psi(G(\k))$ is of order at most
$2$ and the following hold.
\begin{enumerate}
\item\label{3:2}
If $\Psi(G(\k))$ is non-trivial, it is contained in the factor
indexed by $\F_1$ or $C\in J_6$ and there is a $G$-equivariant
birational map $\p^2\rat X$ that conjugates $G$ to a subgroup of
$\Aut(X)$, where $X$ is as in Theorem~\ref{1:1}(\ref{1:61}) or
(\ref{1:62}), respectively.
\item\label{3:3}
Let $X/\p^1$ be a conic fibration as in
Theorem~\ref{thm:1}(\ref{1:6}), which is the blow-up of $\F_n$,
$n\geq2$, or $\Sl^{L,L'}$ in points $p_1,\dots,p_r$. If
$\Psi(\Aut_\k(X))$ is non-trivial, it is generated by the element
whose non-zero entries are indexed by the $\chi_i$ that have $p_i$
as base-point, where $i\in\{1,\ldots,r\}$ is such that
$\deg(p_{i})\geq16$ and $|\{j\in\{1,\ldots,r\}\mid
\deg(p_j)=\deg(p_{i})\}|$ is odd.
\end{enumerate}
\end{prop}

The analogous statement to Proposition~\ref{thm:4} with the homomorphism
from Proposition~\ref{prop:juliaz} for a perfect field $\k$ such
that $[\bk:\k]=2$ can be found in \cite{RZ}.

\subsection*{Acknowledgements}
The authors would like to thank Andrea Fanelli for interesting discussions on
algebraic groups over perfect and imperfect fields, and Michel Brion
for his comments on regularisation of birational group actions. They
would like to thank the first referee for his careful reading and
for pointing out an issue in earlier versions of
Lemma~\ref{lem:size} and Lemma~\ref{lem:linkscbF}, which lead to the
addition of Theorem~\ref{thm:2}(\ref{2:5}).

\section{Surfaces and birational group actions}

\subsection{Birational actions}

Throughout the article, $\k$ denotes a perfect field and $\bk$ an
algebraic closure. By a surface $X$ (or $X_{\k}$) we mean a smooth
projective surface over $\k$ such that
$X_{\bk}:=X\times_{\mathrm{Spec}(\k)}\mathrm{Spec}(\bk)$ is
irreducible. We denote by $X(\k)$ the set of $\k$-rational points of
$X$. The Galois group $\Gal(\bk/\k)$ acts on
$X\times_{\mathrm{Spec}(\k)}\mathrm{Spec}(\bk)$ through the second
factor. By a point of degree $d$ we mean a $\Gal(\bk/\k)$-orbit
$p=\{p_1,\dots,p_d\}\subset X(\bk)$ of cardinality $d\geq1$. The
points of degree one are precisely the $\k$-rational points of $X$.
Let $L/\k$ be an algebraic extension of $\k$ such that all $p_i$ are
$L$-rational points. By the blow up of $p$ we mean the blow up of
these $d$ points, which is a morphism $\pi\colon X'\rightarrow X$
defined over $\k$, with exceptional divisor $E=E_1+\cdots +E_d$
where the $E_i$ are disjoint $(-1)$-curves defined over $L$, and
$E^2=-d$. We call $E$ the exceptional divisor of $p$. More
generally, a birational map $f\colon X\rat X'$ is defined over $\k$
if and only if the birational map $f\times\id\colon X_{\bk}\rat
X_{\bk}'$ is $\Gal(\bk/\k)$-equivariant. In particular, $X\simeq X'$
if and only if there is a $\Gal(\bk/\k)$-equivariant isomorphism
$X_{\bk}\to X'_{\bk}$ (see also \cite[\S2.4]{BorelSerre}).

The surface $X$ being projective and geometrically irreducible implies $\k[X_{\bk}]^*=(\bk)^*$, so
if $X(\k)\neq\emptyset$ we have
$\Pic(X_\k)=\Pic(X_{\bk})^{\Gal(\bk/\k)}$ \cite[Lemma
6.3(iii)]{Sansuc}. This holds in particular if $X$ is $\k$-rational,
because then it has a $\k$-rational point by the Lang-Nishimura
theorem. Since numerical equivalence is $\Gal(\bk/\k)$-stable, also
algebraic equivalence is, and hence
$\NS(X_\k)=\NS(X_{\bk})^{\Gal(\bk/\k)}$. The $\Gal(\bk/\k)$-action
on $\NS(X_{\bk})$ factors through a finite group, that is, its
action factors through a finite group. Indeed, since $\Gal(\bk/\k)$
has only finite orbits on $\bk$, the orbit of any prime divisor of
$X_{\bk}$ is finite. Then each generator of the finitely generated
$\Z$-module $\NS(X_{\bk})$ has a finite $\Gal(\bk/\k)$-orbit, so the
action of $\Gal(\bk/\k)$ on the (finite) union of these orbits
factors through a finite group.

If not mentioned otherwise, any surface, curve, point and rational
map will be defined over the perfect field $\k$. By a geometric
component of a curve $C$ (resp. a point $p=\{p_1,\dots,p_d\}$), we
mean an irreducible component of $C_{\bk}$ (resp. one of
$p_1,\dots,p_d$).

By Ch\^{a}telet's theorem, for $n\geq1$ any smooth projective space
$X$ over $\k$ with $X(\k)\neq\emptyset$ such that
$X_{\bk}\simeq\p^n_{\bk}$ is in fact already isomorphic to $\p^n$
over $\k$. This means in particular that $\p^2$ is the only rational
del Pezzo surface of degree $9$ and that a smooth curve of genus $0$
with rational points is isomorphic to $\p^1$.

For a surface $X$, we denote by $\Bir_\k(X)$ its group of birational
self-maps and by $\Aut_{\k}(X)$ the group of $\k$-automorphisms of
$X$, which is the group of $\k$-rational points of a group scheme
$\Aut(X)$ that is locally of finite type over $\k$ \cite[Theorem
7.1.1]{Brion-notes} with at most countably many connected
components.

An {\em algebraic group} $G$ over a perfect field $\k$ is a (not
necessarily connected) $\k$-group variety. In particular, $G$ is
reduced and hence smooth \cite[Proposition 2.1.12]{Brion-notes}. We
have $G_{\bk}=G\times_{\mathrm{Spec}(\k)}\mathrm{Spec}(\bk)$, on
which $\Gal(\bk/\k)$ acts through the second factor. The definition
of rational actions of algebraic groups on algebraic varieties goes
back to Weil and Rosenlicht, see \cite{Weil_groups,Rosenlicht}.

\begin{defi}\label{defn:alg-subgrp}
We say that an algebraic group $G$ {\em acts birationally} on a
variety $X$ if
\begin{enumerate}
\item\label{alg-subgrp:1} there are open dense subsets $U,V\subset G\times X$ and a birational map
\[G\times X\dashrightarrow G\times X,\quad (g,x)\dashmapsto(g,\rho(g,x))\]
restricting to a isomorphism $U\rightarrow V$ and the projection of
$U$ and $V$ to the first factor is surjective onto $G$, and
\item\label{alg-subgrp:2}
$\rho(e,\cdot)=\id_X$ and $\rho(gh,x)=\rho(g,\rho(h,x ))$ for any
$g,h\in G$ and $x\in X$ such that $\rho(h,x),\rho(gh,x)$ and
$\rho(g,\rho(h,x))$ are well defined.
\end{enumerate}
\end{defi}

\noindent The group $G(\k)$ of $\k$-points of $G$ is the subgroup of
$G_{\bk}$ of elements fixed by the $\Gal(\bk/\k)$-action, so
we have a map $G(\k)\to\Bir_\k(X)$.
Definition~\ref{defn:alg-subgrp}(\ref{alg-subgrp:2}) implies that it
is a homomorphism of groups, and
Definition~\ref{defn:alg-subgrp}(\ref{alg-subgrp:1}) is equivalent
to the induced map $G(\k)\to\Bir_\k(X)$, $g\to f(g,\cdot)$ being a
so-called morphism, see \cite[Definition 2.1, Definition
2.2]{BlancFurter}, usually denoted by $G\to \Bir_\k(X)$ by abuse of
notation. The notion of morphism from a variety to $\Bir_\k(X)$ goes
back to M. Demazure \cite{MR0284446} and J.-P. Serre
\cite{MR2648675}.

We say that $G$ is an {\em algebraic subgroup} of $\Bir_\k(X)$ if
$G$ acts birationally on $X$ with trivial schematic kernel. We say
that $G$ {\em acts regularly} on $X$ if the birational map in
Definition~\ref{defn:alg-subgrp}(\ref{alg-subgrp:1}) is an
isomorphism. In that case, $G$ is a subgroup of $\Aut(X)$ and we
call $X$ a {\em $G$-surface}.

Let $G$ be an algebraic group acting birationally on surfaces $X_1$
and $X_2$ by birational maps $\rho_i\colon G\times X_i\rat X_i$,
$i=1,2$ as in Definition~\ref{defn:alg-subgrp}. A birational map
$f\colon X_1\rat X_2$ is called {\em $G$-equivariant} if the
following diagram commutes
\[
\begin{tikzpicture}[baseline= (a).base]
\node[scale=1](a) at (0,0){
\begin{tikzcd}
G\times X_1\ar[dashed,r,"\rho_1"]\ar[d,dashed,"\mathrm{id}_G\times f",swap]& X_1\ar[d,dashed,"f"]\\
G\times X_2\ar[dashed,r,"\rho_2"]&X_2
\end{tikzcd}
};
\end{tikzpicture}
\]
In particular, if $\tilde{\rho}_i\colon G\to\Bir_\k(X_i)$ denotes
the induced morphism, the following diagram commutes
\[
\begin{tikzpicture}[baseline= (a).base]
\node[scale=1](a) at (0,0){
\begin{tikzcd}
G(\k)\ar[r,"\tilde{\rho}_1"]\ar[rd,"\tilde{\rho}_2",swap]&\Bir_\k(X_1)\ar[d,"f\circ -\circ f^{-1}"]\\
&\Bir_\k(X_2)
\end{tikzcd}
};
\end{tikzpicture}
\]

The following proposition is proven in \cite[\S2.6]{BlancFurter}
over an algebraically closed field and its proof can be generalised
over any perfect field.

\begin{prop}[{\cite[\S2.6]{BlancFurter}}]\label{rem:linear}
Any algebraic subgroup of $\Bir_\k(\p^2)$ is an affine algebraic
group.
\end{prop}

The following proposition was proven separately by A. Weil and M.
Rosenlicht \cite{Weil_groups, Rosenlicht}, but neither of them
needed the new model to be smooth nor projective. Modern proofs can
also be found in \cite{LonjouUrech} over any field and in
\cite{Cornulier19, Kraft18} over algebraically closed fields.

\begin{prop}\label{thm:projective model}
Let $X$ be a surface and $G$ be an affine algebraic group acting
birationally on $X$. Then there exists a $G$-surface $Y$ and a
$G$-equivariant birational map $X\dashrightarrow Y$. Furthermore,
$G(\k)$ has finite action on $\NS(Y)$.
\end{prop}
\begin{proof}
By \cite{Weil_groups, Rosenlicht}, there exists a normal not necessarily projective or smooth $G$-surface
$Y'$ and a $G$-equivariant birational map $X\rat Y'$. The set $Y''$ of
smooth points of $Y'$ is $G$-stable, it is contained in a complete
surface, which can be desingularised \cite{Lipman}, so $Y''$ is
quasi-projective. By \cite[Corollary 2.14]{Brion2017}, $Y''$ has a
$G$-equivariant completion $Y'''$. We now $G$-equivariantly
desingularise $Y'''$ to obtain the smooth projective surface $Y$
\cite{Zar39,LipmanRS} (the sequence of blow-ups and normalisations
over $\k$ can be done $G$-equivariantly).

The second claim is classical and for instance shown in \cite[Lemma
2.10]{RZ} over any perfect field.
\end{proof}

\subsection{Minimal surfaces}\label{ss:min-surf}

\begin{defi}\label{def:Mfs}
Let $X$ be a surface, $B$ a point or a smooth curve and $\pi\colon
X\to B$ a surjective morphism with connected fibres such that $-K_X$
is $\pi$-ample. We call $\pi\colon X\to B$ a {\em rank $r$
fibration}, where $r=\rk\,\NS(X/B)$.
\begin{itemize}
\item If $B=\pt$ is a point, the surface $X$ is called {\em del Pezzo surface}. Then $X_{\bk}$ is isomorphic to $\p^1_{\bk}\times\p^1_{\bk}$ or to the blow-up of $\p^2_{\bk}$ in at most $8$ points in general position. We call $K_X^2$ the {\em degree of $X$}. Note that $1\leq K_X^2\leq9$.
\item If $B$ is a curve, then $\pi\colon X\to B$ is called {\em conic fibration}; the general geometric fibre of $\pi$ is isomorphic to $\p^1_{\bk}$ and a geometric singular fibre of $\pi$ is the union of two secant $(-1)$-curves over $\bk$. Moreover, if $X$ is rational, then $B=\p^1$, see for instance \cite[Lemma 2.4]{Schneider}.
\item If $r=1$, then $\pi\colon X\to B$ is called {\em Mori fibre space}.
\end{itemize}
We may write $X/B$ instead of $\pi\colon X\to B$. Let $X/B$ and
$X'/B'$ be conic fibrations. We say that a birational map
$\varphi\colon X\rat X'$ {\em preserves the fibration} or is a {\em
birational map of conic fibrations} if the diagram
\[
\begin{tikzpicture}[baseline= (a).base]
\node[scale=1](a) at (0,0){
\begin{tikzcd}
X\ar[d]\ar[r,dashed,"\varphi"]&X'\ar[d]\\
B\ar[r,"\simeq"]&B'
\end{tikzcd}
};
\end{tikzpicture}
\]
commutes.
\end{defi}

For a surface $X$, we can run the $\Gal(\bk/\k)$-equivariant Minimal
Model program on $X_{\bk}$, because the action of $\Gal(\bk/\k)$ on
$\NS(X_{\bk})$ is finite. The end result is a $\Gal(\bk/\k)$-Mori
fibre space $Y_{\bk}/B_{\bk}$ as in Definition~\ref{def:Mfs}, which is equivalent to $Y/B$ being a
Mori fibre space.

\begin{ex}\label{ex:cb}
\leavevmode
\begin{enumerate}
\item\label{ex:Fn}
For $n\geq0$, the Hirzebruch surface $\F_n$ is the quotient of the
action of $(\mathbb{G}_m)^2$ on $(\A^2\setminus\{0\})^2$ by
\[(\mathbb{G}_m)^2\times(\A^2\setminus\{0\})^2\to (\A^2\setminus\{0\})^2,\ (\mu,\rho),(y_0,y_1,z_0,z_1)\mapsto(\mu\rho^{-n}y_0,\mu y_1, \rho z_0,\rho z_1).\]
The class of $(y_0,y_1,z_0,z_1)$ is denoted by $[y_0:y_1;z_0:z_1]$.
The projection $\pi_n\colon\F_n\to\p^1$ given by $[y_0:y_1;z_0:z_1]\mapsto
[z_0:z_1]$ is a conic fibration and the special section
$S_{-n}\subset\F_n$ is given by $y_0=0$.
\item\label{ex:Sl}
Let $p$ and $p'$ be two points of degree $2$ in $\p^2$ with
splitting field $L/\k$ and $L'/\k$, respectively, such that their
geometric components are in general position. We denote by
$\Sl^{L,L'}$ a del Pezzo surface obtained by first blowing up
$p,p'$, and then contracting the line passing through one of the two
points. It has a natural conic fibration structure
$\Sl^{L,L'}\to\p^1$; the fibres are the strict transforms of the
conics in $\p^2$ passing through the two points.
\end{enumerate}
\end{ex}

\begin{lem}{\cite[Lemma 6.11]{Schneider}}\label{rmk:cyclicpoint2}
Let $L/\k$ be a finite extension. Let
$p_1,\dots,p_4,q_1,\dots,q_4\in\p^2(L)$ such that the sets
$\{p_1,\dots,p_4\}$ and $\{q_1,\dots,q_4\}$ are $\Gal(\bk/\k)$-invariant and no three of the $p_i$ and no three of the $q_i$ are
collinear. Suppose that for any $g\in\Gal(\bk/\k)$ there exists
$\sigma\in\sym_4$ such that $p_i^g=p_{\sigma(i)}$ and
$q_i^g=q_{\sigma(i)}$ for $i=1,\dots,4$. Then there exists
$\alpha\in\PGL_3(\k)$ such that $\alpha(p_i)=q_i$ for $i=1,\dots,4$.
\end{lem}

\begin{rmk}\label{rmk:cyclicpoint2-dim1}
The argument of \cite[Lemma 6.11]{Schneider} can be applied to show
the following analogue of Lemma~\ref{rmk:cyclicpoint2} on
$\p^1$: let $F/\k$ be a finite extension and
$p_1,p_2,p_3,q_1,q_2,q_3\in\p^1(F)$ such that the sets
$\{p_1,p_3,p_3\}$ and $\{q_1,q_2,q_3\}$ are $\Gal(\bk/\k)$-invariant. Suppose that for any $g\in\Gal(F/\k)$ there exists
$\sigma\in\sym_3$ such that $p_i^g=p_{\sigma(i)}$ and
$q_i^g=q_{\sigma(i)}$ for $i=1,2,3$. Then there exists
$\alpha\in\PGL_2(\k)$ such that $\alpha(p_i)=q_i$ for $i=1,2,3$.
\end{rmk}

\begin{lem}\label{lem:min-cb}{\cite[Remark 6.1, Lemma 6.13]{Schneider}}
Let $\pi\colon X\to \p^1$ be a Mori fibre space and suppose that $X$
is rational. Then $X$ is isomorphic to a Hirzebruch surface, to a
del Pezzo surface $\Sl^{L,L'}$ or to a del Pezzo surface obtained by
blowing up a point of degree $4$ in $\p^2$.
\end{lem}

\begin{prop}\label{lem:sing-fibres}
Let $X/B$ be a Mori fibre space. If
$B$ is a point, then $X$ is rational if and only if $K_X^2\geq5$ and
$X(\k)\neq\emptyset$.
\end{prop}
\begin{proof}
Suppose that $d:=K_X^2\geq5$ and that $X(\k)$ contains a point $r$.
If $d=7$, then $X_{\bk}$ contains three $(-1)$-curves, one of which
must be $\k$-rational, contradicting $\rk\,\NS(X)=1$. If $d=8$ , the
blow-up of $r$ is a del Pezzo surface of degree $7$, which has two
disjoint $(-1)$-curves over $\bk$ that are either both $\k$-rational or
they make up a $\Gal(\bk/\k)$-orbit of curves. Contracting them
induces a birational map over $\k$ to a del Pezzo surface of degree $9$ with a
rational point, which hence is $\p^2$. This argument also holds if
$\rk\,\NS(X)=2$. Let $d=6$. If $r$ is contained in a curve of negative
self-intersection, then that curve is a $\k$-rational $(-1)$-curve,
contradicting $\rk\,\NS(X)=1$. If $r$ is not contained in any curve of
negative self-intersection, the blow-up of $r$ contains a curve with
three pairwise disjoint geometric components of self-intersection
$-1$. Their contraction yields a birational map $X\rat Y$, where $Y$
is a del Pezzo surface of degree $8$ with a rational point, so $Y$ is rational by the argument above.
If $d=5$,
then again $\rk\,\NS(X)=1$ implies that $r$ is not in a $(-1)$-curve.
After blowing up $r$ we can contract a curve with five pairwise
disjoint geometric components and arrive on a del Pezzo surface of
degree $9$, which is $\p^2$ because it has a rational point.

Let's prove the converse implication. If $X$ is a rational del Pezzo
surface, then $X(\k)\neq\emptyset$ by the Lang-Nishimura theorem.
The remaining claim follows from the classification of {\em Sarkisov
links} (see definition in Section~\ref{ss:sarkisov}) between rational Mori fibre spaces over a perfect field
\cite[Theorem 2.6]{iskovskikh_1996}. Indeed, any birational map
between del Pezzo surfaces over $\k$ with Picard rank $1$ decomposes
into Sarkisov links and automorphisms \cite[Theorem
2.5]{iskovskikh_1996}. The list of Sarkisov links implies the
following: for a del Pezzo surface $X$ with $\rk\,\NS(X)=1$ and
$K_X^2\leq4$, any Sarkisov link $X\rat Y$ that is not an isomorphism
is to a del Pezzo surface $Y$, either of degree $K_Y^2\leq4$ and
$\rk\,\NS(Y)=1$, or of degree $K_Y^2=3$ and $Y$ carries moreover the
structure $Y\to\p^1$ of a Mori fibre space. From the latter, any
Sarkisov link $Y\rat Z$ is to a del Pezzo surface $Z$ of degree
$\leq4$, either with $\rk\,\NS(Z)=1$ or it preserves the fibration and
$\rk\,\NS(Z)=2$. In particular, $X$ cannot be joined to $\p^2$ by a
birational map.
\end{proof}

\begin{lem}\label{lem:aut5}
If $X$ is a del Pezzo surface of degree $K_X^2\leq5$, then
$\Aut_{\bk}(X)$ is finite.
\end{lem}
\begin{proof}
It suffices to show the claim for $\k=\bk$. Then $X$ is the blow-up
of $p_1,\dots,p_r\in\p^2$ in general position with $r=9-K_X^2\geq4$.
It has finitely many $(-1)$-curves, say $n$ of them, and the action
of $\Aut_\k(X)$ on the set of the $(-1)$-curves induces an exact
sequence
\[1\rightarrow\Aut_\k(\p^2,p_1,\dots,p_r)\to\Aut_\k(X)\to\sym_{n}.\]
Since $p_1,\dots,p_r$ are in general position and $r\geq4$, the
group $\Aut_\k(\p^2,p_1,\dots,p_r)$ is trivial, which yields the
claim.
\end{proof}

\subsection{Relatively minimal surfaces}

We now generalise the  notion of being a minimal surface to  being
minimal relative to the action of an affine algebraic group.

\begin{defi}\label{def:GMfs}
Let $G$ be an affine algebraic group, let $X$ be a $G$-surface and $\pi\colon
X\to B$ a rank $r$ fibration.
\begin{enumerate}
\item If $\pi$ is $G$-equivariant and $r':=\rk\,\NS(X_{\bk}/B_{\bk})^{G_{\bk}\times\Gal(\bk/\k)}$, we call $\pi$ a {\em $G$-equivariant rank $r'$ fibration}.
If $r'=1$ we call it a {\em $G$-Mori fibre space}.

\item If $\pi$ is $G(\k)$-equivariant and $r'':=\rk\,\NS(X/B)^{G(\k)}$, we call $\pi$ a {\em $G$-equivariant rank $r''$ fibration}.
If $r''=1$ we call it {\em $G(\k)$-Mori fibre space}.
\end{enumerate}
\end{defi}

If a rank $r$ fibration $X\to B$ is $G$-equivariant, we have $r\geq
r''\geq r'$. A $G$-Mori fibre space is not necessarily a
$G(\k)$-Mori fibre space, since $G(\k)$-equivariant does not imply
$G$-equivariant. Examples are, for instance, the del Pezzo surfaces
in Lemma~\ref{prop:DP2} and Lemma~\ref{prop:DP3} (see also
Theorem~\ref{thm:1}(\ref{1:5a})), that are $\Aut(X)$-Mori fibre
spaces but not $\Aut_\k(X)$-Mori fibre spaces.

If $G$ is connected, Blanchard's Lemma \cite[Theorem 7.2.1]{Brion2017}
implies that a $G$-Mori fibre space is a Mori fibre space. However,
the affine algebraic groups we are going to work with are not necessarily
connected. All del Pezzo surfaces $X$ of degree $6$ in
\S\ref{sec:DP6} are $\Aut(X)$-Mori fibre spaces, all but two of them
are also $\Aut_\k(X)$-Mori fibre spaces and only two of them are
Mori fibre spaces. \medskip

Let $G$ be an affine algebraic group and $X$ a $G$-surface. The action
$\rho\colon G\times X\to X$ from Definition~\ref{defn:alg-subgrp}
being defined over $\k$ is equivalent to
$\bar{\rho}:=\rho\times\id\colon G_{\bk}\times X_{\bk}\to X_{\bk}$
being $\Gal(\bk/\k)$-equivariant, {\it i.e.}
$\bar{\rho}(g,x)^h=\bar{\rho}(g^h,x^h)$ for any $h\in\Gal(\bk/\k)$,
$g\in G_{\bk}$, $x\in X_{\bk}$. We can therefore see the $G$-action
on $X$ as the $(\Gal(\bk/\k)\times G_{\bk})$-action on $X_{\bk}$
\[
(\Gal(\bk/\k)\times G_{\bk})\times X_{\bk}\to X_{\bk},\quad
(h,g,x)\mapsto \bar{\rho}(g^h,x^h)
\]
satisfying $\bar{\rho}(g^h,x^h)=\bar{\rho}(g,x)^h$ for any
$h\in\Gal(\bk/\k)$, $g\in G_{\bk}$, $x\in X_{\bk}$.

\begin{rmk}\label{rmk:mmp}
Let $G$ be an affine algebraic group and $X$ a $G$-surface such that
$X_{\bk}$ is rational. By Proposition~\ref{thm:projective model},
the group $G_{\bk}$ and hence also the group $\Gal(\bk/\k)\times
G_{\bk}$ has finite action on $\NS(X_{\bk})$. We can run the
$(\Gal(\bk/\k)\times G_{\bk})$-equivariant Minimal Model program on
$X_{\bk}$, and by \cite[Example 2.18]{KollarMori} the end result is
a $G$-Mori fibre space $Y/B$. We then restrict to the $G(\k)$-action
on $Y$ and recall that $G(\k)$ has finite action on $\NS(Y)$ by
Proposition~\ref{thm:projective model}. Since $Y/B$ is
$G$-equivariant, it is also $G(\k)$-equivariant, and we can run the
$G(\k)$-equivariant Minimal Model Program on $Y$, whose end result
is then a $G(\k)$-Mori fibre space.
\end{rmk}

Let us tidy up the direction for classifying the infinite algebraic
subgroups of $\Bir_\k(\p^2)$.

\begin{prop}\label{prop which cases}
Let $G$ be an infinite algebraic subgroup of $\Bir_\k(\p^2)$. Then
there exists a $G$-equivariant birational map $\p^2\dashrightarrow
X$ to a $G$-Mori fibre space $X/B$ that is one of the following:
\begin{enumerate}
\item $B$ is a point and $X\simeq\p^2$ or $X$ is a del Pezzo surface of degree $6$ or $8$.
\item $B=\p^1$ and there exists a birational morphism of conic fibrations $X\to\Sl^{L,L'}$ or $X\to\F_n$ for some $n\geq0$.
\end{enumerate}
\end{prop}
\begin{proof}
By Proposition~\ref{rem:linear}, $G$ is an affine algebraic group.
By Proposition~\ref{thm:projective model}, there is a $G$-surface
$X'$ and a $G$-equivariant birational map $\phi\colon \p^2\dasharrow
X'$. We now apply the $(G_{\bk}\times\Gal(\bk/\k))$-equivariant
Minimal Model Program and obtain a $G$-equivariant birational
morphism $X'\to X$ to a $G$-Mori fibre space $\pi \colon X\to B$,
see Remark~\ref{rmk:mmp}.

If $B$ is a point, then $X$ is a del Pezzo surface.  Since $G$ is
infinite, Lemma~\ref{lem:aut5} implies that $K_X^2\geq6$. If
$K^2=7$, then $X_{\bk}$ contains exactly three $(-1)$-curves, one of
which is $G_{\bk}\times\Gal(\bk/\k)$-invariant, so $X$ is not a
$G$-Mori fibre space. It follows that $K_X^2\in\{6,8,9\}$, and if
$K_X^2=9$, then $X\simeq\p^2$ by Ch\^{a}telet's Theorem.

Suppose that $B=\p^1$. Then there is a birational morphism $X\to Y$
of conic fibrations onto a Mori fibre space $Y/\p^1$. By
Lemma~\ref{lem:min-cb}, $Y$ is a Hirzebruch surface, $Y\simeq\Sl$ or
$Y$ is the blow-up of $\p^2$ in a point of degree $4$ whose
geometric components are in general position. The latter is a del
Pezzo surface of degree $5$, so by Lemma~\ref{lem:aut5} the group
$\Aut_{\bk}(Y)$ is finite, which does not occur under our
hypothesis. It follows that $Y\simeq\F_n$, $n\geq0$, or
$Y\simeq\Sl^{L,L'}$.
\end{proof}

\begin{lem}\label{lem:aut}
\leavevmode
\begin{enumerate}
\item If $X$ is a del Pezzo surface, then $\Aut(X)$ is an affine algebraic group.
\item Let $\pi\colon X\rightarrow\p^1$ be a conic fibration such that $X_{\bk}$ is rational. Then $\Aut(X,\pi)$ is an affine algebraic group.
\end{enumerate}
\end{lem}
\begin{proof}
(1) Let $N:=h^0(-K_X)$. Then $\Aut(X)$ preserves the ample divisor
$-K_X$, thus it is conjugate via the embedding $|-K_X|\colon
X\hookrightarrow\p^{N-1}$ to a closed subgroup of
$\Aut(\p^{N-1})\simeq\PGL_{N}$ and is hence affine.

(2) Let $G$ be the schematic kernel of $\Aut(X,\pi)\to\Aut(\NS(X))$.
If $D$ is an ample divisor on $X$, it is fixed by $G$ and hence (as
above) $G$ is an affine algebraic group. Since $X_{\bk}$ is rational
and has the structure of a conic fibration, we have
$\NS(X)\simeq\Z^n$ for some $n\geq2$, and it is generated by $-K_X$,
the general fibre and components of the singular fibres. The
(abstract) group $H:=\Aut(X,\pi)/G$ acts faithfully on $\NS(X)$,
fixes $-K_X$ and the general fibre and permutes the components of
the singular fibres. It follows that $H$ is isomorphic (as abstract
group) to a subgroup of $\GL_n(\Z)$ whose elements have entries in $\{0,\pm1\}$.
Therefore, $H$
is finite and hence $\Aut(X,\pi)$ is an affine algebraic group.
\end{proof}

In particular, if $X$ is a del Pezzo surface, the
$\Gal(\bk/\k)$-action on $\Aut_{\bk}(X)$ is a $\k$-structure with fixed
locus $\Aut_\k(X)$. Similarly, if $\pi\colon X\rightarrow\p^1$ is a
conic fibration such that $X_{\bk}$ is rational, then the
$\Gal(\bk/\k)$-action on $\Aut_{\bk}(X,\pi)$ is a $\k$-structure with
fixed locus $\Aut_\k(X,\pi)$.

Our goal is to classify algebraic subgroups of $\Bir_\k(\p^2)$ up to
conjugacy and inclusion. Proposition~\ref{prop which cases} and
Lemma~\ref{lem:aut} imply that it suffices to classify up to
conjugacy and inclusion the automorphism groups of del Pezzo
surfaces of degree $6$ and $8$ and the automorphism groups of
certain conic fibrations.

\section{Del Pezzo surfaces of degree 8}\label{sec:DP89}
We now classify the rational del Pezzo surfaces of degree $8$. Over
an algebraically closed field, any such surface is isomorphic to the
blow-up of $\p^2$ in a point or to $\p^1\times\p^1$. Over $\R$,
there are exactly two rational models of the latter, namely the
quadric surfaces given by $w^2+x^2-y^2-z^2=0$ or $w^2+x^2+y^2-z^2=0$
in $\p^3$. The first is isomorphic to $\p^1_{\R}\times\p^1_{\R}$ and
the second is the $\R$-form of $\p^1_{\C}\times\p^1_{\C}$ given by
$(x,y)\mapsto(y^g,x^g)$, where $\langle g\rangle=\Gal(\C/\R)$. We
now show that the classification is similar over an arbitrary
perfect field $\k$.

\begin{defi}\label{def:Q}
Suppose that $\k$ has a quadratic extension $L/\k$. We denote by
$\Ql^L$ the $\k$-structure on $\p^1_L\times\p^1_L$ given by
$([u_0:u_1],[v_0:v_1])\mapsto([v_0^g:v_1^g],[u_0^g:u_1^g])$, where
$g$ is the generator of $\Gal(L/\k)$.
\end{defi}

The surface $\Ql^L$ is a del Pezzo surface of degree $8$ and it is
rational by Proposition~\ref{lem:sing-fibres} because the point
$([1:1],[1:1])$ $\in \Ql^L(\k)$.

\begin{lem}\label{prop:Q}
Let $X$ be a rational del Pezzo surface of degree $8$.
\begin{enumerate}
\item\label{Q:1} We have $\rk\,\NS(X)=2$ if and only if $X\simeq\F_0$ or $X\simeq\F_1$, and $\rk\,\NS(X)=1$ if and only if $X\simeq \Ql^L$ for some quadratic extension $L/\k$.
\item\label{Q:2} $X\simeq \Ql^L$ if and only if for any $p\in X(\k)$ there is a birational map $X\rat \p^2$ that is the composition of the blow-up of $p$ and the contraction of a curve onto a point of degree $2$ in $\p^2$ whose splitting field is $L$.
\item\label{Q:3} We have $\Ql^L\simeq \Ql^{L'}$ if and only if $L$ and $L'$ are $\k$-isomorphic.
\end{enumerate}
\end{lem}
\begin{proof}
(\ref{Q:1}--\ref{Q:2}) The surface $X_{\bk}$ is a del Pezzo
surface of degree $8$ over $\bk$ and is hence isomorphic to
$\p^1_{\bk}\times\p^1_{\bk}$ or to $(\F_1)_{\bk}$. In the latter
case, the unique $(-1)$-curve is $\Gal(\bk/\k)$-invariant, hence
$X\simeq\F_1$. Suppose that $X_{\bk}$ is isomorphic to
$\p^1_{\bk}\times\p^1_{\bk}$ and consider the blow-up $\pi_1\colon
Y\to X$ of $X$ in a rational point $p\in X(\k)$ (such a point exists
by Proposition~\ref{lem:sing-fibres}). Then $Y$ is a del Pezzo
surface of degree $7$ and $Y_{\bk}$ has three $(-1)$-curves, one of
which is the exceptional divisor over the rational point $p$. The
union of the other two $(-1)$-curves $C_1,C_2\subset Y_{\bk}$ is
preserved by $\Gal(\bk/\k)$, and hence their contraction yields a
birational morphism $\pi_2\colon Y\to\p^2$. If each of $C_1$ and
$C_2$ is preserved by $\Gal(\bk/\k)$, then
$\varphi:=\pi_1\pi_2^{-1}\colon \p^2\rat X$ has two rational
base-points. The pencil of lines through each base-point is sent
onto a fibration of $X$, and Lemma~\ref{lem:min-cb} implies that $X$
is a Hirzebruch surface, so $X\simeq\F_0$. If $C_1\cup C_2$ is a
$\Gal(\bk/\k)$-orbit of curves, then $\varphi$ has a base-point  $q$
of degree $2$. By Remark~\ref{rmk:cyclicpoint2} we can assume that
$q$ is of the form $q=\{[a_1:1:0],[a_2:1:0]\}$, $a_1,a_2\in\bk$. We
consider the projection $\psi\colon \p^2_{\bk}\rat
\p^1_{\bk}\times\p^1_{\bk}$ away from $q$
\begin{align*}
\psi\colon [x:y:z]&\dashmapsto([x-a_1y:z],[x-a_2y:z])\\
\psi^{-1}\colon
([u_0:u_1],[v_0:v_1])&\dashmapsto[-a_2u_0v_1+a_1v_0u_1:
-u_0v_1+v_0u_1: (a_1-a_2)u_1v_1]
\end{align*}
whose inverse $\psi^{-1}$ has base-point $([1:1],[1:1])$. There
exists an isomorphism $\alpha\colon X_{\bk}\stackrel{\simeq}\to
\p^1_{\bk}\times\p^1_{\bk}$ such that  $\alpha\varphi=\psi$. Let
$\rho$ be the canonical action of $\Gal(\bk/\k)$ on $\p^2_{\bk}$.
Then the action $\varphi\rho\varphi^{-1}$ on $X_{\bk}$ corresponds
to the $\k$-structure $X$. It follows that the action of
$\psi\rho\psi^{-1}=\alpha(\varphi\rho\varphi^{-1})\alpha^{-1}$  on
$\p^1_{\bk}\times\p^1_{\bk}$ corresponds to  a $\k$-structure isomorphic
to $X$. For any $g\in \Gal(\bk/\k)$, we have
\[
\psi \rho_g\psi^{-1}\colon ([u_0:u_1],[v_0:v_1])\mapsto
\begin{cases}
([v_0^g:v_1^g],[u_0^g:u_1^g]),&\text{if $a_1^g=a_2$}\\
([u_0^g:u_1^g],[v_0^g:v_1^g]),&\text{if $a_1^g=a_1$.}
\end{cases}
\]
If $L=\k(a_1,a_2)$, which is a quadratic extension of $\k$, then the
generator $g$ of $\Gal(L/\k)$ exchanges the geometric components of
$q$, so $X\simeq \Ql^L$.

(\ref{Q:3}) The surfaces $\Ql^L$ and $\Ql^{L'}$ are isomorphic if
and only if there exist birational maps $\varphi\colon
\Ql^L\rat\p^2$ and $\varphi'\colon \Ql^{L'}\rat\p^2$ as in
(\ref{Q:2}) and $\alpha\in\Aut_\k(\p^2)$ such that
$\varphi^{-1}\alpha\varphi'$ is an isomorphism. This is the case if
and only if the base-points of $\varphi^{-1}$ and $(\varphi')^{-1}$
have the same splitting field. This is equivalent to $L$ and $L'$
being $\k$-isomorphic.
\end{proof}

In order to be complete, we now show an isomorphism from $\Ql^L$ to
a quadratic surface $\Rl^L$ in $\p^3$. Later on, we will choose to
use or announce claims using coordinates in $\Ql^L$ or in $\Rl^L$
according to practicality.

\begin{lem}\label{lem:Q-deg22}
Let $L=\k(a_1)$ be a quadratic extensions of $\k$ and let
$t^2+at+\tilde{a}=(t-a_1)(t-a_2)\in\k[t]$ be the minimal polynomial
of $a_1$. The following hold:
\begin{enumerate}
\item\label{Q-deg22:1} Let $\Rl^L\subset\p^3_{WXYZ}$ be the quadric surface given by $WZ=X^2+aXY+\tilde{a}Y^2$. Then
\[\p^2\rat\Rl^L,\quad [x:y:z]\dashmapsto [x^2+axy+\tilde{a}y^2:xz:yz:z^2]\]
is birational, and $\Rl^L$ is isomorphic to $\Ql^L$.
\item\label{Q-deg22:2} The map $\Ql^L\to \Rl^L$ given by
\begin{align*}
 ([u_0:u_1],[v_0:v_1])&\mapsto[u_0v_0(a_1-a_2):-a_2u_0v_1+a_1u_1v_0: -u_0v_1+u_1v_0:(a_1-a_2)u_1v_1]\\
 [W:X:Y:Z]&\mapsto([X-a_1Y:Z],[X-a_2Y:Z])=([W:X-a_2Y],[W:X-a_1Y])
\end{align*}
is an isomorphism over $\k$.
\item Let $p\in\Ql^L$ be a point of degree $2$ with splitting field $L'=\k(b_1)$ whose components are not on the same ruling of $\Ql_L^L$. Let $t^2+bt+\tilde{b}=(t-b_1)(t-b_2)\in\k[t]$ be the minimal polynomial of $b_1$ over $\k$.
    \begin{enumerate}
    \item\label{Q-deg22:3} Then there is an automorphism of $\Ql^L$ (resp. $\Rl^L$) that sends $p$ respectively onto
    \[\{([b_1:1],[b_1:1]),([b_2:1],[b_2:1])\},\quad \{[b_1^2:b_1:0:1],[b_2^2:b_2:0:1]\}\]
    \item\label{Q-deg22:4} The pencil of $(1,1)$-curves in $\Ql^L$ through $p$ is given in $X^L$ by the pencil of hyperplanes whose equations are $\lambda(W+bX+\tilde{b}Z)+\mu Y=0$ for $[\lambda:\mu]\in\p^1$.
    \end{enumerate}
\end{enumerate}
\end{lem}

\begin{proof}
(\ref{Q-deg22:1}) The given birational map has a single base-point
of degree $2$, namely $q=\{[a_1:1:0],[a_2:1:0]\}$, and it contracts
the line $z=0$. Its image is the quadric surface $\Rl^L$ given by
$WZ=X^2+aXY+\tilde{a}Y^2$, and the inverse map $\Rl^L\rat\p^2$ is
given by the projection from $[1:0:0:0]$. So $\Rl^L\simeq\Ql^L$ by
Lemma~\ref{prop:Q}(\ref{Q:2}).

(\ref{Q-deg22:2}) We compose the birational map from
(\ref{Q-deg22:1}) and the birational map $\psi\colon\p^2\rat\Ql^L$
from the proof of Lemma~\ref{prop:Q}(\ref{Q:2}) whose base-point is
$\{[a_1:1:0],[a_2:1:0]\}$.

(\ref{Q-deg22:3}) We see from the description of $\Aut_\k(\Ql^L)$ in
Lemma~\ref{lem:autQ} that we can assume that $p$ is not in the
ruling of $\Ql^L_L$ passing through $([1:1],[1:1])$. The birational
map $\psi\colon\Ql^L\rat\p^2$ from the proof of
Lemma~\ref{prop:Q}(\ref{Q:1}) sends $p$ onto a point $\psi(p)$ in
$\p^2$ that is not collinear with $\{[a_1:1:0],[a_2:1:0]\}$. By
Lemma~\ref{rmk:cyclicpoint2}, there exists an element
$\alpha\in\Aut_\k(\p^2)$ that sends $\psi(p)$ onto
$\{[b_1:0:1],[b_2:0:1]\}$. Then
$\psi^{-1}\alpha\psi\in\Aut_\k(\Ql^L)$ and sends $p$ onto
$\{([b_1:1],[b_1:1]),([b_2:1],[b_2:1])\}$. We use the isomorphism
from (\ref{Q-deg22:2}) to compute its coordinates in $\Rl^L$.

(\ref{Q-deg22:4}) The pencil of $(1,1)$-curves through $p$ is sent
by $\psi\colon\Ql^L\rat\p^2$ onto the pencil of conics through
through $[a_1:1:0],[a_2:1:0],[b_1:0:1],[b_2:0:1]$. It is given by
$\lambda(x^2+axy+bxz+\tilde{a}y^2+\tilde{b}z^2)+\mu yz$, and
corresponds via $\psi$ to the pencil in the claim.
\end{proof}

\medskip
\begin{rmk}\label{rmk:minpoly_quadext}
    Let $L=\k(a_1)$ be a quadratic extension of $\k$ and let $t^2+at+\tilde{a}=(t-a_1)(t-a_2)\in\k[t]$
    be the minimal polynomial of $a_1$. Depending on the characteristic of $\k$, we can assume the values of $a$ to be $0$ or $1$:
    \begin{itemize}
        \item If the characteristic of $\k$ is not $2$, then we can assume that $a=0$, namely via the $\k$-isomorphism $t\mapsto t-a/2$.
        \item If the characteristic of $\k$ equals $2$, then we can assume that $a=1$. Indeed, as we assume that $\k$ is a perfect field, all elements of $\k$ are squares, and so $a=0$ does not give an irreducible polynomial over $\k$. The $\k$-isomorphism $t\mapsto t/a$ reduces $a\neq0$ to $a=1$.
    \end{itemize}
\end{rmk}

\medskip
\begin{lem}\label{lem:autQ}
Let $L/\k$ be an extension of degree $2$ and let $g$ be the
generator of $\Gal(L/\k)$. The group $\Aut(\Ql^L)\simeq\Aut(\Rl^L)$
is isomorphic to the $\k$-structure on
$\Aut(\p^1_L\times\p^1_L)\simeq\Aut(\p^1_L)^2\rtimes\langle(u,v)\stackrel{\tau}\mapsto(v,u)\rangle$
given by the $\Gal(L/\k)$-action
\[(A,B,\tau)^g=(B^g,A^g,\tau),\]
where $A\mapsto A^g$ is the canonical $\Gal(L/\k)$-action on
$\Aut(\p^1_L)$. Furthermore,
\[\Aut_\k(\Rl^L)\simeq\Aut_\k(Q^L)\simeq\{(A,A^g)\mid A\in\PGL_2(L)\}\rtimes\langle\tau\rangle.\]
\end{lem}
\begin{proof}
Since $\Ql^L$ is the $\k$-structure on
$Q_{L}^L\simeq\p^1_{L}\times\p^1_{L}$,  the $\Gal(L/\k)$-action on
the algebraic group
\[\Aut_L(\Ql^L)=\Aut(\p^1_L\times\p^1_L)\simeq
\Aut(\p^1_L)^2\rtimes\langle\tau\rangle\]
is a $\k$-structure with fixed
points $\Aut_\k(\Ql^L)$. The automorphism $\tau$ commutes with $g$,
and we have
\[
(A,B)^g(q^g,p^g)=(A,B)^g(p,q)^g=\left((A,B)(p,q)\right)^g\\
=(Ap,Bq)^g=(B^gq^g,A^gp^g)
\]
for any $(A,B)\in\Aut(\p^1_L)^2$ and any $(p,q)\in \Ql^L$. It follows that $(A,B)^g=(B^g,A^g)$. The
group $\Aut_\k(\Ql^L)$ is isomorphic to the subgroup of elements of
$\Aut(\p^1_{L}\times\p^1_{L})$ commuting with $\Gal(L/\k)$, which
yields the remaining claim.
\end{proof}

\medskip
By the following lemma, whenever we contract a curve onto a point of
degree $2$ in $\Ql^L$ with splitting field $L$, we can choose the
point conveniently.

\begin{lem}\label{lem:Q-deg2}
\leavevmode
\begin{enumerate}
\item\label{Q-deg2:1} Let $p\in \Ql^L$ be a point of degree $2$ whose geometric components are not on the same ruling of $\Ql^L_{\bk}\simeq\p^1_{\bk}\times\p^1_{\bk}$ and whose splitting field is $L$. Then there exists $\alpha\in\Aut_\k(\Ql^L)$ such that
$\alpha(p)=\{([1:0],[0:1]),$ $([0:1],[1:0])\}$.
\item\label{Q-deg2:2} Let $r,s\in \Ql^L(\k)$ be two rational points not contained in the same ruling of $\Ql_{\bk}^L$. Then there exists $\alpha\in\Aut_\k(\Ql^L)$ such that
$\alpha(r)=([1:0],[1:0])$ and $\alpha(s)=([0:1],[0:1]).$
\end{enumerate}
\end{lem}
\begin{proof}
Let $g$ be the generator of $\Gal(L/\k)$.

(\ref{Q-deg2:2}) We have $r=([a:b],[a^g:b^g])$ and
$s=([c:d],[c^g:d^g])$ for some $a,b,c,d\in L$, and $ad-cd\neq0$
because $r$ and $s$ are not on the same ruling of $Q_L$. It follows
that the map $A\colon[u:v]\mapsto[du-cv:-bu+av]$ is contained in
$\PGL_2(L)$. Then $(A,A^g)\in\Aut_\k(Q)$ and it sends respectively $r$ and $s$
onto $([1:0],[1:0])$ and $([0:1],[0:1])$.

(\ref{Q-deg2:1}) The point $p$ is of the form
$\{([a:b],[c:d]),([c^g:d^g],[a^g:b^g])\}$ for some $a,b,c,d\in L$,
and $ad^g-bc^g\neq0$ because its components are not on the same
ruling of $\Ql_L^L$. It follows that the map $A$ defined by
$[u:v]\mapsto[d^gu-c^gv:-bu+av]$ is contained in $\PGL_2(L)$.
Then $(A,A^g)\in\Aut_\k(\Ql^L)$ and it sends $p$ onto
$\{([1:0],[0:1]),([0:1],[1:0])\}$.
\end{proof}

\begin{lem}\label{lem:cyclicpoint3}
Let $p=\{p_1,p_2,p_3\}$ and $q=\{q_1,q_2,q_3\}$ be points in $\Ql^L$
of degree $3$ such that for any $h\in\Gal(\bk/\k)$ there exists
$\sigma\in\sym_3$ such that $p_i^h=p_{\sigma(i)}$ and
$q_i^h=q_{\sigma(i)}$. Suppose that the geometric components of $p$
$($resp. of $q)$ are in pairwise distinct rulings of $\Ql^L_L$. Then
there exists $\alpha\in\Aut_\k(\Ql^L)$ such that $\alpha(p_i)=q_i$
for $i=1,2,3$.
\end{lem}
\begin{proof}
Let $g$ be the generator of $\Gal(L/\k)$. Since $p$ and $q$ are of
degree $3$, we have $p_i^g=p_i$ and $q_i^g=q_i$ for $i=1,2,3$, and
therefore $p_i=(a_i,a_i^g)$ and $q_i=(b_i,b_i^g)$, $a_i,b_i\in\bk$,
for $i=1,2,3$. By hypothesis, for any $h\in\Gal(\bk/L)$ there exists
$\sigma\in\sym_3$ such that
$(a_i^h,a_i^{gh})=p_i^h=q_{\sigma(i)}=(b_{\sigma(i)},b_{\sigma(i)}^g)$.
We apply Remark~\ref{rmk:cyclicpoint2-dim1} to the
$\Gal(\bar{L}/L)$-invariant sets $\{a_1,a_2,a_3\}$ and
$\{b_1,b_2,b_3\}$ in $\p^1_L$ and to the $\Gal(\bar{L}/L)$-invariant
sets $\{a_1^g,a_2^g,a_3^g\}$ and $\{b_1^g,b_2^g,b_3^g\}$ in
$\p^1_L$. There exist $A,B\in\PGL_2(L)$ such that $Aa_i=b_i$ and
$Ba_i^g=b_i^g$ for $i=1,2,3$. Then $A^ga_i^g=(Aa_i)^g= b_i^g=Ba_i^g$
for $i=1,2,3$, and therefore $B=A^g$. It follows that
$\alpha\in\Aut_\k(\Ql^L)$.
\end{proof}

\section{Del Pezzo surfaces of degree 6}\label{sec:DP6}
In this section, we classify the rational del Pezzo surfaces of
degree $6$ over a perfect field $\k$ and describe their automorphism
groups.

\subsection{Options for rational del Pezzo surfaces of degree $6$}
Let $X$ be a rational del Pezzo surface of degree $6$. Then
$X_{\bk}$ is the blow up of three points in $\p^2_{\bk}$, its
$(-1)$-curves are the three exceptional divisors and strict
transforms of the lines passing through two of the three points, and
they form a hexagon. The hexagon of $X_{\bk}$ is $\Gal(\bk/\k)$-invariant. The Galois group $\Gal(\bk/\k)$ acts on the hexagon by
symmetries, so we have a homomorphism of groups
\[\Gal(\bk/\k)\stackrel{\rho}\to \sym_3\times\Z/2\subseteq \Aut(\NS(X_{\bk})).\]
By {\em hexagon of $X$} we mean the hexagon of $X_{\bk}$ endowed with it canonical $\Gal(\bk/\k)$-action.
The options for the non-trivial action of $\rho(\Gal(\bk/\k))$ on
the hexagon of $X$ are visualised in Figure~\ref{fig:action-dP6}.

\begin{figure}[ht]
\begin{minipage}[ht]{.18\textwidth}
\begin{equation}\label{dp6:1}
\begin{tikzpicture}[scale=.6]
\begin{scope}[every coordinate/.style={shift={(0,3.5)}}]  
\path [c](0,0) pic {hexagon};
\draw[<->, thick] [c](-0.5,1) to [bend right=40] (0.5,1); 
\draw[<->, thick] [c](1,0.8) to [bend right=40] (1.4,0.2); 
\draw[<->, thick] [c](1,-0.8) to [bend left=40] (1.4,-0.2); 
\draw[<->, thick] [c](-0.5,-1) to [bend left=40] (0.5,-1); 
\draw[<->, thick] [c](-1,-0.8) to [bend right=40] (-1.4,-0.2); 
\draw[<->, thick] [c](-1,0.8) to [bend left=40] (-1.4,0.2); 
\end{scope}
\end{tikzpicture}
\end{equation}
\end{minipage}
\begin{minipage}[ht]{.18\textwidth}
\begin{equation}\label{dp6:2}
\begin{tikzpicture}[scale=.6,font=\footnotesize]
\begin{scope}[every coordinate/.style={shift={(0,3.5)}}] 
\path [c](0,0) pic {hexagon}; \draw[<->, thick] (D1) to (D4);
\draw[<->, thick] (D2) to [bend right=15,swap] (D3); \draw[<->,
thick] (D6) to [bend left=15,swap]  (D5);
\end{scope}
\end{tikzpicture}
\end{equation}
\end{minipage}
\begin{minipage}[ht]{.18\textwidth}
\begin{equation}\label{dp6:3}
\begin{tikzpicture}[scale=.6,font=\footnotesize]
\begin{scope}[every coordinate/.style={shift={(0,3.5)}}] 
\path [c](0,0) pic {hexagon}; \draw[<->, thick] (D1) to [bend
right=15,swap] (D3); \draw[<->, thick] (D6) to [bend left=15,swap]
(D4);
\draw[<->, thick] [c](1,0.8) to [bend right=40] (1.4,0.2); 
\draw[<->, thick] [c](-1,-0.8) to [bend right=40] (-1.4,-0.2); 
\end{scope}
\end{tikzpicture}
\end{equation}
\end{minipage}
\begin{minipage}[ht]{.18\textwidth}
\begin{equation}\label{dp6:4}
\begin{tikzpicture}[scale=.6,font=\footnotesize]
\begin{scope}[every coordinate/.style={shift={(0,0)}}] 
\path [c](0,0) pic {hexagon}; \draw[<->, thick] (D1) to (D4);
\draw[<->, thick] (D2) to (D5); \draw[<->, thick] (D3) to (D6);
\end{scope}
\end{tikzpicture}
\end{equation}
\end{minipage}
\begin{minipage}[ht]{.18\textwidth}
\begin{equation}\label{dp6:5}
\begin{tikzpicture}[scale=.6,font=\footnotesize]
\begin{scope}[every coordinate/.style={shift={(0,3.5)}}] 
\path [c](0,0) pic {hexagon}; \draw[<->, thick] (D66) to [bend
right=15,swap] (D2); \draw[<->, thick] (D22) to [bend right=15,swap]
(D3); \draw[<->, thick] (D33) to [bend right=15,swap] (D5);
\draw[<->, thick] (D55) to [bend right=15,swap] (D6); \draw[<->,
thick, shorten <=.1cm, shorten >=.1cm] (D1) to (D4); \draw[<->,
thick, shorten <=.2cm, shorten >=.2cm] (D2) to (D5); \draw[<->,
thick, shorten <=.2cm, shorten >=.2cm] (D3) to (D6);
\end{scope}
\end{tikzpicture}
\end{equation}
\end{minipage}

\begin{minipage}[ht]{.2\textwidth}
\begin{equation}\label{dp6:6}
\begin{tikzpicture}[scale=.6,font=\footnotesize]
\begin{scope}[every coordinate/.style={shift={(0,3.5)}}] 
\path [c](0,0) pic {hexagon}; \draw[->, thick] (D11) to [bend
right=15,swap] (D3); \draw[->, thick] (D22) to [bend right=15,swap]
(D4); \draw[->, thick] (D33) to [bend right=15,swap] (D5); \draw[->,
thick] (D44) to [bend right=15,swap] (D6); \draw[->, thick] (D55) to
[bend right=15,swap] (D1); \draw[->, thick] (D66) to [bend
right=15,swap] (D2);
\end{scope}
\end{tikzpicture}
\end{equation}
\end{minipage}
\begin{minipage}[ht]{.23\textwidth}
\begin{equation}\label{dp6:7}
\begin{tikzpicture}[scale=.6,font=\footnotesize]
\begin{scope}[every coordinate/.style={shift={(0,3.5)}}] 
\path [c](0,0) pic {hexagon}; \draw[->, thick] (D11) to [bend
right=15,swap] (D2); \draw[->, thick] (D22) to [bend right=15,swap]
(D3); \draw[->, thick] (D33) to [bend right=15,swap] (D4); \draw[->,
thick] (D44) to [bend right=15,swap] (D5); \draw[->, thick] (D55) to
[bend right=15,swap] (D6); \draw[->, thick] (D66) to [bend
right=15,swap] (D1);
\end{scope}
\end{tikzpicture}
\end{equation}
\end{minipage}
\begin{minipage}[ht]{.23\textwidth}
\begin{equation}\label{dp6:8}
\begin{tikzpicture}[scale=6,font=\footnotesize]
\begin{scope}[every coordinate/.style={shift={(0,-3.5)}}] 
\path [c](0,0) pic {hexagon}; \draw[<->, thick] (D11) to [bend
right=15,swap] (D3); \draw[<->, thick] (D22) to [bend right=15,swap]
(D4); \draw[<->, thick] (D33) to [bend right=15,swap] (D5);
\draw[<->, thick] (D44) to [bend right=15,swap] (D6); \draw[<->,
thick] (D55) to [bend right=15,swap] (D1); \draw[<->, thick] (D66)
to [bend right=15,swap] (D2);
\end{scope}
\end{tikzpicture}
\end{equation}
\end{minipage}
\begin{minipage}[ht]{.23\textwidth}
\begin{equation}\label{dp6:9}
\begin{tikzpicture}[scale=.6,font=\footnotesize]
\begin{scope}[every coordinate/.style={shift={(0,3.5)}}] 
\path [c](0,0) pic {hexagon}; \draw[<->, thick] (D11) to [bend
right=15,swap] (D3); \draw[<->, thick] (D22) to [bend right=15,swap]
(D4); \draw[<->, thick] (D33) to [bend right=15,swap] (D5);
\draw[<->, thick] (D44) to [bend right=15,swap] (D6); \draw[<->,
thick] (D55) to [bend right=15,swap] (D1); \draw[<->, thick] (D66)
to [bend right=15,swap] (D2); \draw[<->, thick, shorten <=.1cm,
shorten >=.1cm] (D1) to (D4); \draw[<->, thick, shorten <=.2cm,
shorten >=.2cm] (D2) to (D5); \draw[<->, thick, shorten <=.2cm,
shorten >=.2cm] (D3) to (D6);
\end{scope}
\end{tikzpicture}
\end{equation}
\end{minipage}
\caption{The $\Gal(\bk/\k)$-actions on the hexagon of a rational del
Pezzo surface of degree $6$.}\label{fig:action-dP6}
\end{figure}

The groups $\Aut(X)$ and $\Aut_\k(X)$ act by symmetries on the
hexagon of $X_{\bk}$ and $X$, respectively, which induces
homomorphisms
\[\Aut(X)\to\sym_3\times\Z/2,\qquad\Aut_\k(X)\stackrel{\hat{\rho}}\to\sym_3\times\Z/2.\]
We now go through the cases in Figure~\ref{fig:action-dP6}. We will
see that (\ref{dp6:1}), (\ref{dp6:6}), and (\ref{dp6:8}) admit a
birational morphism to $\p^2$ and that
(\ref{dp6:2}), (\ref{dp6:3}), (\ref{dp6:4}), and (\ref{dp6:5}) admit a
birational morphism to $\Ql^L$ or $\F_0$.

\subsection{The del Pezzo surfaces in Figures~\ref{fig:action-dP6}(\ref{dp6:1}), \ref{fig:action-dP6}(\ref{dp6:6}), and \ref{fig:action-dP6}(\ref{dp6:8})}

The following statement is classical over algebraically closed
fields and is proven analogously over a perfect field $\k$.

\begin{lem}\label{prop:DP1}
Let $X$ be a del Pezzo surface of degree $6$ such that
$\rho(\Gal(\bk/\k))=\{1\}$ as indicated in
Figure~\ref{fig:action-dP6}(\ref{dp6:1})
\begin{enumerate}
\item\label{DP1:1} Then $X$ is rational and isomorphic to
\[\{([x_0:x_1:x_2],[y_0:y_1:y_2]) \in \p^2_\k\times\p^2_\k \mid x_0y_0=x_1y_1=x_2y_2\}.\]
\item\label{DP1:2} The action of $\Aut_\k(X)$ on the hexagon of $X$ induces the split exact sequences
\[1\rightarrow T_2\rightarrow \Aut(X) \rightarrow \sym_3\times\Z/2\rightarrow 1,\quad 1\rightarrow T_2(\k)\rightarrow \Aut_\k(X) \stackrel{\hat{\rho}}\rightarrow \sym_3\times\Z/2\rightarrow 1\]
where $T_2$ is a $2$-dimensional split torus, $\Z/2$ is generated by
the image of
\[ ([x_0:x_1:x_2],[y_0:y_1:y_2])\mapsto([y_0:y_1:y_2],[x_0:x_1:x_2])\]
and $\sym_3$ is generated by the image of
\begin{align*} ([x_0:x_1:x_2],[y_0:y_1:y_2])\mapsto([x_1:x_0:x_2], [y_1:y_0:y_2])\\
([x_0:x_1:x_2],[y_0:y_1:y_2])\mapsto([x_0:x_2:x_1], [y_0:y_2:y_1]).
\end{align*}
\item\label{DP1:3} $X\to\ast$ is a $\Aut_\k(X)$-Mori fibre space.
\end{enumerate}
\end{lem}
\begin{proof}
Contracting three disjoint curves in the hexagon of $X$ yields a
birational morphism onto a del Pezzo surface $Z$ of degree $9$, and
since the images of the three contracted curves are rational points,
we have $Z\simeq\p^2$. Choosing the three points to be the
coordinate points yields (\ref{DP1:1}). Any element of
$\ker(\hat{\rho})$ is conjugate via the contraction to an element of
$\Aut_\k(\p^2)$ fixing the coordinate points and vice-versa, so
$\ker(\hat{\rho})\simeq T_2(\k)$. The generators given in
(\ref{DP1:2}) can be verified with straightforward calculations. It
follows that $\Aut_\k(X)$ acts transitively on the sides of the
hexagon, hence $X$ is an $\Aut_\k(X)$-Mori fibre space.
\end{proof}

Over $\bk$, all rational del Pezzo surfaces of degree $6$ are
isomorphic. Therefore, by Lemma~\ref{prop:DP1}, for any del Pezzo
surface $X$ of degree $6$, we have
$\rk\,\NS(X_{\bk})^{\Aut_{\bk}(X)}=1$ and hence $X$ is an
$\Aut(X)$-Mori fibre space. Moreover, $\Aut(X)$  is a $\k$-structure on
$(\bk^*)^2\rtimes(\sym_3\times\Z/2)$. We will however encounter two
rational del Pezzo surfaces of degree $6$ that are not
$\Aut_\k(X)$-Mori fibre spaces, see  Lemma~\ref{prop:DP2} and
Lemma~\ref{prop:DP3}.

\begin{lem}\label{prop:DP5}
Let $X$ be a rational del Pezzo surface of degree $6$ such that
$\rho(\Gal(\bk/\k))=\Z/3$ as indicated in
Figure~\ref{fig:action-dP6}(\ref{dp6:6})
\begin{enumerate}
\item\label{DP5:1}
There exists a point $p=\{p_1,p_2,p_3\}$ in $\p^2$ of degree $3$
with splitting field $L$ such that $\Gal(L/\k)\simeq\Z/3$ and such
that $X$ is isomorphic to the blow-up of $\p^2$ in $p$.
\item\label{DP5:4}  $X$ is isomorphic to the graph of a quadratic involution $\varphi_{p}\in\Bir_\k(\p^2)$ with base-point $p$, and any two such surfaces are isomorphic if and only if the corresponding field extensions are $\k$-isomorphic.
\item\label{DP5:2} The action of $\Aut_\k(X)$ on the hexagon of $X$ induces a split exact sequence
\[1\rightarrow\Aut_\k(\p^2,p_1,p_2,p_3)\rightarrow\Aut_\k(X)\stackrel{\hat{\rho}}\rightarrow \Z/6=\langle\hat{\rho}(\alpha),\hat{\rho}(\beta)\rangle\rightarrow1\]
where $\alpha$ is the lift of an element of
$\Aut_\k(\p^2,\{p_1,p_2,p_3\})$ of order $3$ and $\beta$ is the lift
of $\varphi_p$.
\item\label{DP5:3}
$X\to\ast$ is an $\Aut_\k(X)$-Mori fibre space.
\end{enumerate}
\end{lem}
\begin{proof}
(\ref{DP5:1}) The hexagon of $X$ is the union of two curves $C_1$
and $C_2$, each of whose three geometric components are disjoint.
For $i=1,2$, the contraction of $C_i$ yields a birational morphism
$\pi_i\colon X\rightarrow\p^2$ which contracts the curve onto a
point of degree $3$. By Lemma~\ref{rmk:cyclicpoint2} we can assume
it is the same point for $i=1,2$, which we call $p=\{p_1,p_2,p_3\}$.
It remains to see that $\Gal(L/\k)\simeq\Z/3$, where $L$ is any
splitting field of $p$. Since $\rho(\Gal(\bk/\k))\simeq\Z/3$, the
action of $\Gal(L/\k)$ on $\{p_1,p_2,p_3\}$ induces an exact
sequence $1\to H\to\Gal(L/\k)\to\Z/3\to 1$. The field $L':=\{a\in
L\mid h(a)=a\ \forall \ h\in H\}$ is an intermediate field between
$L$ and $\k$, over which $p_1,p_2,p_3$ are rational. The minimality
of $L$ implies that $L'=L$ and hence $H=\{1\}$ \cite[Corollary
2.10]{Morandi}.

(\ref{DP5:4}) The fact that any two such surfaces $X$ are isomorphic
if and only if the respective field extensions are $\k$-isomorphic
follows from Remark~\ref{rmk:cyclicpoint2}. The map
$\varphi_p:=\pi_2\pi_1^{-1}\in\Bir_\k(\p^2)$ is of degree $2$ and
$p$ is the base-point of $\varphi_p$ and $\varphi_p^{-1}$. By
Lemma~\ref{rmk:cyclicpoint2} we can assume that $\varphi_p$ has a
rational fixed point $r$ and  that it contracts the line through
$p_i,p_j$ onto $p_k$, where $\{i,j,k\}=\{1,2,3\}$. These conditions
imply that $\varphi_p$ is an involution, and by construction of
$\varphi_p$, the surface $X$ is isomorphic to the graph of
$\varphi_p$.

(\ref{DP5:2}) The kernel $\ker(\hat{\rho})$ is conjugate via $\pi_1$
to the subgroup of $\Aut_\k(\p^2)$ fixing $p_1,p_2,p_3$. The only
non-trivial elements of $\sym_3\times\Z/2$ commuting with
$\rho(\Gal(\bk/\k))$ are rotations, so
$\hat\rho(\Aut_\k(X))\subseteq\Z/6$. The involution
$\varphi_{p}\in\Bir_\k(\p^2)$ lifts to an automorphism $\beta$
inducing a rotation of order $2$. If $\langle\sigma\rangle=\Z/3$,
there exists $\tilde{\alpha}\in\Aut_\k(\p^2)$ such that
$\tilde{\alpha}(p_i)=p_{\sigma(i)}$, $i=1,2,3$, and
$\tilde{\alpha}(r)=r$, where $r$ is the fixed point of $\varphi_p$,
see Lemma~\ref{rmk:cyclicpoint2}. Then $\tilde{\alpha}^3$  and
$\tilde{\alpha}\varphi_p\tilde{\alpha}^{-1}\varphi_p$ are linear and
fix $r,p_1,p_2,p_3$, and hence $\tilde{\alpha}$ is of order $3$ and
$\tilde{\alpha}$ and $\varphi_p$ commute. The lift $\alpha$ of
$\tilde\alpha$ is an automorphism commuting with $\beta$ and
inducing a rotation of order $3$.

(\ref{DP5:3}) Since $\Aut_\k(X)$ contains an element inducing a
rotation of order $6$ on the hexagon, we have
$\rk\,\NS(X)^{\Aut_\k(X)}=1$.
\end{proof}

\begin{lem}\label{prop:DP7}
Let $X$ be a rational del Pezzo surface of degree $6$ such that
$\rho(\Gal(\bk/\k))=\sym_3$ as indicated in
Figure~\ref{fig:action-dP6}(\ref{dp6:8})
\begin{enumerate}
\item\label{DP7:1}
There exists a point $p=\{p_1,p_2,p_3\}$ in $\p^2$ of degree $3$
with splitting field $L$ such that $\Gal(L/\k)\simeq\sym_3$ and such
that $X$ is isomorphic to the blow-up of $\p^2$ in $p$.
\item\label{DP7:4}  $X$ is isomorphic to the graph of a quadratic involution $\varphi_{p}\in\Bir_\k(\p^2)$ with base-point $p$, and  any two such surfaces are isomorphic if and only if the corresponding field extensions are $\k$-isomorphic.
\item\label{DP7:2} The action of $\Aut_\k(X)$ on the hexagon of $X$ induces a split exact sequence
\[1\rightarrow\Aut_\k(\p^2,p_1,p_2,p_3)\rightarrow\Aut_\k(X)\stackrel{\hat{\rho}}\rightarrow \Z/2=\langle\hat{\rho}(\alpha)\rangle\rightarrow1\]
where $\alpha$ is the lift of $\varphi_p$ onto $X$.
\item\label{DP7:3}
$X\to\ast$ is an $\Aut_\k(X)$-Mori fibre space.
\end{enumerate}
\end{lem}
\begin{proof}
(\ref{DP7:1}) and (\ref{DP7:4}) are proven analogously to
Lemma~\ref{prop:DP5}(\ref{DP5:1}) and \ref{prop:DP5}(\ref{DP5:4}).

(\ref{DP7:2}) The kernel of $\hat\rho$ is conjugate to
$\Aut_\k(\p^2,p_1,p_2,p_3)$ via the birational morphism $X\to\p^2$
that contracts one curve in the hexagon of $X$ onto $p$. Any element
of $\Aut_\k(X)$ induces a symmetry of the hexagon that commutes with
the $\Gal(\bk/\k)$-action on the hexagon, hence
$\hat{\rho}(\Aut_\k(X))$ is contained in the factor $\Z/2$ generated
by a rotation of order $2$. The quadratic involution $\varphi_p$
lifts to an automorphism $\alpha$ of $X$ and $\hat{\rho}(\alpha)$ is
a rotation of order $2$.

(\ref{DP7:3}) Since $\hat{\rho}(\alpha)$ exchanges the two curves in
the hexagon, we have $\rk\,\NS(X)^{\Aut_\k(X)}=1$.
\end{proof}

\begin{ex}\label{ex:DP5}
A del Pezzo surface as in Lemma~\ref{prop:DP5} exists:  let $|\k|=2$
and $L/\k$ be the splitting field of $p(X)=X^3+X+1$, {\it i.e.}
$|L|=8$. Then $\sigma\colon a\mapsto a^2$ generates $\Gal(L/\k)$
\cite[Theorem 6.5]{Morandi}. If $\zeta$ a root of $P$, then
$\sigma(\zeta^4)=\zeta$ and hence the point
$\{[1:\zeta:\zeta^4],[1:\zeta^2:\zeta],[1:\zeta^4:\zeta^2]\}$ is of
degree $3$, its components are not collinear and they are cyclically
permuted by $\sigma$.
\end{ex}

\begin{ex}\label{ex:DP7}
A del Pezzo surface as in Lemma~\ref{prop:DP7} exists: let $\k=\Q$,
$\zeta:=2^{\frac{1}{3}}$ and $\omega=e^{\frac{2\pi i}{3}}$. Then
$L:=\Q(\zeta,\omega)$ is a Galois extension of $\Q$ of degree $6$
and $\Gal(L/\k)\simeq\sym_3$ is the group of $\k$-isomorphisms of
$L$ sending $(\zeta,\omega)$ respectively to $(\zeta,\omega)$,
$(\omega\zeta,\omega)$, $(\zeta,\omega^2),\ (\omega\zeta,\omega^2)$,
$(\omega^2\zeta,\omega)$, $(\omega^2\zeta,\omega^2)$
 \cite[Example 2.21]{Morandi}. The point $\{[\zeta:\zeta^2:1],[\omega\zeta:\omega^2\zeta^2:1],[\omega^2\zeta:\omega\zeta^2:1]\}$ is of degree $3$, its components are not collinear and any non-trivial element of $\Gal(L/\k)$ permutes them non-trivially.

A del Pezzo surfaces as in Lemma~\ref{prop:DP7} cannot exist over a
finite field, because Galois groups of finite extensions of
finite fields are always cyclic.
\end{ex}

\subsection{The del Pezzo surface in Figures~\ref{fig:action-dP6}(\ref{dp6:7}) and \ref{fig:action-dP6}(\ref{dp6:9})}
Recall that the two del Pezzo surfaces of degree $6$ in
Lemma~\ref{prop:DP5} and Lemma~\ref{prop:DP7} are the blow-up of a
point $p\in\p^2$ of degree $3$.

\begin{lem}\label{prop:DP6}
Let $X$ be a rational del Pezzo surface with
$\rho(\Gal(\bk/\k))=\Z/6$ as in
Figure~\ref{fig:action-dP6}(\ref{dp6:7}). Then $X\to\ast$ is a Mori
fibre space and
\begin{enumerate}
\item\label{DP6:1}
there exists a quadratic extension $L/\k$ such that $X_L$ is
isomorphic to the del Pezzo surface of degree $6$ from
Lemma~\ref{prop:DP5} (see Figure~\ref{fig:action-dP6}$(5)$), which
is the blow-up  $\pi\colon X_L\to\p^2_L$ of a point
$p=\{p_1,p_2,p_3\}$ of degree $3$ with splitting field $F$ such that
$\Gal(F/\k)\simeq\Z/3$.
\item\label{DP6:2}
$\pi\Gal(L/\k)\pi^{-1}$ acts rationally on $\p^2_L$; it is not
defined at $p$, sends a general line onto a conic through $p$ and
acts on $\Aut_L(\p^2,\{p_1,p_2,p_3\})$ by conjugation.
\item\label{DP6:3}
Any two such surfaces are isomorphic if and only if the
corresponding field extensions of degree two and three are
$\k$-isomorphic.
\item\label{DP6:4}
The action of $\Aut_\k(X)$ on the hexagon of $X$ induces a split
exact sequence
\[1\to \Aut_L(\p^2,p_1,p_2,p_3)^{\pi\Gal(L/\k)\pi^{-1}}\to\Aut_\k(X)\to\Z/6=\langle\hat\rho(\alpha),\hat\rho(\pi^{-1}\varphi_p\pi)\rangle\to1\]
where $\alpha$ is the lift of an element in
$\Aut_L(\p^2,\{p_1,p_2,p_3\})^{\pi\Gal(L/\k)\pi^{-1}}$ of order $3$
and $\varphi_p\in\Bir_L(\p^2)$ a quadratic involution with
base-point $p$.
\end{enumerate}
\end{lem}
\begin{proof}
All $(-1)$-curves of $X_{\bk}$ are in the same $\Gal(\bk/\k)$-orbit
and hence $X\to\ast$ is a Mori fibre space.

(\ref{DP6:1}) Since $X$ is rational, it contains a rational point
$r\in X(\k)$, see Proposition~\ref{lem:sing-fibres}, which is in
particular not contained in the hexagon of $X$. Let $\eta_1\colon
Y\to X$ be its blow-up and $E_r$ its exceptional divisor. Then
$Y_{\bk}$ contains an orbit of three $(-1)$-curves $C_1,C_2,C_3$
intersecting $E_r$, each intersecting two opposite sides of the
hexagon. The contraction of $C:=C_1\cup C_2\cup C_3$ yields a
birational morphism $\eta_2\colon Y\to Z$ onto a rational del Pezzo
surface of degree $8$. The birational map $\eta_2\eta_1^{-1}$
conjugates the $\Gal(\bk/\k)$-action on $Z$ to an action that
exchanges the fibrations of $Z_{\bk}$ and hence $Z\simeq \Ql^L$ for
some quadratic extension $L/\k$, by Lemma~\ref{prop:Q}(\ref{Q:1}).
Figure~\ref{fig:rot6} shows the action of $\rho(\Gal(\bk/\k))$ on
the image by $\eta_2\eta_1^{-1}$ of the hexagon of $X$.
\begin{figure}[h]
\[
\begin{tikzpicture}[scale=.8,font=\footnotesize]
\begin{scope}[every coordinate/.style={shift={(0,0)}}]
\path [c](0,0) pic {hexagon}; \draw[thick, shorten <=-.27cm, shorten
>=-.27cm] (D1) to (D4); \draw[thick, shorten <=-.2cm, shorten
>=-.2cm] (D2) to (D5); \draw[thick, shorten <=-.2cm, shorten
>=-.2cm] (D3) to (D6);
\node at (0,0){$\odot$}; \node at (-2,-1){$X$}; \node at
(-0.5,0){$r$}; \node at (1.8,0.8){$C$}; \node at
(3,0.3){$\eta_2\eta_1^{-1}$}; \draw[->,dashed](2,0) to (3.5,0);
\draw[->] (D11) to [bend right=15,swap,shorten >=.4cm,shorten
<=.4cm] (D2); \draw[->] (D22) to [bend right=15,swap,shorten
>=.2cm,shorten <=.4cm] (D3); \draw[->] (D33) to [bend
right=15,swap,shorten >=.2cm,shorten <=.4cm] (D4); \draw[->] (D44)
to [bend right=15,swap,shorten >=.2cm,shorten <=.4cm] (D5);
\draw[->] (D55) to [bend right=15,swap,shorten >=.2cm,shorten
<=.4cm] (D6); \draw[->] (D66) to [bend right=15,swap,shorten
>=.2cm,shorten <=.4cm] (D1);
\end{scope}
\begin{scope}[every coordinate/.style={shift={(5,0)}}]
\path [c](0,0) pic {quadrat}; \draw[shorten <=-.27cm, shorten
>=-.27cm] (B1) to (B3); \draw[shorten <=-.27cm, shorten >=-.27cm]
(B2) to (B4); \draw[dashed, shorten <=-.27cm, shorten >=-.27cm]
[c](A1) to (A3); \node at (A1){$\bullet$}; \node at (A3){$\bullet$};
\node[c] at (-1.25,0){$\bullet$}; \node[c] at (0.3,1.2){$E_r$};
\node[c] at (1,-1){$Z\simeq \Ql^L$}; \draw[->] [c](-0.5,0) to [bend
right=15,swap] (-1,0.5); \draw[->] [c](-1,0.5) to [bend
right=15,swap] (-0.5,1); \draw[->] [c](-0.5,1) to [bend
right=15,swap] (0,0.5); \draw[->] [c](0,-0.5) to [bend left=15,swap]
(0.5,-1); \draw[->] [c](0.5,-1) to [bend left=15,swap] (1,-0.5);
\draw[->] [c](1,-0.5) to [bend left=15,swap] (0.5,0);
\end{scope}
\end{tikzpicture}
\]
\caption{The $\Gal(\bk/\k)$-action on $Z_{\bk}\simeq
\Ql_{\bk}^L$.}\label{fig:rot6}
\end{figure}
Then $\eta_2\eta_1^{-1}$ conjugates the $\Gal(\bk/L)$-action on
$\Ql_L^L$ to an action on the hexagon with $\rho(\Gal(\bk/L))=\Z/3$.
Lemma~\ref{prop:DP5} implies (\ref{DP6:1}).

(\ref{DP6:3}) By Lemma~\ref{lem:cyclicpoint3}, $\Aut_\k(\Ql^L)$ acts
transitively on the set of points of degree $3$ in $\Ql^L$ with
$\k$-isomorphic splitting fields and whose geometric components are
in general position. This yields the claim.

(\ref{DP6:2}) Write $\Gal(L/\k)=\langle g\rangle$. Then $g$
exchanges opposite edges of the hexagon and thus $\rho_g:=\pi
g\pi^{-1}$ acts rationally on $\p^2$; it is not defined at $p$,
contracts the lines through any two of $p_1,p_2,p_3$ onto the third
of these three and it sends a general line onto a conic through $p$.
It follows that for $\beta\in\Aut_L(\p^2,\{p_1,p_2,p_3\})$ the map
$\rho_g\beta\rho_g$ is contained in $\Aut_L(\p^2)$ and preserves
$\{p_1,p_2,p_3\}$.

(\ref{DP6:4})
 The automorphisms of $X$ are the automorphisms of $X_{\bk}$ commuting with the $\Gal(\bk/\k)$-action, hence $\hat\rho(\Aut_\k(X))\subseteq\Z/6$.
 Since $X$ is rational, $\Gal(L/\k)$ has a fixed point $r\in X(\k)$.
Let $\varphi_p\in\Bir_L(\p^2)$ be the quadratic involution from
Lemma~\ref{prop:DP5}(\ref{DP5:2}) such that
$\Phi_p:=\pi^{-1}\varphi_p\pi\in\Aut_L(X)$ induces a rotation of
order $2$ on the hexagon of $X_L$. By Lemma~\ref{rmk:cyclicpoint2},
we can assume that $\varphi_p$ fixes $\pi(r)\in \p^2(L)$. Then
$\Phi_p g\Phi_p g\in\Aut_L(X)$, preserves the edges of the hexagon
and fixes $r$. It therefore descends to an element of
$\Aut_L(\p^2,p_1,p_2,p_3)$ fixing $r$ and is hence equal to the
identity. It follows that $\Phi_p\in\Aut_\k(X)$. By
Lemma~\ref{prop:DP5}(\ref{DP5:2}), there is an element of
$\tilde\alpha\in\Aut_L(\p^2,\{p_1,p_2,p_3\})$ of order $3$ inducing
a rotation of order $3$ on the hexagon of $X_L$, and again we can assume
that it fixes $\pi(r)\in \p^2(L)$. We argue as above that
$\alpha:=\pi^{-1}\tilde\alpha\pi\in\Aut_\k(X)$, and it follows that
the sequence is split. Finally, any element of $\ker(\hat\rho)$
preserves each edge of the hexagon and is therefore conjugate by
$\pi$ to an element of $\Aut_L(\p^2,p_1,p_2,p_3)$ commuting with
$\rho_g$, and any element of $\Aut_L(\p^2,p_1,p_2,p_3)^{\rho_g}$
lifts to an element of $\ker(\hat\rho)$.
\end{proof}

\begin{lem}\label{prop:DP8}
Let $X$ be a rational del Pezzo surface with
$\rho(\Gal(\bk/\k))=\sym_3\times\Z/2$ as in
Figure~\ref{fig:action-dP6}(\ref{dp6:9}). Then $X\to\ast$ is a Mori
fibre space and
\begin{enumerate}
\item\label{DP8:1}
there exists a quadratic extension $L/\k$ such that $X_L$ is
isomorphic to the del Pezzo surface of degree $6$ from
Lemma~\ref{prop:DP7} (see Figure~\ref{fig:action-dP6}$(7)$), which
is the blow-up  $\pi\colon X_L\to\p^2_L$ of a point
$p=\{p_1,p_2,p_3\}$ of degree $3$ with splitting field $F$ such that
$\Gal(F/\k)\simeq\sym_3$.
\item\label{DP8:2}
$\pi\Gal(L/\k)\pi^{-1}$ acts rationally on $\p^2$; it is not defined
at $p$, sends a general line onto a conic through $p$ and acts on
$\Aut_L(\p^2,\{p_1,p_2,p_3\})$ by conjugation.
\item\label{DP8:3}
Any two such surfaces are isomorphic if and only if the
corresponding field extensions of degree two and six are
$\k$-isomorphic.
\item\label{DP8:4}
The action of $\Aut_\k(X)$ on the hexagon of $X$ induces a split
exact sequence
\[1\to \Aut_L(\p^2,p_1,p_2,p_3)^{\pi\Gal(L/\k)\pi^{-1}}\to\Aut_\k(X)\to\Z/2=\langle\hat\rho(\pi^{-1}\varphi_p\pi)\rangle\to1,\]
where $\varphi_p\in\Bir_L(\p^2)$ is a quadratic involution with
base-point $p$.
\end{enumerate}
\end{lem}
\begin{proof}
This is proven analogously to Lemma~\ref{prop:DP6}.
\end{proof}

\begin{ex}\label{ex:4.8}
Rational del Pezzo surfaces of degree $6$ over $\k$ as in
Lemma~\ref{prop:DP6} and Lemma~\ref{prop:Q} exist: in
Example~\ref{ex:DP5} and Example~\ref{ex:DP7}, there is a point
$p\in\p^2$ of degree $3$ with a splitting field $F/\k$ that is
Galois over $\k$ such that $\Gal(F/\k)\simeq\Z/3$ or
$\Gal(F/\k)\simeq\sym_3$, and the blow-up $\pi\colon Y\to\p^2$ of
$p$ is a rational del Pezzo surface of degree $6$ as in
Figure~\ref{fig:action-dP6}(\ref{dp6:6}) or (\ref{dp6:8}). The point
$p$ is also a point of degree $3$ in $\p^2_L$ with splitting field
$FL/L$ because $\Gal(FL/L)\simeq\Gal(F/\k)$ \cite[Theorem
5.5]{Morandi}.

By Lemma~\ref{prop:DP5}(\ref{DP5:4}) and
Lemma~\ref{prop:DP7}(\ref{DP7:4}) there exists a quadratic
involution $\varphi_p\in\Bir_\k(\p^2)$ such that
$\Phi:=\pi^{-1}\varphi_p\pi\in\Aut_\k(Y)$ induces a rotation of
order $2$. By Lemma~\ref{rmk:cyclicpoint2}, we can assume that
$\varphi_p$ has a rational fixed point $r\in\p^2(\k)$. Let $g$ be
the generator of $\Gal(L/\k)$ and define $\psi_g:=\Phi\circ
g=g\circ\Phi$. The group $\langle\psi_g\rangle$ acts on $Y_L$ with
fixed point $\pi^{-1}(r)\in Y_L(L)$ and it induces a rotation of
order $2$ on the hexagon of $Y_L$. It follows that
$\Gal(L/\k)\simeq\langle\psi_g\rangle$ defines a $\k$-structure $X$ on
$Y_L$, which is rational by Proposition~\ref{lem:sing-fibres}. It
follows that the group $\Gal(\bk/\k)$ acts on the hexagon of $Y_L$
by $\Z/6$ or by $\sym_3\times\Z/2$.
\end{ex}

\subsection{The del Pezzo surfaces in Figures~\ref{fig:action-dP6}(\ref{dp6:3}) and \ref{fig:action-dP6}(\ref{dp6:4})}

\begin{lem}\label{prop:DP3}
Let $X$ be a del Pezzo surface of degree $6$ such that
$\rho(\Gal(\bk/\k))$ is generated by a reflection as indicated in
Figure~\ref{fig:action-dP6}(\ref{dp6:3}). Then $X$ is rational and
\begin{enumerate}
\item\label{DP3:1} there is a quadratic extension $L/\k$ and a birational morphism $\eta\colon X\to \Ql^L$ contracting the two $\k$-rational curves in the hexagon onto $p_1=([1:0],[1:0])$ and $p_2=([0:1],[0:1])$.
\item\label{DP3:4} Any two such surfaces are isomorphic if and only if the respective quadratic extensions are $\k$-isomorphic.
\item\label{DP3:2}
The action of $\Aut_\k(X)$ on the hexagon of $X$ induces a split
exact sequence
\[
1\rightarrow T^L(\k)\to \Aut_\k(X)
\stackrel{\hat{\rho}}\to\langle\hat\rho(\alpha)\rangle\times\langle\hat\rho(\beta)\rangle\rightarrow
1,
\]
where $\eta T^L(\k)\eta^{-1}\subseteq\Aut_\k(\Ql^L,p_1,p_2)$ is the
subgroup preserving the ruling of $\Ql^L$, and the automorphisms
$\alpha\colon(u,v)\mapsto (\frac{1}{v},\frac{1}{u})$ and
$\beta\colon(u,v)\mapsto(\frac{1}{u},\frac{1}{v})$.
\item\label{DP3:3} $\rk\,\NS(X)^{\Aut_\k(X)}=2$ and
$\eta\Aut_\k(X)\eta^{-1}=\Aut_\k(\Ql^L,\{p_1,p_2\})$. In particular,
$X\to\ast$ is not an $\Aut_\k(X)$-Mori fibre space.
\end{enumerate}
\end{lem}
\begin{proof}
(\ref{DP3:1}) The hexagon of $X$ has exactly two $\k$-rational curves
$C_1,C_2$, which are moreover disjoint. Their contraction yields a birational morphism
$\eta\colon X\to Z$ onto a del Pezzo surface $Z$ of degree $8$ with
two rational points. By Proposition~\ref{lem:sing-fibres}, $Z$ is
rational and by Lemma~\ref{prop:Q}(\ref{Q:1}) we have $Z\simeq
\Ql^L$. We can assume that $C_1,C_2$ are contracted onto
$p_1=([1:0],[1:0])$ and $p_2=([0:1],[0:1])$ by
Lemma~\ref{lem:Q-deg2}(\ref{Q-deg2:2}).

(\ref{DP3:4}) Any two rational points on $\Ql^L$ that are not on the
same ruling of $\Ql_L^L$ can be sent onto each other by an
element of $\Aut_\k(\Ql^L)$ by
Lemma~\ref{lem:Q-deg2}(\ref{Q-deg2:2}). It follows that any two del
Pezzo surfaces satisfying our hypothesis are isomorphic if and only
if they have a birational contraction to isomorphic del Pezzo
surfaces $\Ql^L$ and $\Ql^{L'}$ of degree $8$. This is the case if
and only if $L$ and $L'$ are $\k$-isomorphic by
Lemma~\ref{prop:Q}(\ref{Q:3}).

(\ref{DP3:2}) The kernel of $\hat\rho$ is the subgroup of
$\Aut_\k(X)$ of elements preserving $C_1,C_2$ and hence its conjugate
$\eta\ker(\hat{\rho})\eta^{-1}\subseteq\Aut_\k(\Ql^L,p_1,p_2)$ is
the subgroup preserving the rulings of $\Ql^L$. The only non-trivial
automorphisms of $X_{\bk}$ commuting with the $\Gal(\bk/\k)$-action
induce a rotation of order $2$ or a reflection that preserves
$C_1\cup C_2$. Let $L/\k$ be an extension of degree $2$ such that
$\Ql^L_L\simeq\p^1_L\times\p^1_L$. The involution
$\alpha\in\Aut_\k(\Ql^L)$ exchanges $p_1,p_2$ and the rulings of
$\Ql^L_L$, it thus lifts to an automorphism of $X$ inducing a
reflection. The involution $\beta\in\Aut_\k(\Ql^L)$ exchanges
$p_1,p_2$ and preserves the rulings of $\Ql^L_L$, it thus lifts to
an involution of $X$ inducing a rotation of order $2$ on the
hexagon. The involutions $\alpha,\beta\in\Aut_\k(\Ql^L)$ commute,
hence their lifts commute, which yields the splitness of the
sequence.

(\ref{DP3:3}) It follows from (\ref{DP3:2}) that any automorphism of
$X$ preserves $C_1\cup C_2$, and since
$\eta^{-1}\alpha\eta\in\Aut_\k(X)$ exchanges $C_1,C_2$, we have
$\rk\,\NS(X)^{\Aut_\k(X)}=2$.
\end{proof}

The $\R$-version of Lemma~\ref{prop:DP3}(\ref{DP3:2}) in
\cite[Proposition 3.4]{RZ} states that the kernel is
$\mathrm{SO}(\R)$, but it should be
$T_Q(\R)\simeq\mathrm{SO}(\R)\times\R_{>0}$.

\begin{lem}\label{prop:DP4}
Let $X$ be a rational del Pezzo surface of degree $6$ such that
$\rho(\Gal(\bk/\k))$ is generated by a rotation of order $2$ as
indicated in Figure~\ref{fig:action-dP6}(\ref{dp6:4}). Then there
exists a quadratic extension $L=\k(a_1)$ of $\k$ such that
\begin{enumerate}
\item\label{DP4:1}
$X$ is isomorphic to the blow-up of $\F_0$ in the point
$\{[a_1:1;a_1:1],[a_2:1;a_2:1]\}$ of degree $2$ and
\[
X\simeq \{([u_0:u_1],[v_0:v_1],[w_0:w_1])\in(\p^1)^3\mid
w_0\tilde{a}(u_0v_0+au_1v_0+\tilde{a}u_1v_1)=w_1(u_0v_1-x_1v_0)\}
\]
where $t^2+at+\tilde{a}=(t-a_1)(t-a_2)\in\k[t]$ is the minimal
polynomial of $a_1$ over $\k$.
\item\label{DP4:4} Any two such surfaces are isomorphic if and only if the respective quadratic extensions are $\k$-isomorphic.
\item\label{DP4:2} The action of $\Aut_\k(X)$ on the hexagon induces an exact sequence,
\[1\rightarrow \Aut_\k(\p^1,p_1,p_2)^2\rightarrow\Aut_\k(X)\stackrel{\hat\rho}\rightarrow \sym_3\times\Z/2\rightarrow 1,\]
which is split if $\mathrm{char}(\k)\neq2$,
$\Z/2=\langle\hat\rho(\tilde\alpha)\rangle$ and
$\sym_3=\langle\hat\rho(\tilde\beta),\hat\rho(\tilde\varphi)\rangle$,
where $\tilde\alpha,\tilde\beta,\tilde\varphi$ are the lifts of the
involutions of $\F_0$
\begin{align*}
\alpha\colon[y_0:y_1;z_0:z_1]\mapsto&[y_0+ay_1:-y_1; z_0+az_1:-z_1],\\
\beta\colon[y_0:y_1;z_0:z_1]\mapsto&[z_0:z_1;y_0:y_1],\\
\psi\colon[y_0:y_1;z_0:z_1]\dashmapsto&[y_0+ay_1:-y_1;
\tilde{a}(y_1z_0-y_0z_1):y_0z_0+ay_0z_1+\tilde{a}y_1z_1].
\end{align*}
\item\label{DP4:3} $X\to\ast$ is an $\Aut_\k(X)$-Mori fibre space.
\end{enumerate}
\end{lem}
\begin{proof}
(\ref{DP4:1}) Let $C_1,C_2,C_3$ be the curves in the hexagon of $X$.
By Lemma~\ref{prop:Q}(\ref{Q:1}), for $i=1,2,3$, there is a
birational morphism $\pi_i\colon X\rightarrow \F_0$ contracting
$C_i$ onto a point of degree $2$. Let $L/\k$ be a quadratic
extension such that $\Gal(L/\k)$ acts by the rotation of order $2$.
Then $\Gal(\bk/L)$ preserves each $C_i$, hence $L$ is the splitting
field of each $C_i$. So, $L$ is also the splitting field of each
$\pi_i(C_i)$. Let $L=\k(a_1)$ for some $a_1\in L$. For $i=1,2,3$ we
write $\pi_i(C_i)=\{[b_{i1}:1;b_{i2}:1],[b_{i3}:1;b_{i4}:1]\}$ for
some $b_{i1},\dots,b_{i4}\in L$. Since the two components of
$\pi_i(C_i)$ are not contained in the same fibre of $\F_0$,
Remark~\ref{rmk:cyclicpoint2-dim1} implies that there is
$A_i\in\PGL_2(\k)$ that sends $[b_{i1}:1],[b_{i3}:1]$ onto
$[a_1:1],[a_2:1]$. Similarly, there is $B_i\in\PGL_2(\k)$ that sends
$[b_{i2}:1],[b_{i4}:1]$ onto $[a_1:1],[a_2:1]$. Up to changing the
rulings on $\F_0$, we can assume that
$\varphi:=\pi_2\pi_1^{-1}\colon\F_0\rat\F_0$ preserves the ruling
given by the first projection, as indicated in the following
commutative diagram.
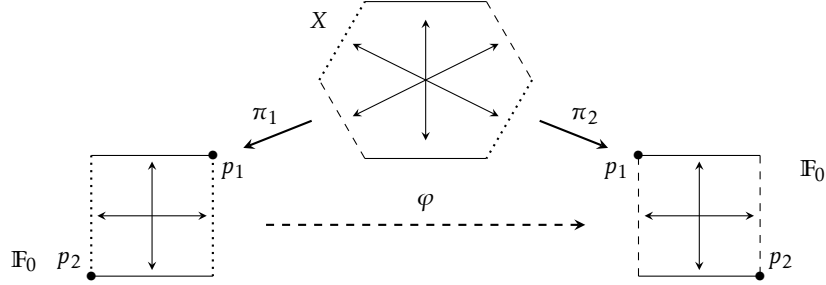
\[
\begin{tikzpicture}[scale=.6,font=\footnotesize]
\begin{scope}
\path (0,0) pic {hexagon-leer}; \node at (n){$X$};
\draw[dashed] (E1) to (E2); \draw[dotted, thick] (E2) to (E3); \draw
(E3) to (E4); \draw[dashed] (E4) to (E5); \draw[dotted, thick]  (E5)
to (E6); \draw (E6) to (E1); 
\draw[<->] (D1) to (D4); \draw[<->] (D2) to (D5); \draw[<->] (D3) to
(D6); 
\draw[->,thick] (-2.5,-0.9) to (-4,-1.5); \draw[->,thick] (2.5,-0.9)
to (4,-1.5); \node at (-3.5,-0.8){$\pi_1$}; \node at
(3.5,-0.8){$\pi_2$};
\end{scope}
\begin{scope}[every coordinate/.style={shift={(-6,-3)}}]
\path [c](0,0) pic {quadrat-leer};
\draw (A4) to (A1); \draw[dotted,thick] (A1) to (A2); \draw (A2) to
(A3); \draw[dotted,thick] (A3) to (A4); \draw[<->] (B1) to (B3);
\draw[<->] (B2) to (B4);
\node at (A1){$\bullet$}; \node at (A3){$\bullet$}; \coordinate (P1)
at (1.8,1){}; \coordinate (P2) at (-1.8,-1){}; \coordinate (F0) at
(-2.8,-1){}; \node at ([c]P1){$p_1$}; \node at ([c]P2){$p_2$}; \node
at ([c]F0){$\F_0$};
\end{scope}
\begin{scope}[every coordinate/.style={shift={(6,-3)}}]
\path [c](0,0) pic {quadrat-leer};
\draw (A4) to (A1); \draw[dashed] (A1) to (A2); \draw (A2) to (A3);
\draw[dashed] (A3) to (A4); \draw[<->] (B1) to (B3); \draw[<->] (B2)
to (B4);
\coordinate (P3) at (-1.8,1){}; \coordinate (P4) at (1.8,-1){};
\node at (A4){$\bullet$}; \node at (A2){$\bullet$}; \node at
([c]P3){$p_1$}; \node at ([c]P4){$p_2$}; \coordinate (F0) at
(2.5,1){}; \node at ([c]F0){$\F_0$};
\end{scope}
\begin{scope}[every coordinate/.style={shift={(0,-3.2)}}]
\coordinate (a1) at (-3.5,0); \coordinate (a2) at (3.5,0);
\coordinate (psi) at (0,0.5); \draw[->,dashed,thick] ([c]a1) to
([c]a2); \node at ([c]psi){$\varphi$};
\end{scope}
\end{tikzpicture}
\]
Up to an isomorphism of the first factor, we can assume that
$\varphi$ induces the identity map on $\p^1$. It then sends a
general fibre $f$ of the second projection onto a curve of bidegree
$(1,1)$ passing through $q$, which is given by
$\lambda(y_0z_1-y_1z_0)+\mu(y_0z_0+ay_1z_0+\tilde{a}y_1z_1)=0$ for
some $[\lambda:\mu]\in\p^1$. So, up to left-composition by an
automorphism of the second factor, $\varphi$ is the involution given
by
\[\varphi\colon [y_0:y_1;z_0:z_1]\dashmapsto[y_0:y_1;\tilde{a}(y_0z_1-y_1z_0):y_0z_0+ay_1z_0+\tilde{a}y_1z_1].\]
By construction of $\varphi$, $X$ is isomorphic to its graph inside
$(\p^1)^4$. The projection forgetting the third factor induces the
isomorphism in (\ref{DP4:1}).

(\ref{DP4:4}) As indicated in (\ref{DP4:1}), any two points of
degree $2$ in $\F_0$ whose geometric components are not in the same
ruling can be sent onto each other by an element of $\Aut_\k(\F_0)$.
It follows that two del Pezzo surfaces $X$ and $X'$ satisfying the
hypothesis of our lemma are isomorphic if and only if there are
contractions $X\to\F_0$ and $X'\to\F_0$ that contract a curve in
each hexagon onto points with $\k$-isomorphic splitting fields. This
is equivalent to contracted curves having $\k$-isomorphic splitting
fields.

(\ref{DP4:2}) The group $\pi_1\ker(\hat\rho)\pi_1^{-1}$ is the
subgroup of $\Aut_\k(\F_0)$ fixing $[a_i:1;a_i:1]$ for $i=1,2$ and
preserving the fibration given by the first projection, hence
$\pi_1\ker(\hat\rho)\pi_1^{-1}\simeq\Aut_\k(\p^1,[a_1:1],[a_2:1])^2$.
The involution $\alpha\in\Aut_\k(\F_0)$ (it is not the identity map
by Remark~\ref{rmk:minpoly_quadext}) preserves the fibrations of
$\F_0$ and exchanges $[a_1:1;a_1:1]$ and $[a_2:1;a_2:1]$. Thus it
lifts to an involution $\tilde{\alpha}\in\Aut_\k(X)$ inducing a
rotation of order $2$ on the hexagon. The involution
$\beta\in\Aut_\k(\F_0)$ exchanges the fibrations of $\F_0$ and fixes
$[a_i:1;a_i:1]$ for $i=1,2$, thus lifts to an involution
$\tilde{\beta}\in\Aut_\k(X)$ inducing the reflection at the axis
through $C_1$. We check that $\psi:=\varphi\circ\alpha$. Since
$\varphi$ induces the reflection on the hexagon that exchanges the
components of $C_3$, $\psi$ induces the reflection preserving each
component of $C_3$. It follows that the sequence is exact. If
$\mathrm{char}(\k)\neq2$, we have $a=0$, and then $\psi$ is an
involution, $\alpha$ commutes with $\beta$ and $\psi$, and
$\beta\circ\psi$ has order $3$. It follows that the sequence is
split.

(\ref{DP4:3}) Since $\Aut_\k(X)$ acts transitively on the edges of
the hexagon, $X\to\ast$ is an $\Aut_\k(X)$-Mori fibre space.
\end{proof}

\subsection{The del Pezzo surfaces in Figures~\ref{fig:action-dP6}(\ref{dp6:2}) and \ref{fig:action-dP6}(\ref{dp6:5})}

Here, we consider the remaining two del Pezzo surfaces of degree $6$
from Figure~\ref{fig:action-dP6}. We will see that none of them is a
$\Aut_\k(X)$-Mori fibre space. However, they carry a conic
fibration, and we will describe the automorphism group preserving
the fibration in this section, which will be used in the
Section~\ref{sec:CB}.

\begin{lem}\label{prop:DP2}
Let $X$ be a rational del Pezzo surface of degree $6$ such that
$\rho(\Gal(\bk/\k))$ is generated by a reflection as indicated in
Figure~\ref{fig:action-dP6}(\ref{dp6:2}). There exists a quadratic
extension $L=\k(a_1)/\k$ such that the following holds:
\begin{enumerate}
\item\label{DP2:4} There is a birational morphism $\eta\colon X\to \Rl^L\simeq\Ql^L$ contracting an irreducible $E$ curve onto the point $\eta(E)=\{[a_1^2:a_1:1:0],[a_2^2:a_2:1:0]\}=\{p_1,p_2\}$ of degree $2$.
\item\label{DP2:5} $X\simeq\{([w:x:y:z],[u:v])\mid v(w+ax+\tilde{a}z)=uy\}\subset\Rl^L\times\p^1$
\item\label{DP2:2} The action of $\Aut_\k(X)$ on the hexagon of $X$ induces a split exact sequence
\[
1\rightarrow
T^{L,L}(\k)\to\Aut_\k(X)\stackrel{\hat{\rho}}\to\langle\hat{\rho}(\alpha)\rangle\times\langle\hat{\rho}(\beta)\rangle\rightarrow1
\]
where $T^{L,L}(\k)\subset\Aut_\k(\Rl^L,p_1,p_2)$ is the subgroup
preserving the rulings of $\Rl^L_L$, and $\hat\rho(\alpha)$ is the
reflection exchanging the singular fibres and $\hat\rho(\beta)$ is a
rotation of order $2$ with
\begin{align*}
\eta\alpha\eta^{-1}\colon&[w:x:y:z]\mapsto[w:x+ay:-y:z]\\
\eta\beta\eta^{-1}\colon&[w:x:y:z]\mapsto[w+a(2x+az+ay):-(x+az):-y:z]
\end{align*}
where $t^2+at+\tilde{a}=(t-a_1)(t-a_2)\in\k[t]$ is the minimal
polynomial of $a_1$ over $\k$.
\item\label{DP2:3} We have $\rk\,\NS(X)^{\Aut_\k(X)}=2$ and $\eta\Aut_\k(X)\eta^{-1}=\Aut_\k(\Rl^L,\{p_1,p_2\})$.
In particular, $X\to\ast$ is not an $\Aut_\k(X)$-Mori fibre space.
\end{enumerate}
\end{lem}
\begin{proof}
(\ref{DP2:4}) By Lemma~\ref{prop:Q}(\ref{Q:1}), contracting $E$
yields a birational morphism $\nu\colon X\to \Ql^L$. The splitting
field of the image of $E$ is $L$, so we can choose
$\nu(E)=\{([1:0],[0:1]),([0:1],[1:0])\}$ by
Lemma~\ref{lem:Q-deg2}(\ref{Q-deg2:1}). Changing the model of
$\Ql^L$ with the isomorphism from
Lemma~\ref{lem:Q-deg22}(\ref{Q-deg22:2}), we get the birational
morphism $\eta\colon X\to\Rl^L$ and
$\eta(E)=\{[a_1^2:a_1:1:0],[a_2^2:a_2:1:0]\}$.

(\ref{DP2:3}) Any element of $\Aut_\k(X)$ preserves $E$. It follows
that $\rk\,\NS(X)^{\Aut_\k(X)}=2$ and that
$\nu\Aut_\k(X)\nu^{-1}=\Aut_\k(\Ql^L,\{p_1,p_2\})$.

(\ref{DP2:2}) The conjugate
$\nu\ker(\hat{\rho})\nu^{-1}\subseteq\Aut_\k(\Ql^L,([1:0],[0:1]),([0:1],[1:0]))$
is the subgroup preserving the rulings of $\Ql^L_L$. The only
non-trivial symmetries in $\sym_3\times\Z/2$ commuting with the
$\rho(\Gal(\bk/\k))$-action are the two reflections preserving $E$
and the rotation of order $2$. By Remark~\ref{rmk:minpoly_quadext},
$\eta\alpha\eta^{-1},\eta\beta\eta^{-1}$ are involutions and they
commute. Moreover, they respectively fix and exchange
$[a_1^2:a_1:1:0],[a_2^2:a_2:1:0]$. Their conjugates by the
isomorphism $\Rl^L\to\Ql^L$ from
Lemma~\ref{lem:Q-deg22}(\ref{Q-deg22:2}) respectively exchange and
preserve the rulings of $\Ql_L^L$. In particular, they induce the
claimed action on the hexagon of $X$, thus the sequence is split.
\end{proof}

\begin{lem}\label{prop:DP9}
Let $X$ be a rational del Pezzo surface of degree $6$ such that
$\rho(\Gal(\bk/\k))\simeq\Z/2\times\Z/2$ is generated by a reflection
and a rotation of order $2$ as in
Figure~\ref{fig:action-dP6}(\ref{dp6:5}). Then there exist quadratic
extensions $L=\k(a_1)$ and $L'=\k(b_1)$ of $\k$ that are not
$\k$-isomorphic, with
\[t^2+at+\tilde{a}=(t-a_1)(t-a_2),\quad t^2+bt+\tilde{b}=(t-b_1)(t-b_2)\in\k[t]\]
the minimal polynomials of $a_1,b_1$ such that the following hold:
\begin{enumerate}
\item\label{DP9:1}
 $X\simeq\Sl^{L,L'}$ and there exists a birational contraction $\eta\colon X\to \Ql^L\simeq\Rl^L$ contracting an irreducible curve onto the point $\{p_1,p_2\}=\{[b_1^2:b_1:0:1],[b_2^2:b_2:0:1]\}$ of degree $2$.
\item\label{DP9:2}
$X\simeq\{([w:x:y:z],[u:v])\mid
v(w+bx+\tilde{b}z)=uy\}\subset\Rl^L\times\p^1$
\item\label{DP9:3} Two surfaces $\Sl^{L,L'}$ and $\Sl^{\tilde{L},\tilde{L}'}$ are isomorphic if and only if $\tilde{L},\tilde{L}'$ are respectively $\k$-isomorphic to $L,L'$.
\item\label{DP9:4} The action of $\Aut_\k(X)$ on the hexagon of $X$ induces a split exact sequence
\[
1\to T^{L,L'}\to \Aut_\k(X)\stackrel{\hat\rho}\to
\langle\hat\rho(\alpha)\rangle\times\langle\hat\rho(\beta)\rangle\to
1
\]
where $ T^{L,L'}\subset\Aut_\k(\Rl^L,p_1,p_2)$ is the subgroup
preserving the rulings of $\Rl^L_L$, and $\hat\rho(\alpha)$ is the
reflection exchanging the singular fibres and $\hat\rho(\beta)$ is a
rotation of order $2$, where
\begin{align*}
\eta\alpha\eta^{-1}\colon&[w:x:y:z]\mapsto[w:x+ay:-y:z]\\
\eta\beta\eta^{-1}\colon&[w:x:y:z]\mapsto[w+b(2x+bz+ay):-(x+bz):-y:z]
\end{align*}
\item\label{DP9:5} $\rk\,\NS(X)^{\Aut_\k(X)}=2$ and $\eta\Aut_\k(X)\eta^{-1}=\Aut_\k(\Rl^L,\{p_1,p_2\})$. In particular, $X\to\ast$ is not an $\Aut_\k(X)$-Mori fibre space.
\end{enumerate}
\end{lem}
\begin{proof}
(\ref{DP9:1}) The hexagon of $X$ contains a unique curve $E$ whose
geometric components are disjoint. The contraction of $E$ yields a
birational morphism $\eta\colon X\to Y$ to a del Pezzo surface $Y$
of degree $8$, and the figure below shows the induced
$\Gal(\bk/\k)$-action on the image of the hexagon, so $Y\simeq\Ql^L$
for some quadratic extension $L/\k$ by
Lemma~\ref{prop:Q}(\ref{Q:1}).
\[
\begin{tikzpicture}[scale=.6,font=\footnotesize]
\begin{scope}
\path (0,0) pic {hexagon-leer}; \node at (n){$X$};
\draw[very thick] (E1) to (E2); \draw (E2) to (E3); \draw (E3) to
(E4); \draw[very thick] (E4) to (E5); \draw (E5) to (E6); \draw (E6)
to (E1); 
\draw[<->] (D1) to [bend right=15,swap] (D3); \draw[<->] (D6) to
[bend left=15,swap] (D4); \draw[<->] (D1) to [bend left=15,swap]
(D6); \draw[<->] (D3) to [bend right=15,swap] (D4); \draw[<->] (D2)
to (D5); \draw[<->,  shorten <=.2cm, shorten >=.2cm] (D1) to (D4);
\draw[<->,  shorten <=.2cm, shorten >=.2cm] (D3) to (D6); 
\draw[->] (2.7,0) to (5,0); \node at (3.5,.5){$\eta$};
\end{scope}
\begin{scope}[every coordinate/.style={shift={(7,0)}}]
\path [c](0,0) pic {quadrat-leer};
\draw (A4) to (A1); \draw(A1) to (A2); \draw (A2) to (A3); \draw
(A3) to (A4); \draw[<->] (B1) to (B3); \draw[<->] (B2) to (B4);
\draw[<->, shorten <=.2cm, shorten >=.2cm] (B1) to (B4);
\draw[<->,shorten <=.2cm, shorten >=.2cm] (B2) to (B3);
\draw[<->,shorten <=.2cm, shorten >=.2cm] (B1) to  (B2);
\draw[<->,shorten <=.2cm, shorten >=.2cm] (B3) to  (B4);
\draw[<->,shorten <=.2cm, shorten >=.2cm] (A1) to  (A3);
\node at (A1){$\bullet$}; \node at (A3){$\bullet$}; \coordinate (F0)
at (2.5,1){}; \node at ([c]F0){$\Ql^L$};
\end{scope}
\end{tikzpicture}
\]
We have $\rho(\Gal(\bk/\k))=\{1,r,s,rs\}$, where $r$ is the rotation
of order $2$ and $s$ is the reflection preserving the components of
$E$. Then $s$ or $sr$ is the image of the generator $g$ of
$\Gal(L/\k)$. It follows that the splitting field of $p$ is a
quadratic extension $L'/\k$ not $\k$-isomorphic to $L$ such that the
generator $g'$ of $\Gal(L'/\k)$ induces the rotation $r$ on the
hexagon. We set $L=\k(a_1)$ and $L'=\k(b_1)$ for some $a_1\in L$,
$b_1\in L'$. We can choose the form of $p$ according to
Lemma~\ref{lem:Q-deg22}(\ref{Q-deg22:3}).

(\ref{DP9:2}) follows from (\ref{DP9:1}) and
Lemma~\ref{lem:Q-deg22}(\ref{Q-deg22:4}).

(\ref{DP9:3}) Consider the birational morphism $\eta'\colon
\Sl^{\tilde{L},\tilde{L}'}\to\Rl^{\tilde{L}}$ with exceptional curve
$E'$. Suppose that we have $\Sl^{\tilde{L},\tilde{L}'}\simeq\Sl^{L,L'}.$
Then $E$ and $E'$ are the unique curves in the hexagon with only two
components. Thus they are defined over the same splitting field over
$\k$, and hence $L'\simeq\tilde{L}'$ over $\k$. It follows that
$\Rl^{L}\simeq\Rl^{\tilde{L}}$, which implies that $L\simeq\tilde L$
over $\k$ by Lemma~\ref{prop:Q}(\ref{Q:3}).

(\ref{DP9:4}--\ref{DP9:5}) The group
$\ker(\hat{\rho})\simeq\eta\ker(\hat{\rho})\eta^{-1}\subset\Aut_\k(\Ql^L,p_1,p_2)$
is the subgroup preserving the rulings of $\Ql^L$. Every element of
$\Aut_\k(X)$ preserves $E$ because it is the only curve in the
hexagon with only two geometric components, so the elements of
$\Aut_\k(X)$ act by symmetries of order $2$, and we have
$\eta\Aut_\k(X)\eta^{-1}=\Aut_\k(\Ql^L,\{p_1,p_2\})$. The only
symmetries of the hexagon that commute with $\rho(\Gal(\bk/\k))$ are
the two reflections preserving $E$ and the rotation of order $2$. By
Remark~\ref{rmk:minpoly_quadext},
$\eta\alpha\eta^{-1},\eta\beta\eta^{-1}$ are involutions and they
commute. Moreover, they respectively fix and exchange
$[b_1^2:b_1:1:0],[b_2^2:b_2:1:0]$. We see that the conjugates of
$\eta\alpha\eta^{-1},\eta\beta\eta^{-1}$ by the isomorphism
$\Rl^L\rat\Ql^L$ from Lemma~\ref{lem:Q-deg22}(\ref{Q-deg22:2})
respectively exchange and preserve the rulings of $\Ql_L^L$. In
particular, they induce the claimed action on the hexagon, thus the
sequence is split.
\end{proof}

\subsection{The fibration on a rational del Pezzo surface of degree $6$ from Figures~\ref{fig:action-dP6}(\ref{dp6:2}) and \ref{fig:action-dP6}(\ref{dp6:5})}\label{ss:fibrationDP}

Let $L/\k$, $L'/\k$ be two extensions of degree $2$. We can obtain
the Mori fibre space $\pi\colon \Sl^{L,L'}\to \p^1$ from
Example~\ref{ex:cb}(\ref{ex:Sl}) as follows: we first blow up the
point $p$, then contract the line passing through it, which yields a
birational map $\p^2\rat\Ql^L$. Since $p,p'$ are not collinear, the
image of $p'$ in $\Ql^L$ is a proper point and blowing it up yields
$\Sl^{L,L'}$. In particular, $\Sl^{L,L'}$ is one of the del Pezzo
surfaces in Figure~\ref{fig:action-dP6}(\ref{dp6:2}) and (\ref{dp6:5}),
which are described in Lemma~\ref{prop:DP2} and
Lemma~\ref{prop:DP9}.

\begin{rmk}\label{rmk:min-cb-sing}
\leavevmode
\begin{enumerate}
\item\label{rmk-min-cb-sing:1}
Let $L=\k(a_1)$ and $L'=\k(b_1)$ be two quadratic extensions of
$\k$, not necessarily non-isomorphic over $\k$, and let
\[t^2+at+\tilde{a}=(t-a_1)(t-a_2),\quad t^2+bt+\tilde{b}=(t-b_1)(t-b_2)\in\k[t]\]
be the minimal polynomials of $a_1$ and $b_1$ over $\k$.
Lemma~\ref{prop:DP2}(\ref{DP2:5}) and
Lemma~\ref{prop:DP9}(\ref{DP9:2}) imply that
\[\Sl^{L,L'}\simeq\{([w:x:y:z],[u:v])\in\p^3\times\p^1\mid wz=x^2+axy+\tilde{a}y^2, \ (w+bx+\tilde{b}z)v=uy\}\]
and the fibration $\pi\colon\Sl^{L,L'}\to\p^1$ is given by the
projection
\[([w:x:y:z],[u:v])\mapsto[u:v]=[w+bx+\tilde{b}z:y].\]
\item\label{rmk-min-cb-sing:2} The group $\Aut(\Sl^{L,L'},\pi)$ preserves a unique irreducible curve $E$ in the hexagon of $X$ that has disjoint geometric components. It induces a morphism
\[\Aut(\Sl^{L,L'},\pi)\to\Z/2,\]
and we denote by $\SlO^{L,L'}\subset\Aut(\Sl^{L,L'},\pi)$ its
kernel.
\item Via the contraction $\eta\colon X\to\Ql^L\simeq\Rl^L$
of $E$ onto a point $\{p_1,p_2\}$ of degree $2$,  the group
$\SlO^{L,L'}$ is conjugate to a subgroup of $T^{L,L'}$, the subgroup
of $\Aut(\Ql^L,p_1,p_2)$ preserving the rulings of $\Ql^L_L$
(see Lemma~\ref{prop:DP2}(\ref{DP2:2}) and
Lemma~\ref{prop:DP9}(\ref{DP9:4})).
\item The image $t,s\in\p^1(L)$ of the singular fibres make up two points of degree $1$ if $L, L'$ are $\k$-isomorphic, and one point of degree $2$ if $L,L'$ are not $\k$-isomorphic.
\end{enumerate}
\end{rmk}

\begin{lem}\label{lem:min-cb-sing1}
Keep the notation of Remark~\ref{rmk:min-cb-sing} and let $g$ be the
generator of $\Gal(L/\k)$. Then the action of $\Aut(\Sl^{L,L'}/\pi)$
on the geometric components of $E$ induces the split exact sequences
\begin{align*}
1\rightarrow \SlO^{L,L'}\to&\Aut(\Sl^{L,L'}/\pi)\to\Z/2\rightarrow 1\\
1\rightarrow
\SlO^{L,L'}(\k)\to&\Aut_\k(\Sl^{L,L'}/\pi)\to\Z/2\rightarrow 1
\end{align*}
where $\Z/2$ is generated by the image of the involution
\[ ([w:x:y:z],[u:v])\mapsto([w+b(2x+ay+bz):-(x+ay+bz):y:z],[u:v]),\]
and $\SlO^{L,L'}\simeq\{(\alpha,\beta)\in T^{L,L'}\mid
\alpha\beta=1\}$, whose $\k$-rational points are given by
    \begin{enumerate}
    \item\label{min-cb-sing1:1} either $\SlO^{L,L}(\k)\simeq\{\alpha\in L^*\mid \alpha\alpha^g=1\}$,
    \item\label{min-cb-sing1:2} or $\SlO^{L,L'}(\k)\simeq\k^*$
    if $L,L'$ are not $\k$-isomorphic.
    \end{enumerate}
\end{lem}
\begin{proof}

The indicated map is the composition of the two commuting
involutions $\alpha,\beta$ from Lemma~\ref{prop:DP2}(\ref{DP2:2})
and Lemma~\ref{prop:DP9}(\ref{DP9:4}). In particular, it is an
involution (it is not the identity by
Remark~\ref{rmk:minpoly_quadext}) that induces a reflection on the
hexagon exchanging the geometric components of the singular fibres.

Let us compute the image of $\SlO^{L,L'}$ in $T^{L,L'}$. Since this
means computing the $\bk$-points of these groups, it suffices to
assume that $L$ and $L'$ are $\k$-isomorphic. We consider $\Ql^L$ as
$\k$-structure on $\p^1_L\times\p^1_L$. By
Lemma~\ref{lem:Q-deg2}(\ref{Q-deg2:1}), we can assume that
$p_1=([0:1],[1:0])$, $p_2=([1:0],[0:1])$. Then $\SlO^{L,L}$ is
conjugate to a subgroup of the group of diagonal maps
$\Aut(\Ql^L,p_1,p_2)$. In these coordinates, the fibration
$\pi\colon\Sl^{L,L}\to\p^1$ is mapped by $\eta$ to the pencil of
curves given by $c u_1v_1-d u_0v_0=0$, $[c:d]\in\p^1$. A diagonal
element $(\alpha,\beta)\in \Aut(\Ql^L,p_1,p_2)$ preserves each fibre
if and only if $\alpha\beta=1$. It follows that
$\SlO^{L,L}=\{(\alpha,\beta)\in T^{L,L}\mid \alpha\beta=1\}$.

(\ref{min-cb-sing1:1}) The $\k$-rational points $\SlO^{L,L}(\k)$ form
the subgroup of elements in $\SlO^{L,L'}(\bk)$ that are fixed by the
$\Gal(L/\k)$-action, see Lemma~\ref{lem:autQ}. The generator
$g\in\Gal(L/\k)$ acts by $(\alpha,\beta)^g=(\beta^g,\alpha^g)$, see
Lemma~\ref{lem:autQ}. It follows that
$\SlO^{L,L}=\{(\alpha,\beta)\in T^{L,L}\mid \alpha\beta=1\}$.

(\ref{min-cb-sing1:2}) Suppose that $L,L'$ are not $\k$-isomorphic.
Let $K:=LL'$. Then $\Gal(K/\k)\simeq\Gal(L/\k)\times\Gal(L'/\k)$.
Lemma~\ref{lem:Q-deg22}(\ref{Q-deg22:3}) tells us that we can assume
that $p_1=([b_1:1],[b_1:1]),p_2=([b_2:1],[b_2:1])$. We now compute
the form of the elements in $\SlO^{L,L'}(K)$: the element
\[
\gamma:=\left(\left(\begin{matrix} b_2&b_1\\ 1&1\end{matrix}\right),
\left(\begin{matrix}
b_1&b_2\\1&1\end{matrix}\right)\right)\in\PGL_2(K)\times\PGL_2(K).
\]
induces a change of coordinates $\gamma\colon\Ql^L_K\to\Q_L^K$
sending $([0:1],[1:0]),([1:0],[0:1])$ onto $p_1,p_2$, respectively.
Then $\SlO^{L,L'}(K)\subset\PGL_2(K)^2$ is the subgroup of elements of the
form
\begin{equation}\tag{AB}\label{eq:AB}
(A,B):=\gamma\circ(\alpha,\beta)\circ\gamma^{-1}=\left(\left(\begin{matrix}
b_2\alpha-b_1 & b_1b_2(1-\alpha)\\ \alpha-1 & b_2-\alpha
b_1\end{matrix}\right), \left(\begin{matrix} b_1\beta-b_2 &
b_1b_2(1-\beta)\\ \beta-1 & b_1-b_2\beta\end{matrix}\right)\right).
\end{equation}
The group $\SlO^{L,L'}(\k)$ is the $\Gal(K/\k)$-invariant subgroup
of $\SlO^{L,L'}(K)$. If $g$ is the generator of $\Gal(L/\k)$, and
$g'$ is the one of $\Gal(L'/g)$, then
\[(A,B)^g=(B^g,A^g),\quad (A,B)^g=(A^{g'},B^{g'})
\]
It follows that
\[
\SlO^{L,L'}(\k)=\{(A,B)\in \PGL_2(L')^2\mid \text{$(A,B)$ of the
form (\ref{eq:AB})},\ \alpha\beta^g=1=\alpha\beta\}
\]
We obtain that $\beta\in\k^*$, and hence that
$\SlO^{L,L'}(\k)\simeq\k^*$.
\end{proof}

\begin{lem}\label{lem:min-cb-sing2}
Keep the notation of Remark~\ref{rmk:min-cb-sing} and let $g$ be the
generator of $\Gal(L/\k)$. Then the action of $\Aut(\Sl^{L,L'},\pi)$
on $\p^1$ induces the exact sequences
\begin{align*}
1\rightarrow\Aut(\Sl^{L,L'}/\pi)\to&\Aut(\Sl^{L,L'},\pi)\to\Aut(\p^1,\{t,s\})\simeq T_1\rtimes\Z/2\rightarrow1\\
1\rightarrow\Aut_\k(\Sl^{L,L'}/\pi)\to&\Aut_\k(\Sl^{L,L'},\pi)\to
D_{\k}^{L,L'}\rtimes\Z/2\rightarrow 1
\end{align*}
where $T_1$ is the $1$-dimensional split torus, $\Z/2$ is generated
by the image of
\[ ([w:x:y:z],[u:v])\mapsto([w+b(2x+ay+bz):-(x+bz):-y:z],[u+abv:-v])\]
and $D^{L,L'}_{\k}\subseteq T_1(\k)$ is the subgroup
    \begin{enumerate}
    \item\label{min-cb-sing2:1} $D^{L,L}_{\k}=\{\delta\in T_1(\k)\mid\delta=\lambda\lambda^g, \lambda\in L^*\}$, where $g$ is the generator of $\Gal(L/\k)$,
    \item\label{min-cb-sing2:2} $D^{L,L'}_{\k}\simeq\{\lambda\lambda^{gg'}\in F\mid \lambda\in K, \lambda\lambda^{g'}=1\}$  if $L$ and $L'$ are not $\k$-isomorphic, where $\k\subset F\subset LL'$ is the intermediate extension such that $\Gal(F/\k)\simeq\langle gg'\rangle\subset\Gal(L/\k)\times\Gal(L'/\k)$, where $g,g'$ are the generators of $\Gal(L/\k),\Gal(L'/\k)$, respectively.
    \end{enumerate}
\end{lem}
\begin{proof}
The birational contraction $\eta\colon\Sl^{L,L'}\to\Ql^L$ induces a
rational map $\hat{\pi}\colon\Ql^L\rat\p^1$ such that
$\hat{\pi}\circ\eta=\pi$. We define
\[
\Aut(\Ql^L,\hat{\pi})=\{\alpha\in\Aut(\Ql^L)\mid\exists
f\in\Aut(\p^1)\ \text{such that}\
\hat{\pi}\circ\alpha=f\circ\hat{\pi}\}
\]
Then $\Aut(\Ql^L,\hat{\pi})=\eta\Aut(\Sl^{L,L'},\pi)\eta^{-1}$. Let
us compute $\Aut_{\bk}(\Ql^L_{\bk},\hat{\pi})$. For this, we can
assume that $p_1=([0:1],[1:0])$, $p_2=([1:0],[0:1])$ (in the
notation of Remark~\ref{rmk:min-cb-sing}), and the fibres of
$\hat{\pi}$ are of the form $c u_1v_1-d u_0v_0=0$, $[c:d]\in\p^1$.
It follows that
\[
\Aut_{\bk}(\p^1_{\bk}\times\p^1_{\bk})\supseteq\Aut_{\bk}(\Ql^L_{\bk},\hat{\pi})=\left\{\left(A_{\lambda},B_{\mu}\right)\mid
\lambda,\mu\in\bk^*\right\}\rtimes\langle\tau\colon(x,y)\mapsto(y,x)\rangle
\]
where
\[
(I)\
A_{\lambda}=\left(\begin{matrix}1&0\\0&\lambda\end{matrix}\right),B_{\mu}=\left(\begin{matrix}1&0\\0&\mu\end{matrix}\right),\qquad\text{or}\qquad
(II)\
A_{\lambda}=\left(\begin{matrix}0&1\\\lambda&0\end{matrix}\right),B_{\mu}=\left(\begin{matrix}0&1\\
\mu&0\end{matrix}\right).
\]
The automorphism $(A_{\lambda},B_{\mu})$ of type (I) induces the
scaling $[c:d]\mapsto[c:\lambda\mu d]$ on $\p^1$, the one of type
(II) induces $[c:d]\mapsto[ d: \lambda\mu c]$, and $\tau$ induces
$\mathrm{id}_{\p^1}$. Hence, the image of $\Aut(\Sl^{L,L'},\pi)$ in
$\Aut(\p^1,\{t,s\})$ is $T_1\rtimes\Z/2\Z$.

Let us compute $\Aut_\k(\Ql^L,\hat{\pi})$, its image in
$\Aut_\k(\p^1,\{t,s\})$  separately for each of the two cases $L=L'$
and $L,L'$ not $\k$-isomorphic. We will use that
$\Ql_K^L\simeq\p^1_K\times\p^1_K$ for $K=LL'$, hence
$(A_\lambda,B_\mu)\in\Aut_K(\Ql_K^L,\hat{\pi})$ exactly if
$\lambda,\mu\in K$.

(\ref{min-cb-sing2:1}) Suppose that $L=L'$. Then $\tau\in
\Aut_{\k}(\Ql^L,\hat{\pi})$. An element $(A_{\lambda},B_{\mu})\in
\Aut_{L}(\Ql^L,\hat{\pi})$ is defined over $\k$ if and only
$\lambda,\mu\in L$ and $A_{\lambda}=B_{\mu}^g$, which is equivalent
to $\mu=\lambda^g$. In that case, $\lambda\mu=\lambda\lambda^g$,
which is contained in $\k$. Therefore, the image
$\Aut_{\k}(\Sl^{L,L'},\pi)$ in $\Aut_{\k}(\p^1,\{t,s\})$ is
isomorphic to $D_\k^{L,L}\rtimes\Z/2$.

(\ref{min-cb-sing2:2}) Suppose that $L$ and $L'$ are not
$\k$-isomorphic. Let $K=LL'$ and
$\Gal(K/\k)\simeq\Gal(L/\k)\times\Gal(L'/\k)=\langle
g\rangle\times\langle g'\rangle$. Let us compute
$\Aut_K(\Ql^L_K,\hat{\pi})$. Observe that we have
$p_i=([b_i:1],[b_i:1])$ for $i=1,2$ and that we can no longer assume
that they are equal to $([1:0],[0:1]),([0:1],[1:0])$. However, the
coordinate change given by
\[
\gamma:=\left(\left(\begin{matrix} b_2&b_1\\ 1&1\end{matrix}\right),
\left(\begin{matrix}
b_1&b_2\\1&1\end{matrix}\right)\right)\in\PGL_2(K)\times\PGL_2(K)
\]
sends $([1:0],[0:1]),([0:1],[1:0])$ onto $p_1,p_2$, respectively.
One can compute that the Galois action on $\Ql^{L}$ (see the proof
of Lemma~\ref{prop:Q}) induced by $\gamma$, namely
$G'=\gamma^{-1}\circ\Gal(K/\k)\circ\gamma$, is given by
\begin{align*}
    ([u_0:u_1],[v_0:v_1])&\mapsto\begin{cases}
        ([v_1^g:v_0^g],[u_1^g:u_0^g])\\
        ([u_1^{g'}:u_0^{g'}],[v_1^{g'}:v_0^{g'}]).
\end{cases}
\end{align*}
Note that $\tau$ is $G'$-invariant and so it remains to study which
$(A_\lambda,A_\mu)$ are $G'$-invariant. So $(A_\lambda,A_\mu)$ is
defined over $\k$ for $\lambda,\mu\in K$ if and only if
\begin{align*}
    (A_\lambda,A_\mu)&=(A_\lambda,A_\mu)^g=( A_{(\mu^{-1})^g}, A_{(\lambda^{-1})^g} ) \\
    (A_\lambda,A_\mu)&=(A_\lambda,A_\mu)^{g'}=(A_{(\lambda^{-1})^{g'}},A_{(\mu^{-1})^{g'}}).
\end{align*}
Hence, the elements of $\Aut_\k(\Ql^{L},\hat{\pi})$ are exactly
those of the form $\gamma\circ(A_\lambda,A_\mu)\circ\gamma^{-1}$
with $\lambda,\mu\in K$ satisfying $\lambda={(\mu^{-1})}^{g}$,
$\mu={(\lambda^{-1})}^g$ $\lambda={(\lambda^{-1})}^{g'}$,
$\mu={(\mu^{-1})}^{g'}$.

Instead of computing the image of $\Aut_\k(\Ql^{L},\hat{\pi})$ in
$\Aut_\k(\p^1)$, we compute the image of
$\gamma^{-1}\Aut_\k(\Ql^{L},\hat{\pi})\gamma$ ({\it i.e.} $(A_\lambda,A_\mu)$) on $\Aut_K(\p^1)$ with the induced Galois
action on $\p^1$, which is given by $([c:d])^g=[d^g:c^g]$ and
$([c:d])^{g'}=[d^{g'}:c^{g'}]$. Again, $(A_{\lambda},A_{\mu})$
induces $[c:d]\mapsto[ c:\lambda\mu d]$ or
$[c:d]\mapsto[d:\lambda\mu c]$, and $\tau'$ induces
$\mathrm{id}_{\p^1}$. We compute the possible $\delta=\lambda\mu$:
On one hand we find
\[\lambda\mu={(\mu^{-1})}^g({\lambda^{-1})}^g=(\mu^{g'})^g(\lambda^{g'})^g=(\lambda\mu)^{gg'},\]
implying $\delta\in F$, where $\k\subset F\subset K$ with
$\Gal(F/\k)=\langle gg'\rangle$. On the other hand, we also have
 \[\lambda\mu=\lambda(\mu^{-1})^{g'}=\lambda\lambda^{g\cdot g'}\]
Hence, $D_\k^{L,L'}$ is conjugated to $\{\lambda\lambda^{gg'}\in
F\mid \lambda\in K,\ \lambda\lambda^{g'}=1\}$.
\end{proof}

In the lemma above, if $L,L'$ are not $\k$-isomorphic, then
$D_{\k}^{L,L'}\simeq\{N_{F/\k}(\lambda)\mid \lambda\in K,\
N_{K/L}(\lambda)=1\}$, where $N_{F/\k}$ and $N_{K/L}$ are the field
norms of $F/\k$ and $K/L$, respectively.

\section{The conic fibration cases}\label{sec:CB}

In this section, we classify the rational conic fibrations
$\pi\colon X\rightarrow\p^1$ that are $\Aut(X,\pi)$-Mori fibre
spaces. Recall that $\pi$ induces a homomorphism
$\Aut(X,\pi)\rightarrow\Aut(\p^1)$ whose kernel we denote by
$\Aut(X/\pi)$ and its $\k$-points by $\Aut_\k(X/\pi)$.

Recall from Lemma~\ref{lem:min-cb} that, for any Mori fibre space
$\pi\colon X\to\p^1$ such that $X$ is rational, we have either
$X\simeq\F_n$ for some $n\geq0$ or $X\simeq\Sl^{L,L'}$ or $X$ is
isomorphic to a del Pezzo surface obtained by blowing up $\p^2$ in a
point of degree $4$. In the latter case, $\Aut(X,\pi)$ is finite by
Lemma~\ref{lem:aut5}, so we do not look at it.

\subsection{Conic fibrations obtained by blowing up a Hirzebruch surface}\label{ss:F}

We study the rational conic fibrations $\pi\colon X\rightarrow\p^1$
that are $\Aut(X,\pi)$-Mori fibre spaces and for which there is a
birational morphism $X\to \F_n$ of conic fibrations for some
$n\geq0$.

\begin{rmk}\label{rmk:cb-F}
Let $n\geq1$ and denote by $\k[z_0,z_1]_n\subset\k[z_0,z_1]$ the
vector space of homogeneous polynomials of degree $n$. In the
coordinates from Example~\ref{ex:cb}(\ref{ex:Fn}) the special
section $S_{-n}\subset\F_n$ is given by $y_0=0$. We denote by
$S_n\subset\F_n$ the section given by $y_1=0$. Since $S_n\cdot
S_{-n}=0$, we have $S_n\sim S_{-n}+nf$ and $S_n^2=n$, where $f$ is
the class of a fibre. The automorphism group of $\F_n$ is
\[\Aut(\F_n)=\Aut(\F_n,\pi_n)\simeq V_{n+1}\rtimes\GL_2/\mu_n,\quad \Aut_\k(\F_n)\simeq \k[z_0,z_1]_n\rtimes\GL_2(\k)/\mu_n(\k),\]
where $V_{n+1}$ is the canonical $\k$-structure on $\bk[z_0,z_1]_n$ and
$\mu_n=\{\lambda\cdot\id\in\GL_2\mid \lambda^n=1\}$. The group
$\Aut_\k(\F_n)$ acts on $\F_n$ by
$$
[y_0:y_1;z_0:z_1]\mapsto[y_0:P(z_0,z_1)y_0+y_1; az_0+bz_1:cz_0+dz_1],
$$
and it has two orbits on $\F_n$, namely $S_{-n}$ and $\F_n\setminus
S_{-n}$.
\end{rmk}

\begin{lem}\label{lem:cb-F}
Let $n\geq0$ and $\eta\colon X\rightarrow\F_n$ be a birational
morphism of conic fibrations that is not an isomorphism, and suppose
that $\Aut_{\bk}(X,\pi)$ contains an element permuting the
components of at least one singular geometric fibre. Let
$G_{\bk}\subset\Aut_{\bk}(X/\pi)$ be the subgroup of elements acting
trivially on $\NS(X_{\bk})$.
\begin{enumerate}
\item\label{cb-F:1a} If $G_{\bk}$ is non-trivial, there exists $N\geq1$ and a birational morphism $X\rightarrow\F_N$ of conic fibrations blowing up
$r\geq1$ points $p_1,\dots,p_{r}$ contained in $S_N$ such that
$\sum_{i=1}^r\deg(p_i)=2N$.
\item If $G_{\bk}=\{1\}$, then $\Aut_{\bk}(X/\pi)\simeq(\Z/2)^r$ for $r\in\{0,1,2\}$.
\end{enumerate}
\end{lem}
\begin{proof}
The claim is proven in \cite[Lemme 4.3.5]{Blanc_alg_subgroups} over
$\C$ and its proof can be repeated word by word over any
algebraically closed field. Over a perfect field $\k$ it suffices to
show that curves contracted by the birational morphism $\nu\colon
X_{\bk}\to (\F_N)_{\bk}$ in (\ref{cb-F:1a}) are already defined over
$\k$. Since $N\geq1$, the surface $X_{\bk}$ contains exactly two
sections of negative self-intersection, namely the strict transforms
$\tilde{S}_{-N}$ and $\tilde{S}_N$ of $S_{-N}$ and $S_N$,
respectively, and $\tilde{S}_{-N}^2=\tilde{S}_N^2=-N$, and every
singular geometric fibre has two components, each intersecting
either $\tilde{S}_{-N}$ or $\tilde{S}_N$. We now show that
$\tilde{S}_{-N}$ and $\tilde{S}_N$ are both defined over $\k$, which
will then imply that the curves contracted by $\eta$ are defined
over $\k$ and we are finished. The birational morphism $\eta\colon
X\to\F_n$ contracts exactly one component in each singular fibre.
This implies that the strict transform $\tilde{S}_{-n}$ of
$S_{-n}\subset\F_n$ has self-intersection $\leq-n$. If $n\geq1$,
then $\tilde{S}_{-n}$ is one of $\tilde{S}_N$ or $\tilde{S}_{-N}$
and hence both $\tilde{S}_N$ or $\tilde{S}_{-N}$ are defined over
$\k$. If $n=0$, then $\eta(\tilde{S}_{-N})$ and $\eta(\tilde{S}_N)$
are sections in $\F_0$ of ruling induced by $\eta$. If they are
permuted by an element of $\Gal(\bk/\k)$, each fibre contains two
points blown-up by $\eta$, which contradicts $X\to\p^1$ being a
conic fibration. It follows that $\eta(\tilde{S}_{-N})$ and
$\eta(\tilde{S}_N)$ are both defined over $\k$ and hence
$\tilde{S}_{-N},\tilde{S}_N$ are defined over $\k$ as well.
\end{proof}

Let us construct a special birational involution of $\F_n$,
$n\geq1$.

\begin{ex}\label{ex:gen-F}
Let $n\geq1$. Let $p_1,\dots,p_r\in S_n\subset\F_n$ be points such
that their geometric components are in pairwise distinct geometric
fibres and $\sum_{i=1}^r\deg(p_i)=2n$, and assume that
$\pi_n(p_i)\neq[0:1],[1:0]$ for $i=1,\dots,r$. Let
$P_i\in\k[z_0,z_1]_{\deg(p_i)}$ be the polynomial defining
$\pi(p_i)\in\p^1$ and define $P:=P_1\cdots P_r\in\k[z_0,z_1]_{2n}$.
Then the map
\[
\varphi\colon\F_n\rat\F_n,\ (y_1,z_1)\dashmapsto(\nicefrac{P(z_1)}{y_1},z_1)
\]
is an involution preserving the fibration, whose base-points are
$p_1,\dots,p_r$, that exchanges $S_n$ and $S_{-n}$ and contracts the
fibres through $p_1,\dots,p_r$.
\end{ex}

We call $\mu_n\subset T_1$ the subgroup of $n^\mathrm{th}$ roots of unity of
the $1$-dimensional standard torus $T_1$.

\begin{lem}\label{prop:cb-F}
Let $n\geq1$ and let $\eta\colon X\rightarrow\F_n$ be a birational
morphism blowing up points $p_1,\dots,p_r\in S_n$ whose geometric
components are on pairwise distinct geometric fibres and such that
$\sum_{i=1}^r\deg(p_i)=2n$. Then $\pi:=\pi_n\eta\colon X\to\p^1$ is
a conic fibration that has exactly two $(-n)$-sections and the
following properties hold.
\begin{enumerate}
\item\label{cb-F:1} There are split exact sequences
\begin{align*}1\rightarrow\Aut(X/\pi)\to&\Aut(X,\pi)\to\Aut(\p^1,\Delta)\rightarrow 1\\
1\rightarrow\Aut_\k(X/\pi)\to&\Aut_\k(X,\pi)\to\Aut_\k(\p^1,\Delta)\rightarrow
1
\end{align*}
where $\Delta\subset\p^1$ is the image of the singular fibres of
$X/\p^1$.
\item\label{cb-F:2} The action of $\Aut(X/\pi)$ on the two $(-n)$-sections induces split exact sequences
\begin{align*}
1\rightarrow H\to&\Aut(X/\pi)\to\Z/2\rightarrow1,\\
1\rightarrow H(\k)\to&\Aut_\k(X/\pi)\to\Z/2\rightarrow1
\end{align*}
where $\eta H\eta^{-1}=\Aut(\F_n/\pi_n, S_n)\simeq T_1/\mu_n$ and
$\eta H(\k)\eta^{-1}\simeq\k^*/\mu_n(\k)$, and
$\Z/2=\langle\eta^{-1}\varphi\eta\rangle$ with
$\varphi\colon\F_n\rat\F_{n}$ the involution from
Example~\ref{ex:gen-F}.
\item\label{cb-F:4} Any element of $\Aut_\k(X/\pi)\setminus H(\k)$ is an involution fixing an irreducible double cover of $\p^1$ branched over $\Delta$ not intersecting $S_{-n}$.
\item\label{cb-F:5} $\pi\colon X\to\p^1$ is an $\Aut(X,\pi)$-Mori fibre space and an $\Aut_\k(X,\pi)$-Mori fibre space.
\end{enumerate}
\end{lem}
\begin{proof}
We denote by $\tilde{S}_n$ and $\tilde{S}_{-n}$ the strict
transforms of the sections $S_n$ and $S_{-n}$ of $\F_n$ in $X$,
which satisfy $\tilde{S}_n^2=\tilde{S}_{-n}^2=-n$ and which are the
only (geometric) sections of negative self-intersection. The
anti-canonical divisor of $X$ is $\pi$-ample because the geometric
components of the $p_i$ are on pairwise distinct geometric fibres,
thus $\pi\colon X\to\p^1$ is a conic fibration with $r$ singular
fibres, each of whose geometric components intersects exactly one of
the sections $\tilde{S}_n$ and $\tilde{S}_{-n}$.

(\ref{cb-F:1}) For any element $\alpha\in\Aut(\p^1,\Delta)$ there
exists $\tilde{\alpha} \in\Aut(\F_n)$ preserving
$\{p_1,\dots,p_r\}$, and we have
$\eta^{-1}\tilde{\alpha}\eta\in\Aut(X,\pi)$. The same argument holds
for the $\k$-points of these groups.

(\ref{cb-F:2}) Up to an element of $\Aut_\k(\F_n)$, we can assume
that $\pi_n(p_i)\neq[1:0],[0:1]$ for $i=1,\dots,r$. Then the
birational involution $\varphi\colon\F_n\rat\F_n$ from
Example~\ref{ex:gen-F} lifts to an element of $\Aut_\k(X/\pi)$ and
exchanges $\tilde{S}_n$ and $\tilde{S}_{-n}$. It follows that the
action of $\Aut(X/\pi)$ on $\{\tilde{S}_n,\tilde{S}_{-n}\}$  induces
split exact sequences
\[
1\rightarrow H\to\Aut(X/\pi)\to\Z/2\rightarrow1,\quad\text{and}\quad 1\rightarrow
H(\k)\to\Aut_\k(X/\pi)\to\Z/2\rightarrow1.
\]
Any element of $H$ fixes $\tilde{S}_n$ and $\tilde{S}_{-n}$
pointwise, so $\eta H\eta^{-1}$ and $\eta H(\k)\eta^{-1}$ are the
subgroups of $\Aut(\F_n/\pi_n)\simeq V_{n+1}\rtimes T_1/\mu_n$ and
$\Aut_\k(\F_n/\pi_n)\simeq \k[z_0,z_1]_n\rtimes \k^*/\mu_n(\k)$,
respectively, fixing $S_n$ pointwise. It follows that $\eta
H\eta^{-1}=T_1/\mu_n$ and $\eta H(\k)\eta^{-1}=\k^*/\mu_n(\k)$.

(\ref{cb-F:5}) The fact that the element
$\eta^{-1}\varphi\eta\in\Aut_\k(X/\pi)$ exchanges the components of
every singular geometric fibre implies that
$\rk\,\NS(X)^{\Aut_\k(X,\pi)}=1$.  It follows that $X/\p^1$ is an
$\Aut_{\k}(X,\pi)$-Mori fibre space and in particular an
$\Aut(X,\pi)$-Mori fibre space.

(\ref{cb-F:4}) For any $\lambda\in \k^*$ the map
\[(\lambda,\varphi)\colon (y_1,z_1)\dashmapsto(\nicefrac{\lambda^nP(z_1)}{y_1}, z_1)
\]
is a birational involution of $\F_n$ and fixes the curve
$y_1^2-\lambda^nP(z_0,z_1)y_0^2=0$, which is a double cover of
$\p^1$ branched over $\Delta$ and does not intersect the section
$S_{-n}$.
\end{proof}

\begin{lem}\label{lem:cb-F-aut}
Let $n\geq1$ and $\eta\colon X\rightarrow\F_n$ be a birational
morphism blowing up points $p_1,\dots,p_r\in S_n$ whose geometric
components are on pairwise distinct geometric fibres and such that
$\sum_{i=1}^r\deg(p_i)=2n$. Let $\pi=\pi_n\eta\colon X\to\p^1$ be
the induced conic fibration on $X$.
\begin{enumerate}
\item\label{cb-F-aut:1} If $n=1$, then $X$ is a del Pezzo surface of degree $6$ as in \ref{fig:action-dP6}(\ref{dp6:1}) or \ref{fig:action-dP6}(\ref{dp6:3}) and $\Aut(X,\pi)\subsetneq\Aut(X)$. Moreover, $\Aut_\k(X,\pi)\subsetneq\Aut_\k(X)$ if $ X$ is as in \ref{fig:action-dP6}(\ref{dp6:1}) and $\Aut_\k(X,\pi)=\Aut_\k(X)$ if $X$ is as in \ref{fig:action-dP6}(\ref{dp6:3}).
\item\label{cb-F-aut:2} If $n\geq2$, then $\Aut(X,\pi)=\Aut(X)$.
\end{enumerate}
\end{lem}
\begin{proof}
(\ref{cb-F-aut:1}) For $n=1$, the conic fibration $X/\p^1$ has two
$(-1)$-sections and $X$ is a del Pezzo surface of degree $6$ as in
Figure~\ref{fig:action-dP6}(\ref{dp6:1}) or
Figure~\ref{fig:action-dP6}(\ref{dp6:3}).
Lemma~\ref{prop:DP1}(\ref{DP1:2}) applied to $X_{\bk}$ implies that
$\Aut(X)$ contains an element inducing a rotation of order $6$ on
the hexagon of $X$, which is not contained in $\Aut(X,\pi)$. The same argument implies that $\Aut_\k(X,\pi)\subsetneq \Aut_\k(X)$ if $X$ is a del Pezzo surface of degree $6$ as in \ref{fig:action-dP6}(\ref{dp6:1}).
However, in the case of Figure~\ref{fig:action-dP6}(\ref{dp6:3}),
any element of $\Aut_\k(X)$ preserves the fibration by
Lemma~\ref{prop:DP3}(\ref{DP3:3}).

(\ref{cb-F-aut:2}) If $n\geq2$, $X$ contains exactly two
$(-n)$-sections $\tilde{S}_n$ and $\tilde{S}_{-n}$, which are the
strict transforms of $S_n$ and $S_{-n}$. Thus the class
$\tilde{S}_n+\tilde{S}_{-n}$ in $\NS(X_{\bk})$ is
$\Aut_{\bk}(X)$-invariant, hence
$K_X+(\tilde{S}_n+\tilde{S}_{-n})=-2f$ is $\Aut_{\bk}(X)$-invariant
as well. It follows that $\Aut(X)=\Aut(X,\pi)$.
\end{proof}

If two conic fibrations as in Lemma~\ref{prop:cb-F} are isomorphic,
they both have a birational morphism to the same Hirzebruch surface
$\F_n$.

\begin{lem}\label{lem:min-cb-singF}
For any fixed $n\geq1$, two conic fibrations as in
Lemma~\ref{prop:cb-F} are isomorphic if and only if the points on
$\p^1$ are the same, up to an element of $\Aut_\k(\p^1)$.
\end{lem}
\begin{proof}

Any element of $\Aut_\k(\p^1)$ lifts to an element of
$\Aut_\k(\F_n)$, so two such conic fibrations are isomorphic, if and
only if the points on the section $S_n$ are the same, up to an
element of $\Aut_\k(\F_n)$. This means that their images on $\p^1$
are the same, up to an element of $\Aut_\k(\p^1)$.
\end{proof}

\subsection{Conic fibrations obtained by blowing up a del Pezzo surface}\label{ss:cbZ}
Let $L=\k(a_1)$ and $L'=\k(b_1)$ be quadratic extensions of $\k$. In
this section, we consider rational conic fibrations $\pi\colon
X\to\p^1$ for which there is a birational morphism $\eta\colon
X/\p^1\to \Sl^{L,L'}/\p^1$ of conic fibrations, where
$\pi_{\Sl^{L,L'}}\colon \Sl^{L,L'}\to\p^1$ is the Mori fibre space
from Example~\ref{ex:cb}(\ref{ex:Sl}). We have described the
fibration $S^{L,L'}\to\p^1$ in Section~\ref{ss:fibrationDP}.

Recall from Lemma~\ref{prop:DP2}(\ref{DP2:4}) and
Lemma~\ref{prop:DP9}(\ref{DP9:1}) that there is a birational
morphism $\nu\colon\Sl^{L,L'}\to \Ql^L$ contracting a curve $E$ onto
a point $p'$ of degree $2$ with splitting field $L'$.

\begin{rmk}\label{rmk:evendeg}
Let $p\in E\subset\Sl^{L,L'}$ be a point whose geometric components
are in distinct smooth geometric fibres of $\Sl^{L,L'}/\p^1$. Any
element of $\Gal(\bk/\k)$ exchanges or preserves the geometric
components of the point $\eta(E)$ and hence of the curve $E$, and
this implies that $\deg(p)$ is even and each geometric component of
$E$ contains $\frac{\deg(p)}{2}$ geometric components of $p$.
\end{rmk}

We now show an analogue of Lemma~\ref{lem:cb-F}, that we prove
similarly to \cite[Lemme 4.3.5]{Blanc_alg_subgroups}.

\begin{lem}\label{lem:cb-Z}
Let $\eta\colon X\to \Sl^{L,L'}$ be a birational morphism of conic
fibrations that is not an isomorphism, and suppose that
$\Aut_{\bk}(X,\pi)$ contains an element exchanging the components of
at least one singular geometric fibre. Let
$G_{\bk}\subset\Aut_{\bk}(X/\pi)$ be the subgroup acting trivially
on $\NS(X_{\bk})$.
\begin{enumerate}
\item\label{cb-Z:1} If $G_{\bk}$ is non-trivial, then $\eta$ is the blow-up of $r\geq1$ points contained in $E\subset \Sl^{L,L'}$ whose geometric components are on pairwise distinct smooth geometric fibres, and each geometric component of $E$ contains half of the geometric components of each point.
\item\label{cb-Z:2} If $G_{\bk}=\{1\}$, then $\Aut_\k(X/\pi)\simeq(\Z/2)^r$ for $r\in\{0,1,2\}$.
\end{enumerate}
\end{lem}
\begin{proof}
(\ref{cb-Z:1}) Suppose that $G_{\bk}$ is nontrivial. It preserves
the geometric components of the singular fibres, so $\eta$ is
$G_{\bk}$-equivariant and $R:=\eta
G_{\bk}\eta^{-1}\subset\Aut_{\bk}(\Sl_{\bk}^{L,L'}/\pi_{\Sl^{L,L'}})$.
The group $R$ fixes the geometric components of $E$ pointwise. Since
$R\subset\PGL_2(\bk(x))$ and since it is non-trivial, it fixes no
other sections of $\Sl_{\bk}^{L,L'}/\p^1$. So, $G_{\bk}$ fixes the
geometric components of the strict transform $\tilde{E}\subset X$ of
$E$ and no other sections of $X_{\bk}/\p^1_{\bk}$. Moreover,
$\Aut_{\bk}(X,\pi)$ contains an element exchanging the components of
at least one singular geometric fibre, so it follows that each
geometric component of $\tilde{E}$ intersects exactly one component
of each geometric singular fibre. In particular, the points blown-up
by $\eta$ are contained in $E$. The hypothesis that $-K_X$ is
$\pi$-ample implies that the geometric components of the blown-up
points are on distinct geometric components of smooth fibres. The
remaining claim follows from Remark~\ref{rmk:evendeg}.

(\ref{cb-Z:2}) If $G_{\bk}$ is trivial, then every non-trivial
element of $\Aut_{\bk}(X/\pi)$ is an involution and the claim
follows from the fact that $\Aut_{\bk}(X/\pi)\subset\PGL_2(\bk(x))$.
\end{proof}

\begin{ex}\label{ex:gen-cb1}
Let us construct a special birational involution of $\varphi_{L,L'}$
of $\Sl^{L,L'}$ that preserves the fibration $\Sl^{L,L'}\to\p^1$ and
induces the identity on $\p^1$.

Let $E_1,E_2$ be the geometric components of $E$. If $g'$ is the
generator of $\Gal(L'/\k)$, then $E_2^g=E_1$. Let $p_1,\dots,p_r\in
E\subset \Sl^{L,L'}$ be points whose geometric components are on
pairwise distinct smooth geometric fibres. We now construct an
involution $\varphi$ of $\Sl^{L,L'}$ whose base-points are
$p_1,\dots,p_r$ and which exchanges $E_1$ and $E_2$. For $i=1,2$, let
$P_i\in L[x,y]$ be homogeneous polynomials defining the set of
components of the $p_1,\dots,p_r$ contained in $E_i$. Consider a
birational morphism $\Sl^{L,L'}\to\Ql^L$ that contracts $E$, and
consider the model of $\Ql^L$ that is a $\k$-structure on
$\p^1_L\times\p^1_L$.
\begin{itemize}
\item If $L$ and $L'$ are $\k$-isomorphic, we can assume that the images of $E_1$ and $E_2$ are respectively $([1:0],[0:1])$ and $([0:1],[1:0])$, by Lemma~\ref{lem:Q-deg2}(\ref{Q-deg2:1}). We define
\begin{multline*}
\tilde{\varphi}_{L,L}\colon([u_0:u_1],[v_0:v_1])\mapsto\\
([v_0P_1(u_0v_0,u_1v_1):v_1P_2(u_0v_0,u_1v_1)],[u_0P_2(u_0v_0,u_1v_1):u_1P_1(u_0v_0,u_1v_1)]).
\end{multline*}
\item If $L$ and $L'$ are not $\k$-isomorphic, we write $L'=\k(b_1)$. By Lemma~\ref{lem:Q-deg22}(\ref{Q-deg22:3}), we can assume that the images of $E_1,E_2$ are $([b_1:1],[b_1:1]),([b_2:1],[b_2:1])$. To compute $\varphi_{L,L'}$, we simply conjugate $\varphi_{L,L}$ over $\bk$ with
\[
\gamma:=\left(\left(\begin{matrix} b_2&b_1\\ 1&1\end{matrix}\right),
\left(\begin{matrix}
b_1&b_2\\1&1\end{matrix}\right)\right)\in\PGL_2(\bk)\times\PGL_2(\bk)
\]
This yields the following form of $\varphi_{L,L'}$
\[\tilde{\varphi}_{L,L'}\colon([u_0:u_1],[v_0:v_1])\mapsto([v_0U+v_1V:v_0W-v_1U],[u_0U+u_1V:u_0W-u_1U])\]
where
\[U:=b_2P_1(t,s)-b_1P_2(t,s),\ V:=b_1^2P_2(t,s)-b_2^2P_1(t,s),\ W:=P_1(t,s)-P_2(t,s)\]
with
\[t:=(u_0-b_1u_1)(v_0-b_2v_1),\ s:=(u_0-b_2v_1)(v_0-b_1v_1).\]
\end{itemize}
In both cases, $\tilde{\varphi}_{L,L'}$ commutes with $\Gal(L/\k)$
and $\Gal(L'/\k)$ and it is an involution. Moreover, it preserves
the image of the fibration $\Sl^{L,L'}\to\p^1$ in $\Ql^L$ and
induces the identity map on $\p^1$. The base-locus of
$\tilde{\varphi}_{L,L'}$ in $\Ql^L$ is the image of $E$, and $\tilde{\varphi}_{L,L'}$
contracts the image of the fibres of $\Sl^{L,L'}\to\p^1$ given by
$P_1P_2=0$. It follows that $\tilde{\varphi}_{L,L'}$ lifts to a birational involution
$\varphi_{L,L'}$ not defined in $p_1,\dots,p_r$.
\end{ex}

\begin{lem}\label{prop:cb1}
Let $\eta\colon X\rightarrow \Sl^{L,L'}$ be the blow-up up of points
$p_1,\dots,p_r\in E$, $r\geq1$, whose geometric components are on
pairwise distinct smooth geometric fibres. Then
$\pi:=\pi_\Sl\eta\colon X\to\p^1$ is a conic fibration and
$\deg(p_i)$ is even and each geometric component of $E$ contains
$\frac{\deg(p_i)}{2}$ geometric components for $i=1,\dots,r$.
Moreover, the following hold.
\begin{enumerate}
    \item\label{cb1:1} The action of $\Aut(X,\pi)$ on $\p^1$ induces the exact sequence
    \begin{align*}1\rightarrow\Aut(X/\pi)\to&\Aut(X,\pi)\to\Aut(\p^1,\Delta)\rightarrow 1\\
1\rightarrow\Aut_\k(X/\pi)\to&\Aut_\k(X,\pi)\to
(D^{L,L'}_{\k}\rtimes\Z/2)\cap\Aut_\k(\p^1,\Delta)\rightarrow 1
    \end{align*}
where $D^{L,L'}_{\k}\rtimes\Z/2$ is the image of
$\Aut_\k(\Sl^{L,L'},\pi)$ in $\Aut_\k(\p^1)$, see
Lemma~\ref{lem:min-cb-sing2}, and $\Delta\subset\p^1$ is the image
of the singular fibres of $X$.
    \item\label{cb1:2}
    The $\Aut(X/\pi)$-action on the components of the strict transform of $E$ induces the split exact sequences
    \begin{align*}
    1\rightarrow H\to&\Aut(X/\pi)\to\Z/2\rightarrow 1,\\
    1\rightarrow H(\k)\to&\Aut_\k(X/\pi)\to\Z/2\rightarrow 1
    \end{align*}
    with $\eta H\eta^{-1}=\SlO^{L,L'}$ from Lemma~\ref{lem:min-cb-sing1} and $\Z/2$ is generated by the involution $\varphi_{L,L'}\colon\Sl^{L,L'}\rat\Sl^{L,L'}$ from Example~\ref{ex:gen-cb1}.
        \item\label{cb1:3}
    Any element of $\Aut_\k(X/\pi)\setminus H(\k)$ is an involution fixing an irreducible double cover of $\p^1$ branched over $\Delta$.
    \item\label{cb1:4} $\pi\colon X\to\p^1$ is an $\Aut(X,\pi)$-Mori fibre space and an $\Aut_\k(X,\pi)$-Mori fibre space.
\end{enumerate}
\end{lem}
\begin{proof}
The first claim follows from Remark~\ref{rmk:evendeg} and the
sequences in (\ref{cb1:1}) are exact by
Lemma~\ref{lem:min-cb-sing1}.

(\ref{cb1:2}) Consider the involution $\varphi_{L,L'}\colon
\Sl^{L,L'}\rat \Sl^{L,L'}$ from Example~\ref{ex:gen-cb1} whose
base-points are $p_1,\dots,p_r$ and that exchanges the geometric
components of $E$. Then
$\hat{\varphi}_{L,L'}:=\eta^{-1}\varphi_{L,L'}\eta$ is contained in
$\Aut_\k(X/\pi)$ and exchanges the geometric components of the
strict transform $\tilde{E}$ of $E$. In particular, the
$\Aut(X/\pi)$-action on the set of geometric components of
$\tilde{E}$ induces split exact sequences as claimed. The groups
$H$ and $H(\k)$ are respectively conjugate by $\eta$ to the subgroups
of $\Aut(\Sl^{L,L'}/\pi_{\Sl})$ and $\Aut_\k(\Sl^{L,L'}/\pi_{\Sl})$
preserving the geometric components of $E$, which are $\SlO^{L,L'}$
and $\SlO^{L,L'}(\k)$ by Lemma~\ref{lem:min-cb-sing1}.

(\ref{cb1:3}) It is enough to show that this is already the case for any
element in $\Aut_{\bk}(X/\pi)\setminus H(\bk)$. Indeed, we have
$\Aut_{\bk}(X/\pi)\simeq H(\bk)\rtimes\Z/2$, and any element of
$\Aut_{\bk}(X/\pi)\setminus H(\bk)$ is of the form
$(\eta^{-1}\alpha\eta,\hat{\varphi}_{L,L'})$, where
$\alpha:=(a,a^{-1})\in\SlO^{L,L'}(\bk)$. Using
Example~\ref{ex:gen-cb1}, we compute that
$(\eta^{-1}\alpha\eta,\hat{\varphi}_{L,L'})$ is an involution. Its
fixed $\bk$-curve in $\Ql^L_{\bk}$ is given by
 \[
a u_0v_1P_2(u_0v_0,u_1v_1)-u_1v_0P_1(u_0v_0,u_1v_1)=0
 \]
 which lifts to the desired curve on $X_{\bk}$.

(\ref{cb1:4}) The involution $\hat{\varphi}$ exchanges the geometric
components of all singular fibres and hence $X\to\p^1$ is a
$\Aut(X,\p^1)$-Mori fibre space and an $\Aut_\k(X,\p^1)$-Mori fibre
space.
\end{proof}

\begin{lem}\label{lem:cb1-aut}
Let $\eta\colon X\rightarrow \Sl^{L,L'}$ be the blow-up up of points
$p_1,\dots,p_r\in E$, $r\geq1$, whose geometric components are on
pairwise distinct smooth geometric fibres. Then
$\Aut(X,\pi)=\Aut(X)$.
\end{lem}
\begin{proof}
By Remark~\ref{rmk:evendeg}, each of the components of $E$ contains
half the geometric components of each $p_i$. It follows that
$n:=\frac{1}{2}\sum_{i=1}^r\deg(p_i)\in\Z$ and $n\geq1$. For
$i=1,\dots,r$, let $E_i$ be the exceptional divisor of $p_i$ and let
$f$ be a general fibre of $X$ and $\tilde{E}$ the strict transform
of $E$. We have $K_{\Sl}=-2f-E$ and hence
$K_X=-2f-\pi^*E+E_1+\cdots+E_r=-2f-\tilde{E}$. The curve $\tilde{E}$
is the unique curve in $X$ with self-intersection
$\tilde{E}^2=-2(1+n)\leq-4$ and hence it is $\Aut(X)$-invariant. In
particular, $K_X+\tilde{E}=-2f$ is $\Aut(X)$-invariant. It follows
that $\Aut(X)=\Aut(X,\pi)$.
\end{proof}

\begin{lem}\label{cor:min-cb-sing}
Two conic fibrations as in Lemma~\ref{prop:cb1} are isomorphic as
conic fibrations if and only if the points on $\p^1$ are the same,
up to an element of $D^{L,L'}_{\k}\rtimes\Z/2$, which is the image
of $\Aut_\k(\Sl^{L,L'},\pi)$ in $\Aut_\k(\p^1)$ (see
Lemma~\ref{lem:min-cb-sing2}).
\end{lem}
\begin{proof}
Let $X\to\Sl^{L,L'}$ and $X'\to\Sl^{L,L'}$ be such conic fibrations
obtained by blowing up $p_1,\dots,p_r\subset E$ and
$p_1',\dots,p_s'\subset E$, respectively, and suppose that they are
isomorphic as conic fibrations. Then this isomorphism sends the
singular fibres of $X$ onto the ones of $X'$, and hence descends to
an automorphism of $\p^1$ that sends the images of the $p_i$ onto
the images of the $p_i'$.

On the other hand, given an automorphism $\alpha$ of $\p^1$
contained in $D_{\k}^{L,L'}\rtimes\Z/2$, we know by
Lemma~\ref{prop:cb1} there exists an automorphism $\psi$ of $X$ that
induces $\alpha$ on $\p^1$. If $\alpha$ sends the $p_i$ onto the
$p_i'$, then either $\psi$ or $\psi\circ\varphi$ sends the $p_i$
onto the $p_i'$, where $\varphi$ is the generator of
$\Z/2\subset\Aut_\k(X/\pi)$ in Lemma~\ref{prop:cb1}(\ref{cb1:2})
exchanging the components of the singular fibres.
\end{proof}

\section{The proof of Theorem~\ref{thm:1}}

In this section, we prove Theorem~\ref{thm:1}.

\begin{lem}\label{lem:cbF-autbig}
Consider a birational morphism of conic fibrations $X\to \F_n$ for
some $n\geq0$, and suppose that $X/\p^1$ has at most two singular
geometric fibres. If there is an element of $\Aut_{\bk}(X,\pi)$ that
permutes the components of at least one singular geometric fibre,
then it has exactly two singular geometric fibres and $X$ is a del
Pezzo surface of degree $6$.
\end{lem}
\begin{proof}
Denote by $\eta\colon X\to\F_n$ the birational morphism. Let
$\tilde{S}_{-n}\subset X$ be the strict transform of the section
$S_{-n}\subset\F_n$. Then $\tilde{S}_{-n}^2\in\{-n,-n-1,-n-2\}$. Let
$\alpha\in\Aut_{\bk}(X,\pi)$ be an element that permutes the
components of at least one singular geometric fibre $f_0$. Then
$\tilde{S}:=\alpha(\tilde{S}_{-n})$ is a section of
$\eta\times\id\colon X_{\bk}\to\p^1_{\bk}$ of self-intersection
$\tilde{S}^2=\tilde{S}_{-n}^2$, and it intersects the other
component of $f_0$. It follows that $S:=\eta(\tilde{S})\subset\F_n$
is a section of self-intersection $S^2\in\{-n+2,-n+1,-n\}$,
depending on how many of the points blown up by $\eta$ are contained in $S_{-n}$.
Since $S^2\geq0$, we have $n\leq2$. If $n=2$, we have
$S^2=0$ and hence $S\sim S_{-2}+f$, which means that $S\cdot
S_{-2}=-1$, which is impossible. It follows that $n=0$ or $n=1$, and so $X$
is a del Pezzo surface of degree $6$ or $7$. In the latter case, no
element of $\Aut_{\bk}(X,\pi)$ permutes the components of the
singular fibre, hence $X$ is a del Pezzo surface of degree $6$.
\end{proof}

\begin{lem}\label{lem:thm1cb}
Let $\pi\colon X\to \p^1$ be a $\Aut(X,\pi)$-Mori fibre space with
at least three singular geometric fibres and suppose that there is a
birational morphism of conic fibrations $X\to Y$, where $Y=\F_n$ for
some $n\geq0$ or $Y=\Sl^{L,L'}$, and that $\Aut_{\bk}(X,\pi)$ is
infinite. Then the pair $(X,\Aut(X))$ is as in
Theorem~\ref{thm:1}$(\ref{1:6})$.
\end{lem}
\begin{proof}
The hypothesis that $X$ is an $\Aut(X,\pi)$-Mori fibre space implies
that $\Aut_{\bk}(X,\pi)$ contains an element permuting the
components of a singular geometric fibre. Moreover, $X/\p^1$ has at
least three singular geometric fibres, the image of the homomorphism
$\Aut_{\bk}(X,\pi)\to\Aut_{\bk}(\p^1)$ is finite and hence the
kernel $\Aut_{\bk}(X/\pi)$ is infinite.

First, suppose that $Y=\F_n$. Since $X/\p^1$ has singular fibres,
$\eta$ is not an isomorphism. Lemma~\ref{lem:cb-F} and the fact that
$\Aut_{\bk}(X/\pi)$ is infinite imply that there exists $N\geq1$ and
a birational morphism $X\to\F_N$ that blows up $p_1,\dots,p_r\in
S_N\subset\F_N$ whose geometric components are in distinct geometric
fibres and such that $\sum_{i=1}^r\deg(p_i)=2N$. Because $\pi$ has
at least three singular geometric fibres,
Lemma~\ref{lem:cb-F-aut}(\ref{cb-F-aut:1}) implies that $N\geq2$,
and now Lemma~\ref{lem:cb-F-aut}(\ref{cb-F-aut:2}) implies that
$\Aut(X,\pi)=\Aut(X)$.
Lemma~\ref{prop:cb-F}(\ref{cb-F:1}--\ref{cb-F:2}) implies that
$(X,\Aut(X))$ is as in Theorem~\ref{thm:1}(\ref{1:61}).

Now, suppose that $Y=\Sl^{L,L'}$. Since $X/\p^1$ has at least three
singular fibres, $\eta$ is not an isomorphism. Since
$\Aut_{\bk}(X/\pi)$ is infinite, Lemma~\ref{lem:cb-Z} implies that
$\eta$ blows up points $p_1,\dots,p_r\in E$ whose geometric
components are on distinct smooth geometric fibres, and
Remark~\ref{rmk:evendeg} implies that they are all of even degree
and each geometric component of $E$ contains half the geometric
components of each $p_i$. Lemma~\ref{lem:cb1-aut} implies that
$\Aut(X,\pi)=\Aut(X)$. Lemma~\ref{prop:cb1} and the description of
$D_\k^{L,L'}$ in Lemma~\ref{lem:min-cb-sing2} imply that the pair
$(X,\Aut(X))$ is as in Theorem~\ref{thm:1}(\ref{1:62}).
\end{proof}

\begin{proof}[Proof of Theorem~\ref{thm:1}]
By Proposition~\ref{prop which cases}, there is a $G$-equivariant
birational map $\p^2\rat X$ to a $G$-Mori fibre space $\pi\colon
X\to B$ that is one of the following:
\begin{itemize}
    \item $B$ is a point and $X\simeq\p^2$ or $X$ is a del Pezzo surface of degree $6$ or $8$,
    \item $B=\p^1$ and there is a (perhaps non-equivariant) birational morphism of conic fibrations $X\to Y$ with $Y=\F_n$ for some $n\geq0$ or $Y=\Sl^{L,L'}$.
    \end{itemize}
By Lemma~\ref{lem:aut}, it suffices to look at the case $G=\Aut(X)$
or $G=\Aut(X,\pi)$, respectively. The pair $(\p^2,\Aut(\p^2))$ is
the one in Theorem~\ref{thm:1}(\ref{1:1}).

If $X$ is a del Pezzo surface of degree $8$, then $X$ is isomorphic
to $\F_0$, to $\F_1$ or to $\Ql^L$ for some quadratic extension
$L/\k$ by Lemma~\ref{prop:Q}(\ref{Q:1}). However, $\F_1$ has a
unique $(-1)$-curve, which is hence $\Aut(\F_1)$-invariant and its
contraction conjugates $\Aut(\F_1)$ to a subgroup of $\Aut(\p^2)$.
It follows that $X=\Ql^L$ or $X=\F_0$, {\it i.e.} the pair
$(X,\Aut(X))$ is as in Theorem~\ref{thm:1}(\ref{1:2})--(\ref{1:3}).

If $X$ is a del Pezzo surface of degree $6$, the
$\Gal(\bk/\k)$-action on the hexagon of $X$ is one of the actions in
Figure~\ref{fig:action-dP6}(\ref{dp6:1})--(\ref{dp6:9}).
Lemma~\ref{prop:DP1}(\ref{DP1:2}--\ref{DP1:3}) applied to
$X_{\bk}$ yields that $\rk\,\NS(X_{\bk})^{\Aut_{\bk}(X_{\bk})}=1$ and
that the action of $\Aut_{\bk}(X)$ on $\NS(X_{\bk})$ induces a split
exact sequence
\[1\to(\bk^*)^2\to\Aut_{\bk}(X)\to\sym_3\times\Z/2\to1.\]

If the $\Gal(\bk/\k)$-action is as in
Figure~\ref{fig:action-dP6}(\ref{dp6:7})and (\ref{dp6:9}),
Lemma~\ref{prop:DP6} and Lemma~\ref{prop:DP8} imply that the pair
$(X,\Aut_\k(X))$ is as in Theorem~\ref{thm:1}(\ref{1:5c}).

If the $\Gal(\bk/\k)$-action is as in
Figure~\ref{fig:action-dP6}(\ref{dp6:2})--\ref{dp6:3}) and (\ref{dp6:5}),
then Lemma~\ref{prop:DP2} and Lemma~\ref{prop:DP3} and
Lemma~\ref{prop:DP9} imply that the pair $(X,\Aut(X))$ is as in
Theorem~\ref{thm:1}(\ref{1:5a}).

If the $\Gal(\bk/\k)$-action is as in
Figure~\ref{fig:action-dP6}(\ref{dp6:1}), Lemma~\ref{prop:DP1}
implies that $(X,\Aut(X))$ is as in
Theorem~\ref{thm:1}(\ref{1:51}).

 If the $\Gal(\bk/\k)$-action is as in Figure~\ref{fig:action-dP6}(\ref{dp6:4}), Lemma~\ref{prop:DP4} implies that $(X,\Aut(X))$ is as in  Theorem~\ref{thm:1}(\ref{1:52}).

If the $\Gal(\bk/\k)$-action is as in
Figure~\ref{fig:action-dP6}(\ref{dp6:6}), Lemma~\ref{prop:DP5}
implies that $(X,\Aut(X))$ is as in
Theorem~\ref{thm:1}(\ref{1:53}).

If the $\Gal(\bk/\k)$-action is as in
Figure~\ref{fig:action-dP6}(\ref{dp6:8}), Lemma~\ref{prop:DP7}
implies that $(X,\Aut(X))$ is as in
Theorem~\ref{thm:1}(\ref{1:54}).

Suppose that $X$ admits a conic fibration $\pi\colon X\to\p^1$ that
is an $\Aut(X,\pi)$-Mori fibre space and there is a birational
morphism $\eta\colon X\to Y$ where $Y=\F_n$ for some $n\geq0$ or
$Y=\Sl^{L,L'}$.

First, suppose that $\eta$ is an isomorphism. If
$X\stackrel{\eta}\simeq Y=\F_n$, recall that $\F_0$ and $\F_1$ have
already been discussed above, and that the family $\Aut(\F_n)$,
$n\geq2$ is the family in Theorem~\ref{thm:1}(\ref{1:4}), see
Remark~\ref{rmk:cb-F}. If $X\stackrel{\eta}\simeq Y=\Sl^{L,L'}$,
then $\Aut(\Sl^{L,L'},\pi)\subseteq\Aut(\Sl^{L,L'})$, and the pair
$(\Sl^{L,L'},\Aut(\Sl^{L,L'}))$ is as in
Theorem~\ref{thm:1}(\ref{1:5a}) by Lemma~\ref{prop:DP2}.

Now, suppose that $\eta$ is not an isomorphism. Since $\pi\colon
X\to\p^1$ is an $\Aut(X,\pi)$-Mori fibre space, there is an element
of $\Aut_{\bk}(X,\pi)$ that permutes the components of at least one
singular geometric fibre. If $X/\p^1$ has at most two singular
fibres, then the fact that $\eta$ is not an isomorphism implies that
$Y=\F_n$, and Lemma~\ref{lem:cbF-autbig} implies that $X$ is a del
Pezzo surface of degree $6$. Then $\Aut(X,\pi)\subseteq\Aut(X)$ and
we have already discussed the pair $(X,\Aut(X))$ above. If $X/\p^1$
has at least three singular fibres, recall that $\Aut_{\bk}(X,\pi)$
is infinite by hypothesis, and now Lemma~\ref{lem:thm1cb} implies
that the pair $(X,\Aut(X))$ is as in Theorem~\ref{thm:1}(\ref{1:6}).
\end{proof}

\section{Classifying maximal algebraic subgroups up to conjugacy}\label{s:6}

In this section we classify up to conjugacy and up to inclusion the
maximal infinite algebraic subgroups of $\Bir_\k(\p^2)$. For this,
we first need to introduce the so-called Sarkisov program. As
before, $\k$ is a perfect field throughout the section.

\subsection{The equivariant Sarkisov program}\label{ss:sarkisov}

The Sarkisov program is an algorithmic way to decompose birational
maps between Mori fibre spaces into nice elementary birational maps
between Mori fibre spaces. In dimension $2$, it is classical and
treated exhaustively in \cite{iskovskikh_1996}, and from a more
modern point of view in \cite{LZ17}. In dimension $3$, it is
developed in \cite{corti1995factoring} over algebraically closed
fields of characteristic zero. A non-algorithmic generalisation to
any dimension $\geq2$ is given in \cite{HMcK} over $\C$.

For surfaces, the Sarkisov program over $\k$ is the
$\Gal(\bk/\k)$-equivariant classical Sarkisov program over $\bk$.
For an affine algebraic group $G$, we can consider two equivariant
Sarkisov programs:
\begin{itemize}
\item The $G(\k)$-equivariant Sarkisov program over $\k$; the links are $G(\k)$-equivariant birational maps between $G(\k)$-Mori fibre spaces. If $G=\Aut(X)$ is one of the groups from Theorem~\ref{thm:1}, it is the tool to give us the conjugacy class of $G(\k)$ inside $\Bir_\k(\p^2)$.

\item The $G$-equivariant Sarkisov program is the $G_{\bk}\times\Gal(\bk/\k)$-equivariant Sarkisov program over $\bk$; the links are $G$-equivariant birational maps between $G$-Mori fibre spaces. If $G=\Aut(X)$ is one of the groups from Theorem~\ref{thm:1}, it is the tool to give us the morphisms $G\to\Bir_\k(\p^2)$ up to conjugation by an element of $\Bir_\k(\p^2)$.
\end{itemize}
As part of Theorem~\ref{thm:2}, we will prove that these two
classifications are not the same if $\k$ has an extension of degree
$2$ or $3$.

Over $\C$ and for connected algebraic groups $G$, the
$G$-equivariant Sarkisov program in dimension $\geq2$ is developed
in \cite{Floris_Sarkisov}. \vskip\baselineskip

\begin{defi}\label{def:links}
Let $G$ be an affine algebraic group. We now define $G(\k)$-equivariant
Sarkisov links. The notion of $G$-equivariant Sarkisov links is
defined analogously by replacing $G(\k)$ with $G$, bearing that by
$G$-orbit we mean a $G_{\bk}\times\Gal(\bk/\k)$-orbit.

A {\em $G(\k)$-equivariant Sarkisov link} (or simply {\em
$G(\k)$-equivariant link}) is a $G(\k)$-equivariant birational map
$\varphi\colon X\rat X'$ between $G(\k)$-Mori fibre spaces
$\pi\colon X\to B$ and $\pi'\colon X'\to B'$ that is one of the
following:
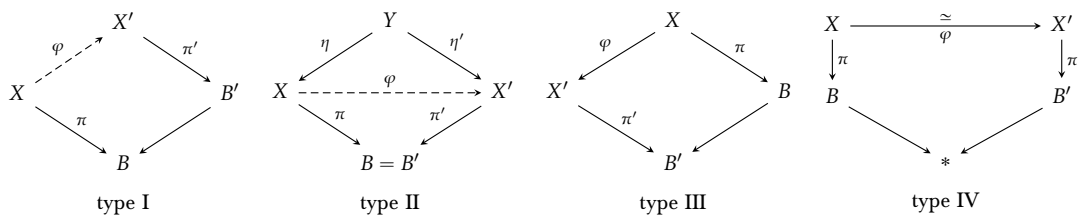
\[
\begin{tikzpicture}[baseline= (a).base]
\node[scale=.75](a) at (0,0){
\begin{tikzcd}
&X'\ar[dr,"\pi' "]&\\
X\ar[dr,"\pi"]\ar[ur,dashed,"\varphi"]&&B'\ar[dl]\\
&B&\\[-15pt]
&\text{type I}&
\end{tikzcd}
};
\end{tikzpicture}
\begin{tikzpicture}[baseline= (a).base]
\node[scale=.75](a) at (0,0){
\begin{tikzcd}
&Y\ar[dl,"\eta",swap]\ar[dr,"\eta'"]&\\
X\ar[dr,"\pi"]\ar[rr,dashed,"\varphi"]&& X'\ar[dl,"\pi'",swap]\\
&B=B'&\\[-15pt]
&\text{type II}&
\end{tikzcd}
};
\end{tikzpicture}
\begin{tikzpicture}[baseline= (a).base]
\node[scale=.75](a) at (0,0){
\begin{tikzcd}
&X\ar[dr,"\pi"]\ar[dl,swap,"\varphi"]&\\
X'\ar[dr,"\pi'"]&&B\ar[dl]\\
&B'&\\[-15pt]
&\text{type III}&
\end{tikzcd}
};
\end{tikzpicture}
\begin{tikzpicture}[baseline= (a).base]
\node[scale=.75](a) at (0,0){
\begin{tikzcd}
X\ar[rr,"\varphi" ', "\simeq"]\ar[d,"\pi"]&&X'\ar[d,"\pi'"]\\
B\ar[dr]&& B'\ar[dl]\\
&\ast&\\[-15pt]
&\text{type IV}&
\end{tikzcd}
};
\end{tikzpicture}
\]
\begin{itemize}
\item[(type I)] $B$ is a point, $B'$ is a curve, $\varphi^{-1}\colon X'\to X$ is the contraction of the $G(\k)$-orbit of a curve in $X'$ and $\pi\varphi^{-1}\colon X'\to B$ is a $G(\k)$-equivariant rank $2$ fibration (see Definition~\ref{def:GMfs}).
We call $\varphi$ a {\em link of type I}.
\item[(type II)] Either $B=B'$ is a curve or a point, both $\eta$ and $\eta'$ are contractions of the $G(\k)$-orbit of a curve and $\pi\eta\colon Y\to B$ is a $G(\k)$-equivariant rank $2$ fibration.
    We call $\varphi$ a {\em link of type II}.
\item[(type III)] $B$ is a curve, $B'$ is a point, $\varphi$ is the contraction of the $G(\k)$-orbit of a curve and $\pi'\varphi \colon X\to B$ is a $G(\k)$-equivariant rank $2$ fibration.
    We call $\varphi$ a {\em link of type III}. Its inverse is a link of type I.
\item[(type IV)] $B'$ and $B'$ are curves, $\varphi$ is an $G(\k)$-equivariant isomorphism not preserving the conic fibrations $X/B$ and $X'/B'$, and $X/\ast$ is a $G(\k)$-equivariant rank $2$ fibration. We call $\varphi$ a {\em link of type IV}.
\end{itemize}
\end{defi}

For $G=\{1\}$ we recover the classical definition of a Sarkisov link
over $\k$.

The statement of Theorem~\ref{thm:sarkisov} for $G=\{1\}$ is
\cite[Theorem 2.5]{iskovskikh_1996}. Its proof can be made
$G(\k)$-equivariant and $G$-equivariant because for a geometrically
rational variety $X$, the $G_{\bk}\times\Gal(\bk/\k)$ has finite
action on $\NS(X_{\bk})$ and $G(\k)$ has finite action on $\NS(X)$.

\begin{thm}[{Equivariant version of \cite[Theorem 2.5]{iskovskikh_1996}}]\label{thm:sarkisov}
Let $G$ be an affine algebraic group. Any $G(\k)$-equivariant
birational map between two geometrically rational surfaces that are
$G(\k)$-Mori fibre spaces is the composition of $G(\k)$-equivariant
Sarkisov links and isomorphisms.

The same statement holds if we replace $G(\k)$ by $G$.
\end{thm}

To study conjugacy classes of the automorphism groups of the
surfaces in Theorem~\ref{thm:1}, it therefore suffices to study
equivariant Sarkisov links between them.

\begin{rmk}\label{rmk:links-opt}
Definition~\ref{def:links} implies the following properties. Let
$\phi\colon X/B\rat X'/B'$ be an equivariant link.
\begin{enumerate}
\item If $\phi$ is a link of type I, then $B$ is a point, $X/B$ is an equivariant rank $1$  fibration above a point and $X'/B$ is an equivariant rank $2$ fibration above a point.
Equivariant rank $s$ fibrations above a point are in particular
(non-equivariant) rank $r$ fibrations above a point for some $r\geq
s$, see Definition~\ref{def:GMfs}, and so they are del Pezzo
surfaces, see Definition~\ref{def:Mfs}. So both $X$ and $X'$ are del
Pezzo surfaces. By symmetry, the same holds for a link of type III.
\item If $\phi$ is a link of type II and $B=B'$ a point, then $X/B$ and $X'/B$ are equivariant rank $1$ fibrations above a point, and $Y/B$ is an equivariant rank $2$ fibration above a point. Again, in particular, $X,X'$ and $Y$ are all del Pezzo surfaces.
\end{enumerate}
\end{rmk}

Many of the surfaces in Theorem~\ref{thm:1} are equivariant Mori
fibre spaces with respect to their automorphism group, as well as to
the group of $\k$-points of their automorphism group, and the
restrictions for the possible $\Aut_\k(X)$-links are also
restrictions on the possibilities of $\Aut(X)$-links.

We now classify the $\Aut_\k(X)$-equivariant links starting from a
surface $X$ from Theorem~\ref{thm:1} in the order
(\ref{1:1}--\ref{1:3}), (\ref{1:5c}), (\ref{1:52}--\ref{1:54}),
(\ref{1:51}), (\ref{1:4}) and (\ref{1:6}).

\subsection{$\Aut_\k(X)$-equivariant links of del Pezzo surfaces of degree $8$ and $9$}

We show that there are no $\Aut_\k(X)$-equivariant links starting
from a $\Aut_\k(X)$-Mori fibre space $X$ that is a rational del
Pezzo surface of degree $8$ or $9$.

\begin{lem}\label{lem:NoInvarPoints}
\begin{enumerate}
\item\label{NoInvarPoints:1} $\Aut_\k(\p^2)$ does not have any orbits in $\p^2$ with $d\in\{1,\dots,8\}$ geometric components that are in general position.
\item\label{NoInvarPoints:2} For $X=\F_0$ and $X=\Ql^L$, $\Aut_\k(X)$ does not have any orbits in $X$ with $d\in\{1,\dots,7\}$ geometric components that are in general position.
\end{enumerate}
\end{lem}
\begin{proof}
(\ref{NoInvarPoints:1}) Lemma~\ref{rmk:cyclicpoint2} implies the
claim for $1\leq d\leq 4$. If $\k$ is infinite and if
$\Aut_\k(\p^2)$ had an orbit with $5\leq d\leq8$ geometric
components, then $\alpha^{d!}=\id$ for any $\alpha\in\Aut_\k(\p^2)$,
which is false. Suppose that $\k$ is finite and let $q:=|\k|\geq2$.
Let $p=\{p_1,\dots,p_e\}$ be a point in $\p^2$ of degree $e\geq 5 $
and $L/\k$ be the smallest field extension such that
$p_1,\dots,p_e\in\p^2(L)$. We view $\Aut_\k(\p^2)$ as an abstract
subgroup of $\Aut_L(\p^2)$, which gives us
\[
1=|\cap_{i=1}^e\mathrm{Stab}_{\Aut_\k(\p^2)}(p_i)|=|\mathrm{Stab}_{\Aut_\k(\p^2)}(p_1)|=\frac{|\Aut_\k(\p^2)|}{|\text{$\Aut_\k(\p^2)$-orbit
of $p_1$ in $\p^2(L)$}|}.
\]
Moreover, we have $|\Aut_\k(\p^2)|=q^3(q^3-1)(q^2-1)> q^3\geq 8$,
and hence the $\Aut_\k(\p^2)$-orbit of $p$ in $\p^2$ has $\geq9$
geometric components.

(\ref{NoInvarPoints:2}) For $X=\F_0$ and $d=1,2$, the claim follows
from Remark~\ref{rmk:cyclicpoint2-dim1}. For $X=\Ql^{L}$, the claim
follows from Remark~\ref{rmk:cyclicpoint2-dim1} for $d=1$, from
Lemma~\ref{lem:Q-deg2} for $d=2$. Let $L/\k$ be a quadratic
extension such that $\Ql^L_L\simeq\p^1_L\times\p^1_L$, and by
Lemma~\ref{lem:autQ} we have
$\Aut_\k(\Ql^L)\simeq\PGL_2(L)\rtimes\Z/2$. For $3\leq d\leq 7$, we
can repeat the argument of (\ref{NoInvarPoints:1}) for $\F_0$ and
$\Ql^L$ by using that for a finite field $\k$ with $q:=|\k|\geq2$ we
have
\begin{align*}
&|\Aut_\k(\F_0)|=2|\PGL_2(\k)|^2=2q^2(q^2-1)^2>8\\
&|\Aut_\k(\Ql^L)|= 2|\PGL_2(L)|=2q^2(q^4-1)>8.
\end{align*}
\end{proof}

\begin{lem}\label{lem:linksPQF}
There is no $\Aut_\k(X)$-equivariant link starting from $X=\p^2$,
$X=\Ql^L$ or $X=\F_0$.
\end{lem}
\begin{proof}
Since $\rk\,\NS(X)^{\Aut_\k(X)}=1$, the only $\Aut_\k(X)$-equivariant
links starting from $X$ are of type I or II. Moreover,
$\Aut_\k(\F_0)$-equivariant links starting from $\F_0$ can be
treated like the ones starting from $\Ql^L$ because
$\NS(\F_0)^{\Aut_\k(\F_0)}=\Z(f_1+f_2)=\NS(\Ql^L)$, where $f_1,f_2$
are the fibres of the two projections of $\F_0$.

By Remark~\ref{rmk:links-opt}, an $\Aut_\k(\p^2)$-equivariant link
of type I or II starting from $\p^2$ blows up an orbit with $\leq8$
geometric components that are in general position, and by
Lemma~\ref{lem:NoInvarPoints}(\ref{NoInvarPoints:1}), there is no
such orbit. An $\Aut_\k(X)$-equivariant link of type I or II
starting from $X=\Ql^L$ or $X=\F_0$ blows up an orbit with $\leq7$
geometric components that are in general position, and by
Lemma~\ref{lem:NoInvarPoints}(\ref{NoInvarPoints:2}), there is no
such orbit.
\end{proof}

\subsection{$\Aut_\k(X)$-equivariant links of del Pezzo surfaces of degree $6$ (\ref{1:5c})}
These del Pezzo surfaces are Mori fibre spaces. We will show that
there are no $\Aut_\k(X)$-equivariant links starting from $X$.

Recall from Lemma~\ref{prop:DP6} and Lemma~\ref{prop:DP8} that there
is a quadratic extension $L/\k$ such that $X_L$ is obtained by
blowing up a point $p=\{p_1,p_2,p_3\}$ in $\p^2$ of degree $3$. We
denote by $\pi\colon X_L\to\p^2_L$ the blow-up of $p$. Recall that
$\pi\Gal(L/\k)\pi^{-1}$ acts rationally on $\p^2$; its generator
$\psi_g$ is not defined at $p$ and sends a general line onto a conic
through $p$. Recall that if $X$ is rational, it has a rational point
by Proposition~\ref{lem:sing-fibres}.

\begin{lem}\label{lem:fixedpoints}
Let $X$ be a del Pezzo surface of degree $6$ from
Theorem~\ref{thm:1}(\ref{1:5c}) and fix $s\in X(\k)$. The map
\[
 \Aut_L(\p^2,p_1,p_2,p_3)^{\langle \psi_g\rangle}\to X(\k),\quad \alpha\mapsto \pi^{-1}(\alpha(\pi(s))=(\pi^{-1}\alpha\pi)(s)
\]
is bijective.
\end{lem}
\begin{proof}
The map is injective, because these automorphisms already fix
$p_1,p_2,p_3$. For any $t\in X(\k)$, we have $\pi(t)\in\p^2_L(L)$,
and by Lemma~\ref{rmk:cyclicpoint2} there exists a unique element of
$\alpha_t\in\Aut_L(\p^2,p_1,p_2,p_3)$ such that
$\alpha_t(\pi(s))=t$. Then $\pi^{-1}\alpha_t\pi\in\Aut_L(X)$ and its
conjugate by the generator of $\Gal(L/\k)$ is still contained in
$\Aut_L(X)$ and preserves each edge of the hexagon, hence
$\psi_g\alpha_t\psi_g\alpha_t^{-1}\in\Aut_L(\p^2,p_1,p_2,p_3)$. The
automorphism $\psi_g\alpha_t\psi_g\alpha_t^{-1}$ fixes
$p_1,p_2,p_3,\pi(t)$, so it is the identity, and therefore
$\alpha_t\in \Aut_L(\p^2,p_1,p_2,p_3)^{\langle \psi_g\rangle}$.
\end{proof}

\begin{lem}\label{lem:size}
Let $X$ be a del Pezzo surface of degree $6$ from
Theorem~\ref{thm:1}(\ref{1:5c}). Then $|X(\k)|\geq7$ if $|\k|\geq3$
and $|X(\k)|=3$ if $|\k|=2$. Moreover, in the latter case the
blow-up of $X(\k)$ is a del Pezzo surface.
\end{lem}
\begin{proof}
If $\k$ is infinite, then $\p^2(\k)$ is dense in $\p^2(\bk)$, and
hence $X(\k)$ is infinite. If $\k$ is finite, pick a rational point
$r\in X(\k)$. There exists a link of type II $\phi\colon X\rat
\Ql^L$ that is not defined at $r$ and contracts a curve with three
geometric components passing through $r$, see Figure~\ref{fig:rot6}.
If $Z\to X$ is the blow-up of $r$ and $L/\k$ a quadratic extension
such that $\Ql_L=\p^1_L\times\p^1_L$, we have
\[
q^2+1=|\p^1(L)|=|\Ql^L(\k)|=|Z(\k)|=|X(\k)|-1+|\p^1(\k)|=|X(\k)|+q
\]
because the exceptional divisor of $r$ is isomorphic to $\p^1_{\k}$.
It follows that $|X(\k)|=q^2-q+1=q(q-1)+1$.

Suppose now that $|\k|=2$ and so $|X(\k)|=3$. Then $X(\k)$ is the
image of the five points $\Ql^L(\k)$ by $\phi$, and it suffices to
show that the blow-up of $\Ql^L(\k)$ is a del Pezzo surface. We
write $L=\k(a)$, where $a^2+a+1=0$. The set $\Ql^L(\k)$ consists of
\[
([1:0],[1:0]),([0:1],[0:1]),([1:1],[1:1]),([1:a],[1:a^2]),([1:a^2],[1:a])
\]
and we check that they are not contained in any fibre of $\Ql^L_L$
nor in any bidegree $(1,1)$-curve. This yields the claim.
\end{proof}

\begin{lem}\label{lem:5c}
Let $X$ be a rational del Pezzo surface as in
Theorem~\ref{thm:1}(\ref{1:5c}).
\begin{enumerate}
\item\label{5c:1} If $|\k|\geq3$, then $X$ does not contain any $\Aut_\k(X)$-orbits with $\leq5$ geometric components.
\item\label{5c:2} If $|\k|=2$, there is exactly one $\Aut_\k(X)$-orbit of $X$ with $\leq5$ geometric components, namely $X(\k)$.
\end{enumerate}
\end{lem}
\begin{proof}
Since $\Gal(\bk/\k)$ acts transitively on the edges of the hexagon,
any orbit with $\leq5$ geometric components is outside of it. Let
$D\subset\p^2_L$ be the image of the hexagon by $\pi$.

Suppose that $|\k|\geq3$. By Lemma~\ref{lem:size}, we have
$|X(\k)|\geq7$, so Lemma~\ref{lem:fixedpoints} implies that the
group $\Aut_L(\p^2,p_1,p_2,p_3)^{\langle \psi_g\rangle}$ has $\geq7$
elements. It acts faithfully on $\p^2\setminus D$, hence any
$\Aut_L(\p^2,p_1,p_2,p_3)^{\langle \psi_g\rangle}$-orbit in
$\p^2\setminus D$ has $\geq7$ geometric components. It follows that
$\Aut_\k(X)$ has no orbits with $\leq5$ geometric components on $X$.

Suppose now that $|\k|=2$ and let $L/\k$ be the extension of degree
$2$. We show that $\pi\Aut_\k(X)\pi^{-1}$-orbit of any point in
$\p^2_L\setminus D$ has either $3$ or $\geq6$ elements, and that
$\pi(X(\k))$ is the only orbit with $3$ elements. Let
$\varphi_p\in\Bir_L(\p^2)$ be the quadratic involution from
Lemma~\ref{prop:DP6}(\ref{DP6:4}) and
Lemma~\ref{prop:DP8}(\ref{DP8:4}) that lifts to an automorphism
$\tilde{\varphi}_p=\pi^{-1}\varphi_p\pi$ on $X$ over $\k$ inducing a
rotation of order $2$ on the hexagon of $X$. By
Lemma~\ref{prop:DP6}(\ref{DP6:4}) (resp.
Lemma~\ref{prop:DP8}(\ref{DP8:4})) the group
\[ \Aut_L(\p^2,p_1,p_2,p_3)^{\langle \psi_g\rangle}\rtimes\langle\varphi_p\rangle
\]
is isomorphic to a subgroup of $\Aut_\k(X)$. Lemma~\ref{lem:size}
and Lemma~\ref{lem:fixedpoints} imply that
$\Aut_L(\p^2,p_1,p_2,p_3)^{\langle \psi_g\rangle}$ has $3$ elements,
and it acts faithfully on $\p^2_L\setminus D$. Over $\bk$, the
involution $\varphi_p$ is conjugate to the involution
$[x:y:z]\rat[yz:xz:xy]$, which has a unique fixed point in
$\p^2_{\bk}$, namely $[1:1:1]$, because $|\k|=2$. Thus $\varphi_p$
has a unique fixed point $r\in\p^2_L$. Then $\tilde{r}:=\pi^{-1}(r)$
is the unique fixed point of $\tilde{\varphi}_p$ on $X$, and it is
$\k$-rational. We have shown that every $\Aut_\k(X)$-orbit in
$X(L)\setminus X(\k)$ has $\geq6$ elements. The set $X(\k)$ is an
$\Aut_\k(X)$-orbit with $3$ elements.
\end{proof}

\begin{lem}
\label{lem:link5a2} Let $|\k|=2$ and let $X$ be a del Pezzo surface
from Theorem~\ref{thm:1}(\ref{1:5c1}). Any $\Aut_\k(X)$-invariant
link $\varphi\colon X\rat Y$ is a link of type II not defined at
$X(\k)$, and $Y$ is a del Pezzo surface as in
Theorem~\ref{thm:1}(\ref{1:52}).
\end{lem}
\begin{proof}
We have $X(\k)=\{r_1,r_2,r_3\}$, see Lemma~\ref{lem:size}, which is
an $\Aut_\k(X)$-orbit by Lemma~\ref{lem:5c}. For a point $s\in
S:=\{\pi(r_1),\pi(r_2)$, $\pi(r_3),p_1,p_2,p_3\}\subset\p^2_L$, we
denote by $C_{s}$ the strict transform of the conic in $\p^2_L$
passing through the five points in $S\setminus\{s\}$, and let
$L_{r_ir_j}$ be the strict transform of the line in $\p^2_L$ through
$\pi(r_i),\pi(r_j)$, $i\neq j$. The curves
\[C_p:=C_{p_1}\cup C_{p_2}\cup C_{p_3},\ D_1:=C_{r_1}\cup L_{r_2r_3},\ D_2:=C_{r_2}\cup L_{r_1r_3},\ D_3:=C_{r_3}\cup L_{r_1r_2}\]
and $L_i:=L_{r_ip_1}\cup L_{r_ip_2}\cup L_{r_ip_3}$, $i=1,2,3$, are
irreducible over $\k$. The curve $C_p$ is $\Aut_\k(X)$-invariant,
while $D_1,D_2,D_3$ and $L_1,L_2,L_3$ make up an $\Aut_\k(X)$-orbit,
see Lemma~\ref{prop:DP6}(\ref{DP6:4}) for the generators of
$\Aut_\k(X)$.

Let $\eta\colon Z\to X$ be the blow-up of $X(\k)$, which is
$\Aut_\k(X)$-equivariant by Lemma~\ref{lem:5c}. The surface $Z$ is a
del Pezzo surface of degree $3$ by Lemma~\ref{lem:size}. There is at
most one way to complete $\eta$ into an $\Aut_\k(X)$-equivariant
link, because $Z$ is an $\Aut_\k(X)$-equivariant rank $2$ fibration,
and hence there are at most two extremal $\Aut_\k(X)$-equivariant
contractions from $Z$. However, any conic fibration $Z\to\p^1$ is
given by the fibres of the strict transforms of conics through four
fixed points in $S$ or the strict transform of lines through one
point in $\p^2_L$, but none of them are $\Aut_\k(X)$-equivariant. So
the link $\varphi$ has to be of type II.

The only $\Aut_\k(X)\times\Gal(\bk/\k)$-orbits of $(-1)$-curves on
$Z_{\bk}$ with $\leq 6$ geometric components which are pairwise
disjoint are the exceptional divisors of $\eta$ and the strict
transform of $C_p$. The contraction $\eta'\colon Z\to Y$ of the
latter induces an $\Aut_\k(X)$-equivariant link $X\rat Y$ to a del
Pezzo surface $Y$ of degree $6$.

Since the strict transforms of $C_{p_i}$ and $C_{r_j}$ on $Z_{\bk}$
are disjoint for $i,j=1,2,3$, the hexagon of $Y$ consists in the curve
$\eta'(D_1)\cup\eta'(D_2)\cup\eta'(D_3)$. Each
component $\eta'(D_i)$ of this union is $\k$-rational, so $\Gal(\bk/\k)$
acts as rotation of order $2$ on the hexagon of $Y$, {\it i.e.} as
in Figure~\ref{fig:action-dP6}(\ref{dp6:4}). By
Lemma~\ref{prop:DP4}, $Y$ is described in
Theorem~\ref{thm:1}(\ref{1:52}).
\end{proof}

\begin{prop}\label{prop:5c}
Let $X$ be a del Pezzo surface from Theorem~\ref{thm:1}(\ref{1:5c}).
Then, if $|\k|\geq3$, there are no $\Aut_\k(X)$-equivariant links
starting from $X$. If $|\k|=2$, the only $\Aut_\k(X)$-equivariant
link is the one from Lemma~\ref{lem:link5a2}.
\end{prop}
\begin{proof}
Since $\rk\,\NS(X)=1$, only $\Aut_\k(X)$-equivariant links of type I
or II can start from $X$. By Remark~\ref{rmk:links-opt}, they are
not defined at an orbit with $\leq5$ geometric components. By
Lemma~\ref{lem:5c}, such an orbit only exists for surfaces $X$ as in
Theorem~\ref{thm:1}(\ref{1:5c}) if $|\k|=2$. The claim now follows
from Lemma~\ref{lem:link5a2}.
\end{proof}

\subsection{$\Aut_\k(X)$-equivariant links of del Pezzo surfaces of degree $6$ (\ref{1:52})--(\ref{1:54})}

Any del Pezzo surface $X$ of degree $6$ from
Theorem~\ref{thm:1}(\ref{1:52})--(\ref{1:54}) is a $\Aut_\k(X)$-Mori
fibre space, and we show that there are no $\Aut_\k(X)$-equivariant
links starting from $X$.

\begin{lem}\label{lem:orbit52}
Let $X$ be a del Pezzo surface of degree $6$ from Theorem~\ref{thm:1}(\ref{1:52}). Then any $\Aut_\k(X)$-orbit on $X$ has at least $6$ geometric components.
\end{lem}
\begin{proof}
Let $\pi\colon X\to\F_0$ be the contraction of a curve in the
hexagon onto the point $p=\{(p_1,p_1),(p_2,p_2)\}$ of degree $2$
with $p_i=[a_i:1]$, $i=1,2$. Since $\Aut_\k(X)$ acts by
$\sym_3\times\Z/2$ on the hexagon of $X$, any orbit with $\leq5$
geometric components is outside of the hexagon. Let $D\subset\F_0$
be the image by $\pi$ of the hexagon, which contains $p$, and
consider the action of $\pi\Aut_\k(X)\pi^{-1}$ on $\F_0\setminus D$.
The elements of $\Aut_\k(\p^1,p_1,p_2)$ are exactly those of the
form
$$
[u:v]\mapsto[(b(a_1+a_2)+c)u -ba_1a_2v : bu+cv],\quad [b:c]\in\p^1(\k)
$$
and thus
\[|\Aut_\k(\p^1,p_1,p_2)|^2=|\p^1(\k)|^2\geq3^2=9.\]
Any non-trivial element of $\Aut_\k(\p^1,p_1,p_2)$ has precisely two
fixed points in $\p^1$. It follows that the stabiliser in
$\Aut_\k(\p^1,p_1,p_2)^2$ of any point $p_3\in(\F_0)_{\bk}\setminus
D_{\bk}$ is trivial and hence
\[|\text{$\Aut_\k(\p^1,p_1,p_2)^2$-orbit of $p_3$ in $(\F_0\setminus D)_{\bk}$}|=|\Aut_\k(\p^1,p_1,p_2)^2|\geq9.\]
We have shown that $\Aut_\k(\p^1,p_1,p_2)^2$ has no orbits on
$\F_0\setminus D$ with $\leq5$ geometric components, and hence that
$\pi\Aut_\k(X)\pi^{-1}$ has not orbits on $\F_0\setminus D$ with
$\leq5$ geometric components.
\end{proof}

\begin{rmk}\label{lem:aut2}
Let $p=\{p_1,p_2,p_3\}$ be a point of degree $3$ in $\p^2$. Fix a
point $r\in\p^2(\k)$. In particular, the point $r$ is not collinear
with any two components of $p$, and so Lemma~\ref{rmk:cyclicpoint2}
implies that the map $\Aut_\k(\p^2,p_1,p_2,p_3)\to\p^2(\k)$,
$\alpha\mapsto\alpha(r)$ is a bijection.
\end{rmk}

\begin{lem}\label{lem:orbit53}
Let $X$ be a del Pezzo surface of degree $6$ from Theorem~\ref{thm:1}(\ref{1:53}). Then any $\Aut_\k(X)$-orbit on $X$ has $\geq6$ geometric components.
\end{lem}
\begin{proof}
Since $\Aut_\k(X)$ contains an element inducing a rotation of order
$6$ on the hexagon of $X$, the hexagon does not contain
$\Aut_\k(X)$-orbits with $\leq5$ geometric components. Consider the
contraction $\pi\colon X\to\p^2$ of a curve in the hexagon of $X$
onto the point $p=\{p_1,p_2,p_3\}$ of degree $3$, let $D\subset\p^2$
be the image of the hexagon and consider the action of
$\Aut_\k(\p^2,p_1,p_2,p_3)\subset\pi\Aut_\k(X)\pi^{-1}$ on
$\p^2\setminus D$. Remark~\ref{lem:aut2} implies that
$|\Aut_\k(\p^2,p_1,p_2,p_3)|=|\p^2(\k)|\geq7$. The stabiliser of
$\Aut_\k(\p^2,p_1,p_2,p_3)$ of any point in $(\p^2\setminus
D)_{\bk}$ is trivial, so in particular all the
$\Aut_\k(\p^2,p_1,p_2,p_3)$-orbits in $\p^2\setminus D$ have $\geq7$
geometric components. It follows that $\pi\Aut_\k(X)\pi^{-1}$ has no
orbits in $\p^2\setminus D$ with $\leq5$ geometric components.
\end{proof}

\begin{lem}\label{lem:orbit54}
Let $X$ be a del Pezzo surface of degree $6$ from
Theorem~\ref{thm:1}(\ref{1:54}). The blow-up of $X$ in any finite
$\Aut_\k(X)$-orbit is not a del Pezzo surface.
\end{lem}
\begin{proof}
Let $\pi\colon X\to\p^2$ be the contraction of a curve $C$ in the
hexagon of $X$ onto the point $p=\{p_1,p_2,p_3\}$ of degree $3$. By
hypothesis, the splitting field $L/\k$ of $p$ satisfies
$\Gal(L/\k)\simeq\sym_3$, so $\k$ is not finite \cite[Theorem
6.5]{Morandi}. Remark~\ref{lem:aut2} implies that
$\Aut_\k(\p^2,p_1,p_2,p_3)$ is infinite. Let $D\subset\p^2$ be the
image by $\pi$ of the hexagon and consider the action of
$\Aut_\k(\p^2,p_1,p_2,p_3)\subset\pi\Aut_\k(X)\pi^{-1}$ on
$\p^2\setminus D$. The stabiliser of $\Aut_\k(\p^2,p_1,p_2,p_3)$ of
any point in $(\p^2\setminus D)_{\bk}$ is trivial, and hence any
$\Aut_\k(\p^2,p_1,p_2,p_3)$-orbit on $\p^2\setminus D$ has
infinitely many geometric components. It follows that any
$\Aut_\k(X)$-orbit with finitely many geometric components is
contained in the hexagon of $X$, and so its blow-up is not a del
Pezzo surface.
\end{proof}

\begin{prop}\label{pro:links-list2-4}
There is no $\Aut_\k(X)$-equivariant link starting from a del  Pezzo
surface $X$ of degree $6$ as in
Theorem~\ref{thm:1}$(\ref{1:52})-(\ref{1:54})$.
\end{prop}
\begin{proof}
Since $\rk\,\NS(X)^{\Aut_\k(X)}=1$, the only $\Aut_\k(X)$-equivariant
links starting from $X$ are of type I or II, and by
Remark~\ref{rmk:links-opt}, they are not defined in an
$\Aut_\k(X)$-orbit with $\leq5$ geometric components and its blow-up
is a del Pezzo surface. If $X$ is as in
Theorem~\ref{thm:1}(\ref{1:52})--(\ref{1:53}) no such orbit exists
respectively by Lemma~\ref{lem:orbit52} and Lemma~\ref{lem:orbit53}.
If $X$ is as in Theorem~\ref{thm:1}(\ref{1:54}), then the blow-up of
any such orbit is not a del Pezzo surface by
Lemma~\ref{lem:orbit54}.
\end{proof}

\subsection{$\Aut_\k(X)$-equivariant links of del Pezzo surfaces of degree $6$ (\ref{1:51})}
Studying $\Aut_\k(X)$-equivariant links for such a del Pezzo surface
is a bit more involved. We will show that there are equivariant
links starting from $X$ only if $|\k|=2$ and provide examples.
Recall Lemma~\ref{prop:DP1} for a description of $X$.

\begin{lem}\label{lem:permutations}
Fix homogeneous coordinates in $\p^2$ and consider the subgroup
$H\subset\PGL_3(\k)$ of permutation matrices. If the $H$-orbit $O$
of a point in $\{xyz\neq0\}\subset\p^2$ has $\leq5$ geometric
components, it is one of the following:
\begin{enumerate}
\item\label{perm:1} $O=\{[1:1:1]\}$,
\item\label{perm:2} $O=\{[1:a:a^2],[1:a^2:a]\}$ with $a^3=1$,
\item\label{perm:3} $O=\{[1:a:a],[a:a:1],[a:1:a]\}$ for some $a\in\k^*$.
\end{enumerate}
\end{lem}
\begin{proof}
The $H$-orbit $O_{\bk}$ of a point $p:=[1:a:b]\in\{xyz\neq0\}_{\bk}$
is contained in the set {\small
\begin{align*}&\{[1:a:b],[1:b:a],[a:b:1],[b:a:1],[a:1:b],[b:1:a]\}\\
=&\{[1:a:b],[1:b:a],[1:a^{-1}b:a^{-1}],[1:ab^{-1}:b^{-1}],[1:a^{-1}:a^{-1}b],[1:b^{-1}:ab^{-1}]\}
\end{align*}
} If $p$ is an $H$-fixed point, we have $O_{\bk}=O=\{[1:1:1]\}$. We
check that if $|O_{\bk}|=2$, then we have $O_{\bk}=\{[1:a:a^2],[1:a^2:a]\}$
with $a^3=1$. If $|O_{\bk}|=3$, then
$O_{\bk}=\{[1:1:c],[1:c:1],[1:c^{-1}:c^{-1}]\}$ for some $c\in\k^*$.
We also check that $4\leq|O_{\bk}|\leq5$ is not possible.
\end{proof}

\begin{lem}\label{lem:orbit51}
Let $X$ be the del Pezzo surface of degree $6$ from
Theorem~\ref{thm:1}(\ref{1:51}).
\begin{enumerate}
\item\label{orbit51:2}  If $|\k|\geq4$, then $X$ contains no $\Aut_\k(X)$-orbits with $\leq5$ geometric components.
\item\label{orbit51:4}  If $|\k|=3$, then $\Aut_\k(X)$ has exactly one orbit on $X$ with $\leq5$ geometric components, namely the orbit $\{([1:\pm1:\mp1],[1:\pm1:\mp1])\}$ with $4$ elements. Its blow-up is not a del Pezzo surface.
\item\label{orbit51:3}  If $|\k|=2$, then $\Aut_\k(X)$ has exactly two orbits on $X$ with $\leq5$ geometric components, namely the fixed point $([1:1:1],[1:1:1])$ and the point $\{\left([1:\zeta:\zeta^2],[1:\zeta^2:\zeta]\right)$, $ \left([1:\zeta^2:\zeta],[1:\zeta:\zeta^2]\right)\}$ of degree $2$, where $\zeta\notin\k$, $\zeta^3=1$.
\end{enumerate}
\end{lem}
\begin{proof}
By Lemma~\ref{prop:DP1}(\ref{DP1:2}), the group $\Aut_\k(X)$ acts
transitively on the edges of the hexagon, so the hexagon does not
contain $\Aut_\k(X)$-orbits with $\leq5$ geometric components. We
pick three disjoint edges of the hexagon and consider their
contraction $\pi\colon X\to\p^2$ onto the coordinate points, which
maps the hexagon onto the curve $\{xyz=0\}$. It remains to study the
$\pi\Aut_\k(X)\pi^{-1}$-action on $\{xyz\neq0\}$. The stabiliser
subgroup of the subgroup $(\k^*)^2\subset\pi\Aut_\k(X)\pi^{-1}$ of
diagonal elements of any point in $\{xyz\neq0\}$ is trivial. It
follows that the $(\k^*)^2$-orbit of any point in $\p^2$ has $\geq9$
geometric components if $|\k^*|\geq3$, proving (\ref{orbit51:2}).

Let $2\leq|\k|\leq3$ and recall from
Lemma~\ref{prop:DP1}(\ref{DP1:2}) that
$\pi\Aut_\k(X)\pi^{-1}\simeq(\k^*)^2\rtimes(H\times\Z/2)$, where
$H=\pi\sym_3\pi^{-1}$ is the group of permutation matrices in
$\Aut_\k(\p^2)$ and $\Z/2$ is generated by the involution
$(x,y)\dashmapsto(\frac{1}{x},\frac{1}{y})$.

If a $\pi\Aut_\k(X)\pi^{-1}$-orbit in $\{xyz\neq0\}$ has $\leq5$
geometric components, then this holds in particular for an $H$-orbit
$O$, which is one of the following by Lemma~\ref{lem:permutations}
\begin{enumerate}
\item[(i)] $O=\{[1:1:1]\}$,
\item[(ii)] $O=\{[1:a:a^2],[1:a^2:a]\}$ with $a^3=1$,
\item[(iii)] $O=\{[1:a:a],[1:1:a^{-1}],[1:a^{-1}:1]\}$ for some $a\in\k^*$.
\end{enumerate}

(\ref{orbit51:3}) If $|\k|=2$, then
$\pi\Aut_\k(X)\pi^{-1}\simeq(H\times\Z/2)$ and the point $[1:1:1]$
is a fixed point and is equal to (iii) and (ii) for $a=1$. If
$a\notin\k$ and $a^3=1$, the point $\{[1:a:a^2],[1:a^2:a]\}$ of
degree $2$ is a $\pi\Aut_\k(X)\pi^{-1}$-fixed point.

(\ref{orbit51:4}) If $|\k|=3$, then the
$\pi\Aut_\k(X)\pi^{-1}$-orbit of $[1:1:1]$ is the set
$O=\{[1:\pm1:\pm1]\}$, which has $4$ elements. The
$\pi\Aut_\k(X)\pi^{-1}$-orbit of a point in (ii) or (iii) is either
the orbit of $[1:1:1]$ or has $\geq6$ geometric components. The line
$\{y=z\}\subset\p^2$ contains $[1:0:0],[1:-1:-1],[1:1:1]$, so the
blow-up of $X$ in $\pi^{-1}(O)$ is not a del Pezzo surface.
\end{proof}

\begin{lem}\label{lem:rk1-2}
Let $|\k|=2$ and let $X$ be the del Pezzo surface of degree $6$ from
Theorem~\ref{thm:1}(\ref{1:51}). The blow-up of $X$ in any
$\Aut_\k(X)$-orbit with $\leq5$ geometric components does not admit
a $\Aut_\k(X)$-equivariant fibration over $\p^1$.
\end{lem}
\begin{proof}
Let $\pi\colon X\to\p^2$ be the blow-up of the coordinate points
$p_1,p_2,p_3$. By Lemma~\ref{lem:orbit51}(\ref{orbit51:3}), the only
$\Aut_\k(X)$-orbits on $X$ with $\leq5$ geometric components are a
fixed-point $r\in X(\k)$ and a point $q\in X$ of degree $2$, both
not on the hexagon.

Let $Y\to X$ be the blow-up of $r$ and let $Y/\p^1$ be a conic
fibration. Its fibres are either the strict transform of the lines
through one of $p_1,p_2,p_3,r$, or the strict transform of the
conics through $p_1,p_2,p_3,r$. Since
$\Aut_\k(X)\simeq\sym_3\times\Z/2$ acts transitively on the edges of
the hexagon of $X$ by Lemma~\ref{prop:DP1} and the quadratic
involution in $\pi\Aut_\k(X)\pi^{-1}$ sends a general line through
$r$ onto a conic through $p_1,p_2,p_3,r$, it follows that $Y/\p^1$
is not $\Aut_\k(X)$-equivariant.

Let $Y\to X$ be the blow-up of $q$ and $Y/\p^1$ a conic fibration.
Its fibres are the strict transforms of the conics through $q$ and
two of $p_1,p_2,p_3$ or of a line through one of $p_1,p_2,p_3$.
Again, as $\Aut_\k(X)$ acts transitively on the edges of the hexagon
of $X$, it follows that $Y/\p^1$ is not $\Aut_\k(X)$-equivariant.
\end{proof}

\begin{ex}\label{ex:links6}
Let $\pi\colon X\to\p^2$ be the blow-up of the coordinate points
$p_1,p_2,p_3$ of $\p^2$. If $|\k|=2$, then by
Lemma~\ref{prop:DP1}(\ref{DP1:2}) the group
$\pi\Aut_\k(X)\pi^{-1}\simeq\sym_3\times\Z/2$ is generated by
\[
\alpha\colon[x:y:z]\mapsto[x:z:y],\quad
\beta\colon[x:y:z]\mapsto[z:y:x],\quad\sigma\colon(x,y)\dashmapsto(\frac{1}{x},\frac{1}{y})
\]
\begin{enumerate}
\item\label{links6:1}
If $\mathrm{char}(\k)=2$, the birational map $\psi_1\colon
\p^2\rat\F_0$
\begin{align*}
\psi_1\colon [x:y:z]&\dashmapsto([x-z:y-z],[y(x-z):x(y-z)]),\\
\psi_1^{-1}\colon([u_0:u_1],[v_0:v_1])&\dashmapsto[u_0(u_0+u_1)v_1:
u_1(u_0+u_1)v_0 : u_0u_1(v_0+v_1)]
\end{align*}
is not defined at $p_1,p_2,p_3,[1:1:1]$ and contracts the
$\pi\Aut_\k(X)\pi^{-1}$-orbit $\{(y-z)(x-z)(x-y)=0\}$. If $|\k|=2$,
it lifts to an $\Aut_\k(X)$-birational map
\[\varphi_1:=\psi_1\pi\colon X\rat\F_0\]
not defined at $\pi^{-1}([1:1:1])$, because
\begin{align*}
\psi_1\alpha\psi_1^{-1}\colon &([u_0:u_1],[v_0:v_1])\mapsto([u_0+u_1:u_1],[v_0+v_1:v_1]),\\
\psi_1\beta\psi_1^{-1}\colon &([u_0:u_1],[v_0:v_1])\mapsto([u_0:u_0+u_1],[v_0:v_0+v_1]),\\
\psi_1\sigma\psi_1^{-1}\colon&([u_0:u_1],[v_0:v_1])\mapsto([v_0:v_1],[u_0:u_1])
\end{align*}
are automorphisms of $\F_0$. So $\varphi_1\colon X\rat \F_0$ is an
$\Aut_\k(X)$-equivariant link of type II.

\item\label{links6:2}
Let $\mathrm{char}(\k)=2$ and $\zeta\in\bk\setminus\k$, $\zeta^3=1$
and $q:=\{[1:\zeta:\zeta^2],[1:\zeta^2:\zeta]\}$. The birational map
$\psi_2\colon \p^2\rat\F_0$
\begin{align*}
\psi_2\colon [x:y:z]\dashmapsto([xy+xz+yz:y(x+y+z)],[xy+xz+yz:z(x+y+z)],\\
\psi_2^{-1}\colon([u_0:u_1],[v_0:v_1])\dashmapsto[u_0v_0(u_1v_0+u_0v_1+u_1v_1): u_1v_0(u_1v_0+u_0v_1+u_0v_0):\\
u_0v_1(u_1v_0+u_0v_1+u_0v_0)]
\end{align*}
is not defined at $p_1,p_2,p_3,q$ and contracts the rational
curves $\{(x+y+z)(xy+xz+yz)=0\}$, and the conic $\{y^2+yz+z^2=0\}$
onto $q':=\{([1:\zeta],[1:\zeta^2]),([1:\zeta^2],[1:\zeta])\}$. Let
$\eta\colon X'\to\F_0$ be the blow-up of $q'$, which is a del Pezzo
surface of degree $6$ as in Lemma~\ref{prop:DP4}
(Figure~\ref{fig:action-dP6}(\ref{dp6:4})). If $|\k|=2$, the
contracted curves are $\Aut_\k(X)$-invariant and $\psi_2$ lifts to
an $\Aut_\k(X)$-equivariant birational map
\[\varphi_2:=\eta^{-1}\psi_2\pi\colon X\rat X'\]
not defined at $\pi^{-1}(q)$. Consider the conjugates
\begin{align*}
\psi_2\alpha\psi_2^{-1}\colon &([u_0:u_1],[v_0:v_1])\mapsto([v_0:v_1],[u_0:u_1]),\\
\psi_2\beta\psi_2^{-1}\colon &([u_0:u_1],[v_0:v_1])\dashmapsto([u_0:u_1],[u_0v_0+(u_1v_0+u_0v_1) : u_1v_1+(u_0v_1+u_1v_0)]),\\
\psi_2\sigma\psi_2^{-1}\colon&([u_0:u_1],[v_0:v_1])\mapsto([u_1:u_0],[v_1:v_0]).
\end{align*}
Then
$\psi_2\alpha\psi_2^{-1},\psi_2\sigma\psi_2^{-1}\in\Aut_\k(\F_0)$
exchange the geometric components of $q'$ and exchange or preserve
the rulings of $\F_0$, hence lift to elements of $\Aut_\k(X')$. The
birational involution $\psi_2\beta\psi_2^{-1}$ preserves the first
ruling of $\F_0$ and exchanges its sections through the components
of $q'$, and it contracts the fibre above
$\{[1:\zeta],[1:\zeta^2]\}$ onto $q'$, so it lifts to an
automorphism of $X'$. So $\varphi_2\colon X\rat X'$ is an
$\Aut_\k(X)$-equivariant link of type II.
\end{enumerate}
\end{ex}

\begin{lem}\label{lem:links6}
Let $|\k|=2$ and let $X$ be the del Pezzo surface of degree $6$ from
Theorem~\ref{thm:1}$(\ref{1:51})$. Any $\Aut_\k(X)$-equivariant link
of type II starting from $X$ is one of the links
$\varphi_1,\varphi_2$ in Example~\ref{ex:links6}, up to
automorphisms of the target surface.
\end{lem}
\begin{proof}
Let $\varphi$ be an $\Aut_\k(X)$-equivariant link of type II
starting from  $X$ and let $\eta\colon Y\to X$ be the blow-up of its
base-locus. Then $Y\to\ast$ is an $\Aut_\k(X)$-equivariant rank $2$
fibration, and by Remark~\ref{rmk:links-opt} the orbit blown-up by
$\eta$ has $\leq 5$ components. Since $\rk\,\NS(Y)^{\Aut_\k(X)}=2$,
there are exactly two extremal $\Aut_\k(X)$-equivariant contractions
starting from $Y$, namely the birational morphisms $\eta$ and
$\eta'$. It follows that the orbit blown up by $\eta$ determines
$\varphi$ up to automorphisms of $X'$. By
Lemma~\ref{lem:orbit51}(\ref{orbit51:3}), the only
$\Aut_\k(X)$-orbits on $X$ are $p:=([1:1:1],[1:1:1])$ and
$q:=\{([1:\zeta:\zeta^2],[1:\zeta^2:\zeta]),([1:\zeta^2:\zeta],[1:\zeta:\zeta^2])\}$,
$\zeta\notin\k$, $\zeta^3=1$. The birational maps $\varphi_1\colon
X\rat\F_0$ and $\varphi_2\colon X\rat X'$ in Example~\ref{ex:links6}
are $\Aut_\k(X)$-equivariant links of type II with base-points $p$
and $q$, respectively.
\end{proof}

\begin{prop}\label{pro:links-list}
Let $X$ be the del  Pezzo surface of degree $6$ from
Theorem~\ref{thm:1}$(\ref{1:51})$.
    \begin{enumerate}
    \item\label{links-list:21}  If $|\k|\geq3$, there is no $\Aut_\k(X)$-equivariant  link starting from $X$.
    \item\label{links-list:22}  If $|\k|=2$, any $\Aut_\k(X)$-equivariant link starting from $X$ is one of the $\Aut_\k(X)$-equivariant links of type II in Example~\ref{ex:links6}, up to automorphisms of the target surface.
    \end{enumerate}
\end{prop}
\begin{proof}
Since $\rk\,\NS(X)^{\Aut_\k(X)}=1$, the only $\Aut_\k(X)$-equivariant
links starting from $X$ are of type I or II, and by
Remark~\ref{rmk:links-opt}, they are not defined in an
$\Aut_\k(X)$-orbit with $\leq5$ geometric components and the blow-up
of this orbit is a del Pezzo surface.

If $|\k|\geq4$, no such orbits exist by
Lemma~\ref{lem:orbit51}(\ref{orbit51:2}). If $|\k|=3$, the blow-up
of any $\Aut_\k(X)$-orbit $X$ with $\leq5$ geometric components is
not a del Pezzo surface by Lemma~\ref{lem:orbit51}(\ref{orbit51:4}).

If $|\k|=2$, Lemma~\ref{lem:rk1-2} implies that the blow-up of any
$\Aut_\k(X)$-orbit on $X$ with $\leq5$ geometric components does not
admit an $\Aut_\k(X)$-equivariant conic fibration. In particular,
there is no $\Aut_\k(X)$-equivariant link of type I starting from
$X$. By Lemma~\ref{lem:links6}, any $\Aut_\k(X)$-equivariant link of
type II starting from $X$ is one of the birational maps in
Example~\ref{ex:links6}.
\end{proof}

\subsection{$\Aut_\k(X,\pi)$-equivariant links of conic fibrations}
We compute all $\Aut_\k(X,\pi)$-equivariant links starting from the
conic fibrations listed in Theorem~\ref{thm:1}.

\begin{lem}\label{lem:linkscbF}
Let $\pi\colon X\to\p^1$ be a conic fibration from
Theorem~\ref{thm:1}$(\ref{1:61})$ such that $\k^*/\mu_n(\k)$ is
trivial. Let $\pi'\colon Y\to\p^1$ be a conic fibration such that
$\Aut(Y/\pi')$ is infinite. Suppose that there is a
$\Aut_\k(X,\pi)$-equivariant link $\psi\colon X\rat Y$ of type II.
Then $Y\simeq X$.
\end{lem}
\begin{proof}
The link $\psi$ preserves the set of singular fibres, of which there
are at least $4$, and it commutes with the $\Gal(\bk/\k)$-action on
the set of geometric components of the singular fibres. It follows
from Lemma~\ref{lem:min-cb} that $Y$ is obtained by blowing up a
Hirzebruch surface. Since $Y$ is an $\Aut_\k(X,\pi)$-Mori fibre
space by definition of an equivariant link, the subgroup
$\Aut_\k(X,\pi)\subseteq\Aut_\k(Y,\pi')$ contains an element
exchanging the components of a singular geometric fibre. Moreover,
since $\Aut(Y/\pi)$ is infinite by hypothesis, Lemma~\ref{lem:cb-F}
implies that there is a birational morphism $\eta'\colon Y\to\F_m$
blowing up points $q_1,\dots,q_s\in S_m$ such that
$\sum_{i=1}^s\deg(q_i)=2m$. By Lemma~\ref{prop:cb-F}(\ref{cb-F:2})
and since $\k^*/\mu_n(\k)$ is trivial, we have
$\Aut_\k(X/\pi)=\langle\varphi\rangle\simeq\Z/2$ for some involution
$\varphi$. By Lemma~\ref{prop:cb-F}(\ref{cb-F:4}) it has a fixed
curve in $X$, which is the strict transform $C$ of a hyperelliptic
curve $C'$ in $\F_n$ (the irreducible double cover of $\p^1$)
ramified at $p_1,\dots,p_s$ and disjoint from $S_{-n}$.  It follows
that $C'\sim 2S_{-n}+2nf=2S_n$ and hence $C^2=-4n$ since the strict
transform of $S_n$ is a $(-n)$-curve on $X$.
 An $\Aut_\k(X,\pi)$-orbit contains either $1$ or $2$ points in the same fibre.
The base-points of the $\Aut_\k(X,\pi)$-equivariant link $\psi$ are
therefore necessarily contained in the $\Aut_\k(X,\pi)$-fixed curve
$C$. Since $C$ is a double cover of $\p^1$, it follows that
$C^2=\psi(C)^2$. The map $\psi\varphi\psi^{-1}\in\Aut_\k(Y/\pi')$
exchanges the components of each singular fibre, so it also
exchanges the two special sections of $Y$. By
Lemma~\ref{prop:cb-F}(\ref{cb-F:4}) it fixes a curve $D\subset Y$,
which satisfies $D^2=-4m$ with the same argument as above. It
follows that $C=\psi^{-1}(D)$, and now $-4n=C^2=D^2=-4m$ implies
$n=m$. Since $\psi$ induces the identity on $\p^1$, we conclude that
$\{q_1,\dots,q_s\}=\{p_1,\dots,p_r\}$.
\end{proof}

\begin{lem}\label{lem:linkscbFS}
Suppose that $\pi\colon X\to\p^1$ is a conic fibration as in
Theorem~\ref{thm:1}$(\ref{1:4})$ or $(\ref{1:6})$. Then there are no
$\Aut_\k(X,\pi)$-equivariant links of type $I$, $III$ and $IV$
starting from $X$. Moreover,
\begin{enumerate}
\item\label{linkscbFS:2} if $X=\F_n$, $n\geq2$, there are no $\Aut_\k(\F_n,\pi_n)$-equivariant links of type II starting from $\F_n$.
\item\label{linkscbFS:3} If $X$ is as in Theorem~\ref{thm:1}$(\ref{1:61})$ and $\k^*/\mu_n(\k)$ is non-trivial, there are no $\Aut_\k(X,\pi)$-equivariant links of type II starting from $X$.
\item\label{linkscbFS:4} If $X$ is as in Theorem~\ref{thm:1}$(\ref{1:62})$, there are no $\Aut_\k(X,\pi)$-equivariant links of type II starting from $X$.
\end{enumerate}
\end{lem}
\begin{proof}
Since $\NS(X)^{\Aut_\k(X,\pi)}\simeq\Z^2$, no
$\Aut_\k(X,\pi)$-equivariant links of type I can start from $X$. An
$\Aut_\k(X,\pi)$-link of type III can only start from a del Pezzo
surface (see Remark~\ref{rmk:links-opt}), so not from $X$. Since
$\Aut_\k(X,\pi)=\Aut_\k(X)$, any automorphism of $X$ preserves the
conic bundle structure, so there are no $\Aut_\k(X,\pi)$-equivariant
links of type IV starting from $X$.

(\ref{linkscbFS:2}) Suppose that there is a
$\Aut_\k(\F_n)$-equivariant link $\psi\colon\F_n\rat Y$ of type II,
and let $B\subset\F_n$ be the orbit of base-points and $d\geq1$ its
number of geometric components. We have
$|\Aut_\k(\F_n/\pi_n)|=|\k^{n+1}|\geq2^3$ by Remark~\ref{rmk:cb-F},
so the $\Aut_\k(\F_n/\pi_n)$-orbit of any point outside the special
section has at least two geometric components in the same geometric
fibre. If follows that $B\subset S_{-n}$ and hence $\psi$ is a
birational map from $\F_n$ to $\F_{n+d}$ and sends $S_{-n}$ onto
$S_{-(n+d)}$. Let $P\in\k[z_0,z_1]_d$ be a homogeneous polynomial
defining $B$. Then $\psi$ is of the form
\[
\psi\colon\F_n\rat\F_{n+d},\ [y_0:y_1;z_0:z_1]\dashmapsto
[Q(z_0,z_1)y_0:R(z_0,z_1)y_0+P(z_0,z_1)y_1;z_0:z_1]
\]
for some homogeneous $Q,R\in\k[z_0,z_1]$ of degree $d$. For any
$\alpha\in\Aut_\k(\F_n/\pi_n)\simeq\k[z_0,z_1]_n$ we have
$\psi\alpha\psi^{-1}\in\Aut_\k(\F_{n+d}/\pi_{n+d})$, and we compute
that it implies $\lambda:=\frac{P}{Q}\in\k^*$ and hence
$\lambda\alpha\in\k[z_0,z_1]_{n+d}$ (see Remark~\ref{rmk:cb-F}),
contradicting $d\geq1$.

(\ref{linkscbFS:3}) If $\pi\colon X\to\p^1$ is a conic fibration as
in Theorem~\ref{thm:1}(\ref{1:61}) and the torus subgroup
$\k^*/\mu_n(\k)\subset\Aut_\k(X/\pi)$ is non-trivial, then the
$\Aut_\k(X/\pi)$-orbit of a point on a smooth fibre outside the two
$(-n)$-sections has at least two geometric components in the same
smooth fibre. Since $\Z/2\subset\Aut_\k(X/\pi)$ exchanges the two
$(-n)$-sections, the same holds for any point contained in them. It
follows that there are no $\Aut_\k(X,\pi)$-equivariant links of type
II starting from $X$.

(\ref{linkscbFS:4}) Let $\pi\colon X\to\p^1$ be a conic fibration as
in Theorem~\ref{thm:1}(\ref{1:62}). Consider the subgroup
$\SlO^{L,L'}(\k)$ of $\Aut_\k(X/\pi)$ fixing the geometric
components of the special double section $E$ from
Lemma~\ref{prop:cb1}(\ref{cb1:2}). Let us show that
$|\SlO^{L,L'}(\k)|\geq2$. From Lemma~\ref{lem:min-cb-sing1} we
obtain:
\begin{itemize}
\item If $L,L'$ are not $\k$-isomorphic, then $\k$ is infinite, and so $\SlO^{L,L'}(\k)\simeq\k^*$ is infinite.
\item If $L=L'$, then
$\SlO^{L,L}(\k)\simeq\{\alpha\in L^*\mid \alpha\alpha^g=1\}$, where
$g$ is the generator of $\Gal(L/\k)$. If $|\k|\geq3$, then $\pm1\in
\SlO^{L,L}(\k)$, and if $|\k|=2$, then $|\SlO^{L,L}(\k)|=|L^*|=3$.
\end{itemize}
In any case, it follows that the $\Aut_\k(X/\pi)$-orbit of a point
on a smooth fibre outside $E$ has at least two geometric components
in the same smooth fibre. Since $\Aut_\k(X/\pi)$ contains an
involution exchanging the geometric components of  $E$ by
Lemma~\ref{prop:cb1}(\ref{cb1:2}), the same holds for any point in
$E$. It follows that there are no $\Aut_\k(X,\pi)$-equivariant links
of type II starting from $X$.
\end{proof}

\subsection{Proof of Theorem~\ref{thm:2}, Corollary~\ref{cor:1} and Theorem~\ref{thm:3}}

Let $G$ be an affine algebraic group and let $X/B$ be a $G$-Mori fibre space
that is also a $G(\k)$-Mori fibre space. A $G$-equivariant
birational map is in particular $G(\k)$-equivariant, hence if $X$ is
$G(\k)$-birationally (super)rigid it is also $G$-birationally
(super)rigid.

On the other hand, $G$-birationally (super)rigid does not imply
$G(\k)$-birationally (super)rigid: the next lemma shows that the del
Pezzo surface $X$ of degree $6$ obtained by blowing up $\p^2$ in
three rational points is $\Aut(X)$-birationally superrigid and
Example~\ref{ex:links6} shows that $X$ is not even
$\Aut_\k(X)$-birationally rigid if $|\k|=2$.

\begin{lem}\label{lem:DP6-Autsuper}
Any del Pezzo surface $X$ of degree $6$ is $\Aut(X)$-birationally
superrigid.

\end{lem}
\begin{proof}
The surface $X_{\bk}$ is isomorphic to the del Pezzo surface
obtained by blowing up three rational points in $\p^2_{\bk}$. In
particular, $\rk\,\NS(X_{\bk})^{\Aut_{\bk}(X)}=1$ by
Lemma~\ref{prop:DP1}(\ref{DP1:3}), hence $X$ is an $\Aut(X)$-Mori
fibre space and there are no $\Aut(X)$-equivariant links of type III
or IV starting from $X$. The base-locus of an $\Aut(X)$-equivariant
link of type I or II is an $\Aut_{\bk}(X)\times\Gal(\bk/\k)$-orbit
on $X_{\bk}$, and by Remark~\ref{rmk:links-opt} it has $\leq5$
elements. Lemma~\ref{lem:orbit51}(\ref{orbit51:2}) implies that
$\Aut_{\bk}(X)=\Aut(X_{\bk})$ has no such orbits. By
Theorem~\ref{thm:sarkisov}, any $\Aut(X)$-equivariant birational map
starting from $X$ decomposes into isomorphisms and
$\Aut(X)$-equivariant links. As there are no $\Aut(X)$-equivariant
links starting from $X$, it follows that $X$ is
$\Aut(X)$-birationally superrigid.
\end{proof}

\begin{proof}[Proof of Theorem~\ref{thm:2}]
(\ref{2:2})--(\ref{2:4}) Any surface $X$ as in
Theorem~\ref{thm:1}(\ref{1:1})--(\ref{1:3}), (\ref{1:5c}), and (\ref{1:5b})
is a del Pezzo surface that is at the same time a $\Aut_\k(X)$-Mori
fibre space and an $\Aut(X)$-Mori fibre space. Any conic fibration
$\pi\colon X\to\p^1$ as in
Theorem~\ref{thm:1}(\ref{1:4}) and (\ref{1:6}) has $\Aut(X)=\Aut(X,\pi)$
and $\Aut_\k(X)=\Aut_\k(X,\pi)$, and it is at the same time a
$\Aut_\k(X)$-Mori fibre space and an $\Aut(X)$-Mori fibre space. By
Theorem~\ref{thm:sarkisov}, any equivariant birational map between
equivariant Mori fibre spaces decomposes into equivariant Sarkisov
links, hence in order to show that an equivariant Mori fibre space
$X/B$ is equivariantly birationally superrigid, it suffices to show
that there are no equivariant links starting from $X$.

(\ref{2:2}) For $X=\p^2$, $X=\Ql^L$ and $X=\F_0$ the claim follows
from Lemma~\ref{lem:linksPQF} and for $X=\F_n$, $n\geq2$, from
Lemma~\ref{lem:linkscbFS}(\ref{linkscbFS:2}). For $X$ a del Pezzo
surface of degree $6$ as in (\ref{1:52})--(\ref{1:54}) the claim
follows from Proposition~\ref{pro:links-list2-4}, and for a conic
fibration $X/\p^1$ as in (\ref{1:62}) from
Lemma~\ref{lem:linkscbFS}.

(\ref{2:5}) For $X$ a del Pezzo surface of degree $6$ as in
(\ref{1:5c}) the claim is Proposition~\ref{prop:5c}.

(\ref{2:3}) The claim follows from Proposition~\ref{pro:links-list}.

(\ref{2:4}) The claim follows from Lemma~\ref{lem:linkscbF} and
Lemma~\ref{lem:linkscbFS}(\ref{linkscbFS:3}).

(\ref{2:1}) It follows from (\ref{2:2})--(\ref{2:4})  that for any
surface $X$ in Theorem~\ref{thm:1} there is an algebraic extension
$L/\k$ such that $X_L$ is $\Aut_L(X)$-birationally superrigid.
Therefore, $X$ is also $\Aut(X)$-birationally superrigid.
\end{proof}

\begin{proof}[Proof of Corollary~\ref{cor:1}]
Theorem~\ref{thm:1} implies (\ref{cor1:1}). By
Theorem~\ref{thm:2}(\ref{2:1}), the surfaces $X$ in
Theorem~\ref{thm:1} are $\Aut(X)$-birationally superrigid, so the
groups $\Aut(X)$ are maximal and they are conjugate if and only if
their surfaces are isomorphic. Theorem~\ref{thm:1} now implies
(\ref{cor1:3}).

By Theorem~\ref{thm:2}(\ref{2:2})--(\ref{2:4}), the surfaces $X$
from
Theorem~\ref{thm:1}(\ref{1:1})--(\ref{1:4}) and (\ref{1:52})--(\ref{1:54}),
(\ref{1:62}) are $\Aut_\k(X)$-birationally superrigid.  The surface
$X$ from (\ref{1:61}) are $\Aut_\k(X)$-birationally rigid within the
set of classes of surfaces from Theorem~\ref{thm:1}.
 The del Pezzo surfaces $X$ from (\ref{1:5c}) and (\ref{1:51}) are $\Aut_\k(X)$-birationally superrigid if $|\k|\geq3$.
Hence the listed groups $\Aut_\k(X)$ are maximal and they are
conjugate by a birational map if and only if their surfaces are
isomorphic. Theorem~\ref{thm:1} now implies (\ref{cor1:4}).
\end{proof}

\begin{lem}\label{rmk:galoisextensions}
  Let $\k$ be a perfect field and let $F/\k$ be a field extension. The following are equivalent:
  \begin{enumerate}
    \item\label{item:point} There exists a point $p$ of degree $3$ in $\p^2$, not all irreducible components collinear, such that $F$ is the splitting field of $p$.
    \item \label{item:polynomial} $F$ is the splitting field of an irreducible polynomial of degree $3$ over $\k$.
    \item\label{item:Galois} The field extension $F/\k$ is Galois and $\Gal(F/\k)$ is isomorphic to a transitive subgroup of $\sym_3$ (that is to $\Z/3\Z$ or $\sym_3$).
  \end{enumerate}
\end{lem}

\begin{proof}
  (\ref{item:point}) implies (\ref{item:polynomial}): Since the irreducible components $p_i$ of $p$ are not collinear, there is an irreducible conic defined over $\k$ that contains $p$. With a linear transformation defined over $\k$ this conic can be assumed to be given by $x^2-yz=0$, and so $p_i=[a_i:a_i^2:1]$ for some $a_i\in F$ for $i=1,2,3$, and $\{a_1,a_2,a_3\}$ is a Galois orbit.
  Hence $q(t)=(t-a_1)(t-a_2)(t-a_3)\in\k[t]$ is irreducible.
  The splitting field $L$ of $q(t)$ is $\k(a_1,a_2,a_3)=F$.

  (\ref{item:polynomial}) implies (\ref{item:point}): Similar to above.

  (\ref{item:polynomial}) implies (\ref{item:Galois}):  By assumption $F$ is the splitting field of an irreducible and hence separable polynomial $f$. Therefore, $F/\k$ is normal and hence Galois.
  So $\Gal(F/\k)$ acts transitively on the three roots of $f$, hence $\Gal(F/\k)$ is isomorphic to a transitive subgroup of $\sym_3$.

  (\ref{item:Galois}) implies (\ref{item:polynomial}):
  Note that by the Primitive element Theorem, there exists $a\in F$ such that $F=\k(a)$. Let $f$ be the minimal polynomial of $a$ over $\k$, hence $\deg(f)=[F:\k]=|\Gal(F/\k)|\in\{3,6\}$.
  Let $L$ be the splitting field of $f$, which is a normal extension of $\k$. In particular, $F=\k(a)=L$.
  Hence, if $\deg(f)=3$ we are done.

  In the other case we have $\Gal(F/\k)\simeq \sym_3$, so $\deg(f)=6$. The roots of $f$ form one Galois-orbit.
  After fixing an isomorphism $\Gal(F/\k)\simeq \sym_3$, we write $\sigma_{ij}=(ij)$, and we write $\tau=(123)$.
  So we can write the six roots of $f$ as $a_i=\tau^i(a)$ for $i=1,2,3$ (so $a_3=a$), and $a_4=\sigma_{13}(a)$, $a_5=\sigma_{23}(a)$, $a_6=\sigma_{12}(a)$.
  Set
  \[b_1= a_1a_4, b_2 = a_2a_5, b_3=a_3a_6\]
  and note that the $\sigma_{ij}$ act as transposition of $b_i, b_j$, and that that $\tau$ is the translation $b_1\mapsto b_2\mapsto b_3$.
  So $\{b_1,b_2,b_3\}$ is a $\Gal(F/\k)$-orbit of size $3$ with minimal polynomial $g=(t-b_1)(t-b_2)(t-b_3)\in\k[t]$.
  So the splitting field $L'$ of $g$ is contained in $F$ and its Galois group is isomorphic to $\sym_3$. Hence
  \[6=|\Gal(L'/\k)|=[L':\k]\leq[F:\k]=6,\]
  which implies $F=L'$ is the splitting field of an irreducible polynomial of degree $3$.
\end{proof}

\begin{proof}[Proof of Theorem~\ref{thm:3}]
By Corollary~\ref{cor:1}(\ref{cor1:4}) it suffices to list the
isomorphism classes of the surfaces in
Theorem~\ref{thm:1}(\ref{1:1})--(\ref{1:4}),
(\ref{1:52})--(\ref{1:54}), (\ref{1:6}), and for (\ref{1:5c}) and
(\ref{1:51}) if $|\k|\geq3$.

The plane $\p^2$ is unique up to isomorphism by Ch\^{a}telet's
Theorem, $\F_0$ is unique up to isomorphism by
Lemma~\ref{prop:Q}(\ref{Q:1}), and for any $\k$-isomorphism class of
quadratic extensions $L/\k$ we have a unique isomorphism class of
$\Ql^L$, also by Lemma~\ref{prop:Q}(\ref{Q:1}). Hirzebruch surfaces
are determined by their special section. The parametrisation of the
classes of del Pezzo surfaces from (\ref{1:5c}) follows from
Lemma~\ref{prop:DP6}(\ref{DP6:3}),
Lemma~\ref{prop:DP8}(\ref{DP8:3}) and Lemma~\ref{rmk:galoisextensions}. The parametrisation of the
classes of del Pezzo surfaces from (\ref{1:5b}) follows from
Lemma~\ref{prop:DP1}(\ref{DP1:1}),
Lemma~\ref{prop:DP5}(\ref{DP5:4}), Lemma~\ref{prop:DP7}(\ref{DP7:4}), Lemma~\ref{prop:DP4}(\ref{DP4:4}) and Lemma~\ref{rmk:galoisextensions}. The parametrisations for the
conic fibrations from (\ref{1:61}) and (\ref{1:62}) follow from
Lemma~\ref{lem:min-cb-singF} and Lemma~\ref{cor:min-cb-sing}.
\end{proof}

\section{The image by a quotient homomorphism}\label{s:hom}

We call two Mori fibre spaces $X_1/\p^1$ and $X_2/\p^1$ equivalent
if there is a birational map $X_1\rat X_2$ that preserves the
fibrations. In particular, if $\varphi\colon X_1\rat X_2$ is a link
of type II between Mori fibre spaces $X_1/\p^1$ and $X_2/\p^1$, then
these two are equivalent. There is only one class of Mori fibre
spaces birational to $\F_1$ \cite[Lemma]{Schneider}, because all
rational points in $\p^2$ are equivalent up to $\Aut(\p^2)$. We
denote by $J_6$ the set of classes of Mori fibre spaces birational
to some $\Sl^{L,L'}$, and by $J_5$ the set of classes birational to a
blow-up of $\p^2$ in a point of degree $4$ whose geometric
components are in general position. We call two Sarkisov links
$\varphi$ and $\varphi'$ of type II between conic fibrations
equivalent if the conic fibrations are equivalent and and if the
base-points of $\varphi$ and $\varphi'$ have the same degree. For a
class $C$ of equivalent rational Mori fibre spaces, we denote by
$M(C)$ the set of equivalence classes of links of type II between
conic fibrations in the class $C$ whose base-points have degree
$\geq16$.

\begin{proof}[Proof of Proposition~\ref{prop:juliaz}]
First, suppose that $[\bk:\k]=2$. Then every non-trivial algebraic
extension of $\k$ is $\bk$ by \cite[Satz 4]{ArtinSchreier} and $\k$
is of characteristic zero \cite[p.231]{ArtinSchreier}. In
particular, $\p^2$ contains no points of degree $\geq3$, and so the
only rational Mori fibre spaces are Hirzebruch surfaces and
$\Sl^{\bk,\bk}\to\p^1$. Moreover, $M(\F_1)$ is empty. By
\cite[Theorem 1.3]{Zimmermannb}, there is a surjective homomorphism
$\Bir_\R(\p^2)\to\bigoplus_I\Z/2$, where $|I|=|\R|$. In fact, by
construction of the homomorphism, there is a natural bijection $I\to
\{\frac{|a|}{a^2+b^2}\mid a,b\in\R,b\neq0\}$. The whole article
\cite{Zimmermann} can be translated word-by-word over a field $\k$
with $[\bk:\k]=2$, and consequently we have a surjective
homomorphism $\Bir_\k(\p^2)\to\bigoplus_I\Z/2$, where
$I=\{\frac{a^2}{a^2+b^2}\mid a,b\in\k,b\neq0\}$ (we replace $|a|$ by
$a^2$), and $I$ has the cardinality of $\k$. If $[\bk:\k]>2$, the
result is \cite[Theorem 3, Theorem 4.]{Schneider}.
\end{proof}

\begin{defi}\label{def:psi}
Let $\mathrm{BirMori}(\p^2)$ be the groupoid of birational maps
between Mori fibre spaces birational to $\p^2$. It is generated by
Sarkisov links by Theorem~\ref{thm:sarkisov}. The homomorphism
$\tilde{\Psi}$ of groupoids from \cite[Theorem 3, Theorem
4]{Schneider}
\[
\begin{tikzcd}[link2]
\mathrm{BirMori(\p'^2)}\ar[r,"\tilde{\Psi}"] & (\bigoplus_{\chi\in M(\F_1)}\Z/2)\bigast_{C\in J_5}(\bigoplus_{\chi\in M(C)}\Z/2)\ast(\bigast_{C\in J_6}\bigoplus_{\chi\in M(C)}\Z/2)\\[-8pt]
\rotatebox{90}{$\subseteq$}&\\[-7pt]
\Bir_\k(\p^2)\ar[uur,"\Psi",swap] &
\end{tikzcd}
\]
sends any Sarkisov link of type II between conic fibrations and
whose base-point has degree $\geq16$ onto the generator indexed by
its class, and it sends all other Sarkisov links and all
isomorphisms between Mori fibre spaces to zero.
\end{defi}

\begin{rmk}\label{rmk:psi-nontrivial}
The homomorphism $\Psi$ is non-trivial. Indeed, the surjective
homomorphism $\Bir_\k(\p^2)\to
(\bigoplus_{I_0}\Z/2)\ast(\bigast_{J_5}\bigoplus_{I}\Z/2)\ast(\bigast_{J_6}\bigoplus_{I}\Z/2)$
from \cite[Theorem 4]{Schneider} is obtained by composing $\Psi$
with suitable projections within each abelian factor of the free
product, see \cite[Proof of Theorem 4 in \S6]{Schneider}.
\end{rmk}

We now compute the images by $\Psi$ of $\k$-points of the maximal
algebraic subgroups of $\Bir_\k(\p^2)$ listed in
Theorem~\ref{thm:1}.

\begin{rmk}\label{rmk:homomorphism}
By definition of the groupoid homomorphism $\tilde\Psi$
(Definition~\ref{def:psi}), it maps automorphism groups of Mori
fibre spaces onto zero, so the groups $\Psi(\Aut_\k(\p^2))$,
$\tilde{\Psi}(\Aut_\k(\Ql^{L}))$, $\tilde{\Psi}(\Aut_\k(\F_n))$,
$n\neq1$, and $\tilde{\Psi}(\Aut(\Sl^{L,L'},\pi))$ are trivial. A
del Pezzo surface $X$ of degree $6$ as in
Theorem~\ref{thm:1}(\ref{1:5c}) is a Mori fibre space by
Lemma~\ref{prop:DP6} and Lemma~\ref{prop:DP8}, so
$\tilde\Psi(\Aut_\k(X))$ is trivial as well.

If $X$ is a del Pezzo surface from Theorem~\ref{1:1}(\ref{1:5a}),
there exists a birational morphism $\eta\colon X\to \Ql^L$ such that
$\eta\Aut_\k(X)\eta^{-1}\subset\Aut_\k(\Ql^{L})$, so in particular
$\tilde\Psi(\eta\Aut_\k(X)\eta^{-1})$ is trivial as well.
\end{rmk}

\begin{lem}\label{lem:links-dP-hom}
Let $X$ be a del Pezzo surface of degree $6$ from
Theorem~\ref{thm:1}$(\ref{1:5b})$, which is equipped with a
birational morphism $\eta\colon X\to Y$ to $Y=\p^2$ or $Y=\F_0$.
Then $\tilde{\Psi}(\eta\Aut_\k(X)\eta^{-1})$ is trivial.
\end{lem}
\begin{proof}
Let $X$ be a del Pezzo surface of degree $6$ from
Theorem~\ref{thm:1}(\ref{1:51}), (\ref{1:53}), and (\ref{1:54}), which is
the blow-up $\eta\colon X\to \p^2$ in three rational points or in a
point of degree $3$. By Lemma~\ref{prop:DP1}(\ref{DP1:2}),
Lemma~\ref{prop:DP5}(\ref{DP5:2}) and
Lemma~\ref{prop:DP7}(\ref{DP7:2}), the group
$\eta\Aut_\k(X)\eta^{-1}$ is generated by subgroups of
$\Aut_\k(\p^2)$ and a quadratic involution of $\p^2$ that has either
three rational base-points or is a Sarkisov link of type II with a
base-point of degree $3$. It follows from the definition of
$\tilde{\Psi}$ (Definition~\ref{def:psi}) that
$\tilde{\Psi}(\eta\Aut_\k(X)\eta^{-1})$ is trivial.

The del Pezzo surface  $X$ of degree $6$ from
Theorem~\ref{thm:1}(\ref{1:52}) is the blow-up of $\eta\colon
X\to\F_0$ in a point of degree $2$. By
Lemma~\ref{prop:DP4}(\ref{DP4:2}), the group
$\eta\Aut_\k(X)\eta^{-1}$ is generated by subgroups of
$\Aut_\k(\F_0)$ and a birational involution of $\F_0$ that is a link
of type II of conic fibrations with a base-point of degree $2$.
Again it follows that $\tilde{\Psi}(\eta\Aut_\k(X)\eta^{-1})$ is
trivial.
\end{proof}

\begin{lem}\label{lem:links-F-hom}
Let $n\geq2$ and let $\varphi\colon\F_n\rat\F_n$ be the involution
from Example~\ref{ex:gen-F} with base-points $p_1,\dots,p_r\in\F_n$.
Then there exist links $\varphi_1,\dots,\varphi_r$ of type II
between Hirzebruch surfaces such that $\varphi_i$ has a base-point
of degree $\deg(p_i)$ and $ \varphi=\varphi_r\cdots\varphi_1$.
\end{lem}
\begin{proof}
Recall from Example~\ref{ex:gen-F} that $p_1,\dots,p_r$ are
contained in the section $S_n\subset\F_n$ and that the homogeneous polynomials
$P_i\in\k[z_0,z_1]_{\deg(p_i)}$ define $\pi(p_i)$. The involution $\varphi$ is given by
\[\
\varphi\colon(y_1,z_1)\rat(\nicefrac{P(z_1)}{y_1},z_1)
\]
We define $d_0:=0$ and $d_i:=\sum_{j=1}^i\deg(p_j)$. For
$i=1,\dots,r$, the birational maps
\begin{align*}
\varphi_i\colon\F_{n-d_{i-1}}\rat\F_{n-d_i},\
&(y_1,z_1)\dashmapsto(\nicefrac{y_1}{P_i(z_1)},z_1)\quad d_i\leq n,\\
\varphi_i\colon\F_{n-d_{i-1}}\rat\F_{d_{i}-n},\
&(y_1,z_1)\dashmapsto(\nicefrac{P_i(z_1)}{y_1},z_1)\quad d_{i-1}\leq n,d_i>n\\
\varphi_i\colon\F_{d_{i-1}-n}\rat\F_{d_i-n},\
&(y_1,z_1)\dashmapsto(P_i(z_1)y_1,z_1),\quad d_{i-1}>n
\end{align*}
are links of type II with a base-point of degree $\deg(p_i)$, and we
compute that $\varphi=\varphi_r\cdots\varphi_1$.
\end{proof}

\begin{lem}\label{cor:links-F-hom}
Let $\pi\colon X\to\p^1$ be a conic fibration from
Theorem~\ref{thm:1}(\ref{1:61}) and let $\eta\colon X\to \F_n$,
$n\geq2$, be the birational morphism blowing up $p_1,\dots,p_r$. Let
$\varphi\colon\F_n\rat\F_n$ be the involution from
Example~\ref{ex:gen-F} and $\varphi=\varphi_r\cdots\varphi_1$ the
decomposition into links of type II from
Lemma~\ref{lem:links-F-hom}. Then
$\tilde\Psi(\eta\Aut_\k(X,\pi)\eta^{-1})$ is generated by the
element
$\tilde\Psi(\varphi)=\tilde\Psi(\varphi_r)+\cdots+\tilde\Psi(\varphi_1)$.
\end{lem}
\begin{proof}
Let $\Delta\subset\p^1$ be the image of the singular fibres of $X$.
By Lemma~\ref{prop:cb-F}(\ref{cb-F:1}--\ref{cb-F:2}), we have
\[\Aut_\k(X,\pi)\simeq\Aut_\k(X/\pi)\rtimes \Aut_\k(\p^1,\Delta)\quad\text{and}\quad \Aut_\k(X/\pi)\simeq H\rtimes\langle\eta^{-1}\varphi\eta\rangle\]
where $\eta H\eta^{-1}\subset\Aut_\k(\F_n)$. Moreover, any
$\alpha\in\Aut_\k(\p^1,\Delta)$ lifts to an element
$\tilde\alpha\in\Aut_\k(\F_n,p_1,\dots,p_r)$, which lifts via $\eta$
to an element of $\Aut_\k(X,\pi)$. It follows from the definition of
$\tilde{\Psi}$ that
$\tilde{\Psi}(\eta\Aut_\k(\p^1,\Delta)\eta^{-1})$ and
$\tilde{\Psi}(\eta H\eta^{-1})$ are trivial, and that
$\tilde{\Psi}(\eta\Aut_\k(X,\pi)\eta^{-1})$ is generated by
$\tilde{\Psi}(\varphi)=\tilde{\Psi}(\varphi_r)+\cdots+\tilde{\Psi}(\varphi_1)$.
\end{proof}

\begin{lem}\label{lem:links-cb1-hom}
Let $\varphi\colon\Sl^{L,L'}\rat\Sl^{L,L'}$ be the involution from
Example~\ref{ex:gen-cb1} with base-points
$p_1,\dots,p_r\in\Sl^{L,L'}$. Then there exist links
$\varphi_1,\dots,\varphi_r\colon \Sl^{L,L'}\rat\Sl^{L,L'}$ of type
II over $\p^1$ and $\alpha\in\Aut_\k(\Sl^{L,L'}/\pi)$ such that
$\varphi_i$ has base-point $p_i$ and such that $
\varphi=\alpha\varphi_r\cdots\varphi_1$.
\end{lem}
\begin{proof}
It suffices to construct the $\varphi_i$ for the involution
$\varphi$ in the case that $L=L'$, since the involution for the
other case is obtained by conjugating $\varphi$ with a suitable
element of $\gamma\in\PGL_2(\bk)\times\PGL_2(\bk)$, see
Example~\ref{ex:gen-cb1}. Let $E_1,E_2$ be the geometric components
of the unique irreducible curve contracted by any birational
contraction $\eta\colon\Sl^{L,L'}\to \Ql^L$. For $i=1,\dots,r$, let
$T_{i1},T_{i2}\in L[x,y]$ be the homogeneous polynomials defining
the fibres through the geometric components of the $p_i$ contained
in $E_1,E_2$, respectively. Let $P_1:=T_{11}\cdots T_{r1}$ and
$P_2:=T_{12}\cdots T_{r2}$. Recall from Example~\ref{ex:gen-cb1}
that $\psi:=\eta\phi\eta^{-1}$ is of the form
\[
\psi\colon([u_0:u_1],[v_0:v_1])\dashmapsto([v_0P_1(u_0v_0,u_1v_1) :
v_1P_2(u_0v_0,u_1v_1)],[u_0P_2(u_0v_0,u_1v_1) :
u_1P_1(u_0v_0,u_1v_1)])
\]
For $i=1,\dots,r$, define
\begin{multline*}
\psi_i\colon([u_0:u_1],[v_0:v_1])\dashmapsto([u_0T_{i2}(u_0v_0,u_1v_1)
: u_1T_{i1}(u_0v_0,u_1v_1)],\\
[v_0T_{i1}(u_0v_0,u_1v_1) : v_1T_{i2}(u_0v_0,u_1v_1)])
\end{multline*}
and let
\[
\tilde{\alpha}\colon([u_0:u_1],[v_0:v_1])\dashmapsto([v_0:v_1],[u_0:u_1]).
\]
Then $\alpha\psi_r\cdots\psi_1=\psi$. We take
$\varphi_i:=\eta^{-1}\psi_i\eta$ and
$\alpha:=\eta^{-1}\tilde{\alpha}\eta$.
\end{proof}

\begin{lem}\label{cor:links-cb1-hom}
Let $\pi\colon X\to\p^1$ be a conic fibration from
Theorem~\ref{thm:1}$(\ref{1:62})$ and let $\eta\colon
X\to\Sl^{L,L'}$ be the birational morphism blowing up
$p_1,\dots,p_r$. Let $\varphi\colon\Sl^{L,L'}\rat\Sl^{L,L'}$ be the
involution from Example~\ref{ex:gen-cb1} and
let $\varphi=\alpha\varphi_r\cdots\varphi_1$ be the decomposition into
links $\varphi_i$ of type II and an automorphism
$\alpha\in\Aut_\k(\Sl^{L,L'},\pi)$ from
Lemma~\ref{lem:links-cb1-hom}. Then
$\tilde{\Psi}(\eta\Aut_\k(X,\pi)\eta^{-1})$ is generated by the
element
$\tilde\Psi(\varphi)=\tilde{\Psi}(\varphi_r)+\cdots+\tilde{\Psi}(\varphi_1)$.
\end{lem}
\begin{proof}
Let $\Delta\subset\p^1$ be the image of the singular fibres of $X$.
By Proposition~\ref{prop:cb1}(\ref{cb1:1}--\ref{cb1:2}), we have
\[\Aut_\k(X,\pi)\simeq\Aut_\k(X/\pi)\rtimes ((D^{L,L'}_{\k}\rtimes\Z/2)\cap\Aut_\k(\p^1,\Delta)),\quad \Aut_\k(X/\pi)\simeq H\rtimes\langle\eta^{-1}\varphi\eta\rangle\]
where $\eta H\eta^{-1}\subset\Aut_\k(\Sl^{L,L'}/\pi)$. Moreover, any
element of $G:=D^{L,L'}_{\k}\rtimes\Z/2\cap\Aut_\k(\p^1,\Delta)$
lifts to an element of $\Aut_\k(\Sl^{L,L'},\pi)$, which lifts via
$\eta$ to an element of $\Aut_\k(X,\pi)$. It follows from the
definition of $\tilde{\Psi}$, that $\tilde{\Psi}(\eta G\eta^{-1})$,
$\tilde{\Psi}(\eta H\eta^{-1})$ and $\tilde\Psi(\alpha)$ are
trivial, and hence that $\tilde{\Psi}(\eta\Aut_\k(X,\pi)\eta^{-1})$
is generated by
$\tilde{\Psi}(\varphi)=\tilde{\Psi}(\varphi_r)+\cdots+\tilde{\Psi}(\varphi_1)$.
\end{proof}

\begin{proof}[Proof of Proposition~\ref{thm:4}]
Let $G$ be an infinite algebraic subgroup of $\Bir_\k(\p^2)$. By
Theorem~\ref{thm:1}, it is conjugate by a birational map to a
subgroup of $\Aut(X)$, where $X$ is one of the surfaces listed in
Theorem~\ref{thm:1}. We now compute
$\Psi(\theta\Aut_\k(X)\theta^{-1})$ for some birational map
$\theta\colon \p^2\rat X$. For any birational morphism $\eta\colon
X\to Y$ to a Mori fibre space $Y/B$, we have
\[
\Psi(\theta\Aut_\k(X)\theta^{-1})=\tilde{\Psi}(\theta^{-1}\eta^{-1})\tilde{\Psi}(\eta\Aut_\k(X)\eta^{-1})\tilde{\Psi}(\eta\theta).
\]
For the surfaces $X$ from
Theorem~\ref{thm:1}(\ref{1:1})--(\ref{1:5}), there exists such a
birational morphism $\eta$ such that
$\tilde{\Psi}(\eta\Aut_\k(X)\eta^{-1})$ is trivial by
Remark~\ref{rmk:homomorphism} and Lemma~\ref{lem:links-dP-hom}, and
hence $\Psi(\theta\Aut_\k(X)\theta^{-1})$ is trivial. Hence, if
$\Psi(G(\bk))$ is not trivial then $X$ is as in
Theorem~\ref{thm:1}(\ref{1:6}) and (\ref{3:2}) follows.

Let $X/\p^1$ be a conic fibration from
Theorem~\ref{thm:1}(\ref{1:6}), which is the blow-up $\eta\colon
X\to Y$ of points $p_1,\dots,p_r\in Y$ and $Y=\F_n$, $n\geq2$ or
$Y=\Sl^{L,L'}$. By Lemma~\ref{cor:links-F-hom} and
Lemma~\ref{cor:links-cb1-hom} the image
$\tilde{\Psi}(\eta\Aut_\k(X)\eta^{-1})$ is generated by the element
$\tilde{\Psi}(\varphi_r)+\cdots+\tilde{\Psi}(\varphi_1)$, where
$\varphi_i$ is a link of type II between conic fibrations in the
respective class and whose base-point is of degree $\deg(p_i)$. In
particular, since each factor of the free product is abelian, it
follows that $\Psi(\theta\Aut_\k(X)\theta^{-1})$ is generated by
$\tilde{\Psi}(\varphi_r)+\cdots+\tilde{\Psi}(\varphi_1)$.

By definition of $\tilde{\Psi}$ the image $\tilde\Psi(\varphi_i)$ is
non-trivial if and only if $\deg(p_i)\geq16$. Therefore, if
$\tilde\Psi(\varphi_r)+\cdots+\tilde\Psi(\varphi_1)$ is non-trivial,
it is the element indexed by the $i_1,\dots,i_s$ such that
$\deg(p_{i_k})\geq16$ and we infer that $|\{j\in
\{1,\ldots,r\}\mid\deg(p_j)=\deg(p_{i_k})\}|$ is odd for
$k=1,\dots,s$. This proves (\ref{3:3}). In particular,
$\Psi(G(\k))\simeq\Z/2\Z$.
\end{proof}



\pagebreak

\newpage
\appendix

\newcommand{\hiddensection}[1]{
    \refstepcounter{section}
    \subsection*{\Alph{section}.\hspace{1em}{#1}}
}

\begin{center}
{\LARGE \bf Corrigendum to: ``Algebraic subgroups \\ of the plane Cremona group over a perfect field''}
\end{center}

\hiddensection{Correction of the main theorems}

Below in Theorems~\ref{thm:1n} and~\ref{thm:3n} we correct case~\eqref{1:5c} of the classification from Theorems~\ref{thm:1} and~\ref{thm:3}.

\begin{thm}[\textit{cf.~}Theorem~\ref{thm:1}]\label{thm:1n}
Let $\k$ be a perfect field and $G$ an infinite algebraic subgroup
of\, $\Bir_\k(\p^2)$. Then there is a $\k$-birational map $\p^2\rat X$
that conjugates $G$ to a subgroup of $\Aut(X)$, with $X$ one of the
following surfaces, where no indication of the $\Gal(\bk/\k)$-action
means the canonical action.
\begin{enumerate}
\item\label{1:1n} $X=\p^2$ and $\Aut(\p^2)\simeq\PGL_3$.
\item\label{1:2n} $X=\F_0$ and $\Aut(\F_0)\simeq\Aut(\p^1)^2\rtimes\Z/2\simeq\PGL_2^2\rtimes\Z/2$.
\item\label{1:3n} $X=\Ql^L$ and  $\Aut(\Ql^L)$ is the $\k$-structure on $\Aut(\p^1_L)^2\rtimes\Z/2$ given by the $\Gal(L/\k)$-action $(A,B,\tau)^g=(B^g,A^g,\tau)$, where $L/\k$ is a quadratic extension.
\item\label{1:4n} $X=\F_n$, $n\geq2$, and the action of $\Aut(\F_n)$ on $\p^1$ induces a split exact sequence
\[
 1\to V_{n+1}\to \Aut(\F_n)\to \GL_2/\mu_n\to 1,
\]
where $\mu_n=\{a\id\mid a^n=1\}$ and $V_{n+1}$ is a vector space of dimension $n+1$.

\item\label{1:5n} $X$ is a del Pezzo surface of degree $6$ with $\NS(X_{\bk})^{\Aut_{\bk}(X)}=1$. The action of $\Aut_{\bk}(X)$ on $\NS(X_{\bk})$ induces the split exact sequence
\[1\rightarrow (\bk^*)^2\to \Aut_{\bk}(X)\to \sym_3\times\Z/2\rightarrow 1.\]
Moreover, we are in one of the following cases: 
    \begin{enumerate}
    \item\label{1:5cn}
        $\rk\,\NS(X)=1$ and there is a quadratic extension $L/\k$ and a birational morphism $\pi\colon X_L\to\p^2_L$ blowing up a point $p=\{p_1,p_2,p_3\}$ of degree $3$ with splitting field $F$ over $\k$ containing $L$, and one of the following cases holds:
        \begin{enumerate}
        \item\label{1:5c1n} $\Gal(F/L)\simeq\Z/3$ and $\Gal(F/\k)\simeq\Z/6$ and the action of $\Aut_\k(X)$ on $\NS(X)$ induces the split exact sequence
            \[1\rightarrow \Aut_L(\p^2,p_1,p_2,p_3)^{\pi\Gal(L/\k)\pi^{-1}}\to \Aut_\k(X) \to \Z/6\rightarrow 1,\]
        \item\label{1:5c2n}
        $\Gal(F/L)\simeq\sym_3$ and $\Gal(F/\k)\simeq\sym_3\times\Z/2$
        and the action of $\Aut_\k(X)$ on $\NS(X)$ induces the split exact sequence
            \[1\rightarrow \Aut_L(\p^2,p_1,p_2,p_3)^{\pi\Gal(L/\k)\pi^{-1}}\to \Aut_\k(X) \to \Z/2\rightarrow 1,\]
        \item\label{1:5c3n}$\Gal(F/L)\simeq\Z/3$ and $\Gal(F/\k)\simeq\sym_3$ acts on the hexagon of $X$ by a rotation of order $3$ and a reflection at an axis through two vertices. The action of $\Aut_\k(X)$ on $\NS(X)$ induces the split exact sequence
                       \[1\rightarrow \Aut_L(\p^2,p_1,p_2,p_3)^{\pi\Gal(L/\k)\pi^{-1}}\to \Aut_\k(X) \to \Z/2\rightarrow 1.\]
        \end{enumerate}
    
    \item\label{1:5bn} $\rk\,\NS(X)\geq2$, $\rk\,\NS(X)^{\Aut_\k(X)}=1$ and $X$ is one of the following:
        \begin{enumerate}
        \item\label{1:51n} $X$ is the blow-up of\, $\p^2$ in the coordinate points, and the action of $\Aut_\k(X)$ on $\NS(X)$ induces the split exact sequence
    \[ 1\rightarrow (\k^*)^2\to \Aut_\k(X)\to\sym_3\times\Z/2\rightarrow 1.\]
        \item\label{1:52n} $X$ is the blow-up of\, $\F_0$ in a point $p=\{(p_1,p_1),(p_2,p_2)\}$ of degree $2$. The action of $\Aut_\k(X)$ on $\NS(X)$ induces the exact sequence
    \[1\rightarrow \Aut_\k(\p^1,p_1,p_2)^2\to \Aut_\k(X)\to \sym_3\times\Z/2\rightarrow 1,\]
    which is split if $\mathrm{char}(\k)\neq2$.
        \item\label{1:53n} $X$ is the blow-up of\, $\p^2$ in a point $p=\{p_1,p_2,p_3\}$ of degree $3$ with splitting field $L$ such that $\Gal(L/\k)\simeq\Z/3$. The action of $\Aut_\k(X)$ on $\NS(X)$ induces the split exact sequence
    \[1\rightarrow \Aut_\k(\p^2,p_1,p_2,p_3)\to \Aut_\k(X) \to \Z/6\rightarrow 1\]
        \item\label{1:54n} $X$ is the blow-up of\, $\p^2$ in a point $p=\{p_1,p_2,p_3\}$ of degree $3$ with splitting field $L$ such that $\Gal(L/\k)\simeq\sym_3$. The action of $\Aut_\k(X)$ on $\NS(X)$ induces the split exact sequence
    \[1\rightarrow \Aut_\k(\p^2,p_1,p_2,p_3)\to \Aut_\k(X)\to \Z/2\rightarrow 1,\]
    where $\Z/2$ is generated by a rotation.
        \end{enumerate}
    \item\label{1:5an} $\rk\,\NS(X)^{\Aut_\k(X)}=2$ and there is a quadratic extension $L/\k$ and a birational morphism $\nu\colon X\to \Ql^L$ contracting two curves onto rational points $p_1,p_2$ or one curve onto a point $\{p_1,p_2\}$ of degree $2$ with splitting field $L'/\k$. The action of $\Aut_\k(X)$ on $\NS(X)$ induces the split exact sequence
    \[1\rightarrow T^{L,L'}(\k)\to\Aut_\k(X)\to \Z/2\times\Z/2\rightarrow1,\]
    where $\nu\Aut_\k(X)\nu^{-1}=\Aut_\k(\Ql^L,\{p_1,p_2\})$ and $T^{L,L'}$ is the subgroup of $\Aut_\k(\Ql^L,p_1,p_2)$ preserving the rulings of $\Ql^L_L$.
    \end{enumerate}
\item\label{1:6n} $\pi\colon X\to\p^1$ is one of the following conic fibrations with
\[\rk\,\NS(X_{\bk}/\p^1)^{\Aut_{\bk}(X,\pi)}=\rk\,\NS(X/\p^1)^{\Aut_\k(X,\pi)}=1\!\!:\]
    \begin{enumerate}
    \item\label{1:61n} $X/\p^1$ is the blow-up of points $p_1,\dots,p_r\in\F_n$, $n\geq2$, contained in a section
    $S_n\subset \F_n$ with $S_n^2=n$.
    The geometric components of the $p_i$ are on pairwise distinct geometric fibres and $\sum_{i=1}^r\deg(p_i)=2n$. There are split exact sequences
            \[
\begin{tikzcd}[link2]
& (T_1/\mu_n)\rtimes\Z/2&\Aut(X)&&\\[-5pt]
& \simeqv&\vert\vert&&\\[-5pt]
1\ar[r] & \Aut(X/\pi_X)\ar[r]&\Aut(X,\pi_X)\ar[r] & \Aut(\p^1,\Delta)\ar[r] &1\\
1\ar[r]& \Aut_\k(X/\pi_X)\ar[r]&\Aut_\k(X,\pi_X)\ar[r] & \Aut_\k(\p^1,\Delta)\ar[r] &1\\[-5pt]
& \simeqv&\vert\vert&&\\[-5pt]
& (\k^*/\mu_n(\k))\rtimes\Z/2&\Aut_\k(X)\rlap{,}&&
\end{tikzcd}
\]
    where $\Delta=\pi(\{p_1,\dots,p_r\})\subset\p^1$, $T_1$ is the split one-dimensional torus and $\mu_n$ its subgroup of $\supth{n}$ roots of unity.
    \item\label{1:62n} There exist quadratic extensions $L$ and $L'$ of $\k$ such that $X/\p^1$ is the blow-up of $\Sl^{L,L'}$ in points $p_1,\dots,p_r\in E$, $r\geq1$. The $p_i$ are all of even degree, their geometric components are on pairwise distinct geometric components of smooth fibres, and each geometric component of $E$ contains half of the geometric components of each $p_i$. There are exact sequences
        \[
\begin{tikzcd}[link2]
& \SlO^{L,L'}\rtimes\Z/2&\Aut(X)&&\\[-5pt]
& \simeqv&\vert\vert&&\\[-5pt]
1\ar[r] & \Aut(X/\pi_X)\ar[r]&\Aut(X,\pi_X)\ar[r] & \Aut(\p^1,\Delta)\ar[r] &1\\
1\ar[r]& \Aut_\k(X/\pi_X)\ar[r]&\Aut_\k(X,\pi_X)\ar[r] & (D^{L,L'}_\k\rtimes\Z/2)\cap\Aut_\k(\p^1,\Delta)\ar[r] &1\\[-5pt]
& \simeqv&\vert\vert&&\\[-5pt]
& \SlO^{L,L'}(\k)\rtimes\Z/2&\Aut_\k(X)&&
\end{tikzcd}
\]
with $\Delta=\pi(\{p_1,\dots,p_r\})\subset\p^1$ and
$\SlO^{L,L'}=\{(a,b)\in T^L\mid ab=1\}$, and
\begin{itemize}
\item if $L,L'$ are $\k$-isomorphic, then $\SlO^{L,L'}(\k)\simeq\{a\in L^*\mid aa^g=1\}$\\
 and $D^{L,L'}_\k\simeq\{\alpha\in k^*\mid \alpha=\lambda\lambda^g,\lambda\in L\}$, where $g$ is the generator of $\Gal(L/\k)$,
\item if $L,L'$ are not $\k$-isomorphic, then
$\SlO^{L,L'}(\k)\simeq\k^*$ and \\
$D^{L,L'}_{\k}\simeq\{\lambda\lambda^{gg'}\in F\mid \lambda\in
K,\lambda\lambda^{g'}=1\}$, where $\k\subset F\subset LL'$ is the
intermediate extension such that $\Gal(F/\k)\simeq\langle
gg'\rangle\subset\Gal(L/\k)\times\Gal(L'/\k)$, where $g,g'$ are the
generators of $\Gal(L/\k),\Gal(L'/\k)$, respectively.
\end{itemize}
        
    \end{enumerate}
\end{enumerate}
\end{thm}

The statements of Theorem~\ref{thm:2} and Corollary~\ref{cor:1} remain correct, and we complete their proofs further below.

\begin{thm}[\textit{cf.~}Theorem~\ref{thm:3}]\label{thm:3n}
Let $\k$ be a perfect field. The conjugacy classes of the maximal
subgroups $\Aut_\k(X)$ of\, $\Bir_\k(\p^2)$  from Theorem~\ref{thm:1n}
are parametrised by
    \begin{itemize}
    \item (\ref{1:1n}), (\ref{1:2n}): one point
    \item (\ref{1:3n}): one point for each $\k$-isomorphism class of quadratic extensions of $\k$
    \item (\ref{1:4n}): one point for each $n\geq2$
    \item (\ref{1:5c1n}) if\, $|\k|\geq3$, one point for any pair of extensions $F\supset L\supset \k$, where $L/\k$ is quadratic and  $F/L$ such that $\Gal(F/L)\simeq\Z/3$ and $\Gal(F/\k)\simeq\Z/6$, up to the following equivalence class: $\k$-isomorphisms $F\simeq F'$ that induce an isomorphism $L\simeq L'$.
    \item (\ref{1:5c2n}): one point for any pair $F\supset L\supset \k$, where $L/\k$ is quadratic and $F/L$ such that $\Gal(F/L)\simeq\sym_3$ and $\Gal(F/\k)\simeq\sym_3\times\Z/2$, up to the following equivalence class: $\k$-isomorphisms $F\simeq F'$ that induce an isomorphism $L\simeq L'$.
    \item (\ref{1:5c3n}): one point for any pair $F\supset L\supset \k$, where $L/\k$ is quadratic and $F/L$ such that $\Gal(F/L)\simeq\Z/3$ and $\Gal(F/\k)\simeq \sym_3$, up to the following equivalence class: $\k$-isomorphisms $F\simeq F'$ that induce an isomorphism $L\simeq L'$.
    \item (\ref{1:51n}): one point if $|\k|\geq3$
    \item (\ref{1:52n}): one point for each $\k$-isomorphism class of quadratic extensions of\, $\k$
    \item (\ref{1:53n}): one point for each $\k$-isomorphism class of Galois extensions $F/\k$ with $\Gal(F/\k)\simeq\Z/3$.
    \item (\ref{1:54n}): one point for any $\k$-isomorphism class of Galois extensions $F/\k$ with $\Gal(F/\k)\simeq\sym_3$.
    \item (\ref{1:61n}): for each $n\geq2$ the set of points $\{p_1,\dots,p_r\}\subset\p^1$ with $\sum_{i=1}^r\deg(p_i)=2n$ up to the action of $\Aut_\k(\p^1)$
    \item (\ref{1:62n}): for each $n\geq1$ and for each pair of $\k$-isomorphism classes of quadratic extensions $(L,L')$, the set of points $\{p_1,\dots,p_r\}\subset\p^1$ of even degree with $\sum_{i=1}^r\deg(p_i)=2n$ up to the action of $D^{L,L'}_{\k}(\k)\rtimes\Z/2$
    \end{itemize}
\end{thm}

\hiddensection{Further corrections}



By the following lemma, whenever we contract a curve onto a point of
degree $2$ in $\Ql^L$ with splitting field~$L$, we can choose the
point conveniently.
It's proof contained contained a gap that we now close.

\begin{lem}[\textit{cf.}~Lemma~\ref{lem:cyclicpoint3}]\label{lem:cyclicpoint3n}
Let $p=\{p_1,p_2,p_3\}$ and $q=\{q_1,q_2,q_3\}$ be points in $\Ql^L$
of degree $3$ such that for any $h\in\Gal(\bk/\k)$ there exists
$\sigma\in\sym_3$ such that $p_i^h=p_{\sigma(i)}$ and
$q_i^h=q_{\sigma(i)}$. Suppose that the geometric components of $p$
$($resp. of $q)$ are in general position on $\Ql^L$. Then
there exists $\alpha\in\Aut_\k(\Ql^L)$ such that $\alpha(p_i)=q_i$
for $i=1,2,3$.
\end{lem}

\begin{proof}
	The assumption on $p$ and $q$ implies that the residue fields of $p$ and $q$ are $\k$-isomorphic; \textit{cf.} \cite[Lemma 3.3]{LS21}.
	Since they are points of degree $3$, $p$ and $q$  therefore have the same splitting field $F/\k$.
	Let $g$ be the generator of the Galois group $\Gal(L/\k)$ of order $2$. We consider the composite field $FL$.
	For $i=1,\ldots,4$, we can write $p_i=(a_i,a_i')$, $q_i=(b_i,b_i')$ with $a_i,a_i',b_i,b_i'\in\p^1_{FL}$.
	By hypothesis, for any $h\in\Gal(\bk/L)$ there exists
	$\sigma\in\sym_3$ such that one has $a_i^h=a_{\sigma(i)}$; similarly for $a_i'$ and $b_i,b_i'$.
	We apply Remark~\ref{rmk:cyclicpoint2-dim1} to the $\Gal(\overline{\k}/FL)$-invariant sets $\{a_i\}_i$ and $\{a_i'\}$ in $\p^1\times\p^1$, respectively $\{b_i\}$ and $\{b_i'\}$, and find a unique $(A,B)\in\PGL_2(FL)\times\PGL_2(FL)$ such that $(A,B)p_i=q_i$ for $i=1,2,3$.
	It remains to see that $(A,B)$ gives an automorphism of $\Ql^L$, that is, $(A,B)$ commutes with the action of $\Gal(LF/\k)$ induced on $\p^1\times\p^1$.
	Let $h\in\Gal(LF/\k)$ and let $\sigma\in\sym_3$ be the permutation induced by $h$.
	We compute $((A,B)p_i)^h= q_i^h=q_{\sigma(i)}$ and $(A,B)p_i^h=(A,B)p_{\sigma(i)}=q_{\sigma(i)}$ for $i=1,2,3$. This concludes the proof since matrices in $\PGL_2(\bar \k)$ are uniquely determined by their action on three points.
\end{proof}


The hexagon of $X_{\bk}$ is $\Gal(\bk/\k)$-invariant. The Galois group $\Gal(\bk/\k)$ acts on the hexagon by symmetries, so we have a group homomorphism 
\[\Gal(\bk/\k)\stackrel{\rho}\to \sym_3\times\Z/2\subseteq \Aut(\NS(X_{\bk})).\]
By the {\em hexagon of}\, $X$ we mean the hexagon of $X_{\bk}$ endowed with its canonical $\Gal(\bk/\k)$-action.
Since the group $\Aut_\k(X)$ acts by symmetries on the hexagon of $X$, it induces a homomorphism
\[\widehat\rho\colon \Aut_\k(X)\to \sym_3\times\Z/2.\]

In our classification of subgroups of $\sym_3\times\Z/2$ there is one missing case, namely $\sym_3$ acting transitively on the edges of a hexagon as (\ref{dp6:10n}) in Figure~\ref{fig:action-dP6n}, that is generated by a rotation of degree $3$ and a reflection at an axis through two vertices; this is because the dihedral group $\sym_3\times\Z/2$ contains two non-conjugate embeddings of $\sym_3$. (The other one is Figure~\ref{fig:action-dP6n}(\ref{dp6:8n}).)

In Figure~\ref{fig:action-dP6n} below, we redraw the hexagons of Figure~\ref{fig:action-dP6} in a slightly different manner, but the numbering remains the same. Here, we choose generators of the respective subgroup of $\sym_3\times\Z/2$ as in the proof of Lemma~\ref{lem:Hexagon} and draw the image of each edge under every generator.
\setcounter{equation}{0}
\begin{figure}[ht]
\begin{minipage}[ht]{.18\textwidth}
\begin{equation}\label{dp6:1n}
\begin{tikzpicture}[scale=.6]
\begin{scope}[every coordinate/.style={shift={(0,3.5)}}]  
\path [c](0,0) pic {hexagon};
\draw[->, thick] [c](-0.2,1) to [bend right=100, looseness=3] (0.2,1); 
\draw[->, thick] [c](1.1,0.7) to [bend right=100,looseness=3] (1.3,0.3); 
\draw[<-, thick,swap] [c](1.1,-0.7) to [bend left=100,looseness=3] (1.3,-0.3); 
\draw[<-, thick] [c](-0.2,-1) to [bend left=100, looseness=3] (0.2,-1); 
\draw[->, thick] [c](-1.1,-0.7) to [bend right=100, looseness=3] (-1.3,-0.3); 
\draw[<-, thick,swap] [c](-1.1,0.7) to [bend left=100, looseness=3] (-1.3,0.3); 
\end{scope}
\end{tikzpicture}
\end{equation}
\end{minipage}
\begin{minipage}[ht]{.18\textwidth}
\begin{equation}\label{dp6:2n}
\begin{tikzpicture}[scale=.6,font=\footnotesize]
\begin{scope}[every coordinate/.style={shift={(0,3.5)}}] 
\path [c](0,0) pic {hexagon}; \draw[<->, thick] (D1) to (D4);
\draw[<->, thick] (D2) to [bend right=15,swap] (D3); \draw[<->,
thick] (D6) to [bend left=15,swap]  (D5);
\end{scope}
\end{tikzpicture}
\end{equation}
\end{minipage}
\begin{minipage}[ht]{.18\textwidth}
\begin{equation}\label{dp6:3n}
\begin{tikzpicture}[scale=.6,font=\footnotesize]
\begin{scope}[every coordinate/.style={shift={(0,3.5)}}] 
\path [c](0,0) pic {hexagon}; \draw[<->, thick] (D2) to [bend
right=-15,swap] (D6); \draw[<->, thick] (D3) to [bend left=-15,swap]
(D5);
\draw[->, thick] [c](-0.2,1) to [bend right=100, looseness=3] (0.2,1); 
\draw[<-, thick] [c](-0.2,-1) to [bend left=100, looseness=3] (0.2,-1); 
\end{scope}
\end{tikzpicture}
\end{equation}
\end{minipage}
\begin{minipage}[ht]{.18\textwidth}
\begin{equation}\label{dp6:4n}
\begin{tikzpicture}[scale=.6,font=\footnotesize]
\begin{scope}[every coordinate/.style={shift={(0,0)}}] 
\path [c](0,0) pic {hexagon}; \draw[<->, thick] (D1) to (D4);
\draw[<->, thick] (D2) to (D5); \draw[<->, thick] (D3) to (D6);
\end{scope}
\end{tikzpicture}
\end{equation}
\end{minipage}
\begin{minipage}[ht]{.18\textwidth}
\begin{equation}\label{dp6:5n}
\begin{tikzpicture}[scale=.6,font=\footnotesize]
\begin{scope}[every coordinate/.style={shift={(0,3.5)}}] 
\path [c](0,0) pic {hexagon};  \draw[<->, thick] (D22) to [bend right=15,swap]
(D3);
\draw[<->, thick] (D55) to [bend right=15,swap] (D6); \draw[<->,
thick, shorten <=.1cm, shorten >=.1cm] (D1) to (D4); \draw[<->,
thick, shorten <=.2cm, shorten >=.2cm] (D2) to (D5); \draw[<->,
thick, shorten <=.2cm, shorten >=.2cm] (D3) to (D6);
\end{scope}
\end{tikzpicture}
\end{equation}
\end{minipage}

\begin{minipage}[ht]{.2\textwidth}
\begin{equation}\label{dp6:6n}
\begin{tikzpicture}[scale=.6,font=\footnotesize]
\begin{scope}[every coordinate/.style={shift={(0,3.5)}}] 
\path [c](0,0) pic {hexagon}; \draw[->, thick] (D11) to [bend
right=15,swap] (D3); \draw[->, thick] (D22) to [bend right=15,swap]
(D4); \draw[->, thick] (D33) to [bend right=15,swap] (D5); \draw[->,
thick] (D44) to [bend right=15,swap] (D6); \draw[->, thick] (D55) to
[bend right=15,swap] (D1); \draw[->, thick] (D66) to [bend
right=15,swap] (D2);
\end{scope}
\end{tikzpicture}
\end{equation}
\end{minipage}
\begin{minipage}[ht]{.23\textwidth}
\begin{equation}\label{dp6:7n}
\begin{tikzpicture}[scale=.6,font=\footnotesize]
\begin{scope}[every coordinate/.style={shift={(0,3.5)}}] 
\path [c](0,0) pic {hexagon}; \draw[->, thick] (D11) to [bend
right=15,swap] (D2); \draw[->, thick] (D22) to [bend right=15,swap]
(D3); \draw[->, thick] (D33) to [bend right=15,swap] (D4); \draw[->,
thick] (D44) to [bend right=15,swap] (D5); \draw[->, thick] (D55) to
[bend right=15,swap] (D6); \draw[->, thick] (D66) to [bend
right=15,swap] (D1);
\end{scope}
\end{tikzpicture}
\end{equation}
\end{minipage}
\begin{minipage}[ht]{.23\textwidth}
\begin{equation}\label{dp6:8n}
\begin{tikzpicture}[scale=.6,font=\footnotesize]
\begin{scope}[every coordinate/.style={shift={(0,3.5)}}] 
\path [c](0,0) pic {hexagon};
\draw[->, thick] (D11) to [bend
right=15,swap] (D3); \draw[->, thick] (D2) to [bend right=15,swap]
(D4); \draw[<->, thick] (D33) to [bend right=15,swap] (D5); \draw[->,
thick] (D44) to [bend right=15,swap] (D6); \draw[->, thick] (D55) to
[bend right=15,swap] (D1); \draw[<->, thick] (D666) to [bend
right=15,swap] (D2222);
\end{scope}
\end{tikzpicture}
\end{equation}
\end{minipage}
\begin{minipage}[ht]{.23\textwidth}
\begin{equation}\label{dp6:9n}
\begin{tikzpicture}[scale=.6,font=\footnotesize]
\begin{scope}[every coordinate/.style={shift={(0,3.5)}}] 
\path [c](0,0) pic {hexagon};  \draw[->, thick] (D11) to [bend
right=15,swap] (D2); \draw[<->, thick] (D22) to [bend right=15,swap]
(D3); \draw[->, thick] (D33) to [bend right=15,swap] (D4); \draw[->, thick] (D44) to [bend right=15,swap] (D5); \draw[<->, thick] (D55) to
[bend right=15,swap] (D6); \draw[->, thick] (D66) to [bend
right=15,swap] (D1);
\draw[<->, thick,shorten <=.2cm, shorten >=.2cm] (D1) to [bend
right=0,swap] (D4);
\end{scope}
\end{tikzpicture}
\end{equation}
\end{minipage}
\begin{minipage}[ht]{.23\textwidth}
\begin{equation}\label{dp6:10n} 
\begin{tikzpicture}[scale=6,font=\footnotesize]
\begin{scope}[every coordinate/.style={shift={(0,-3.5)}}] 
\path [c](0,0) pic {hexagon}; \draw[->, thick,shorten <=.2cm, shorten >=.2cm] (D1) to [bend
right=15,swap] (D3); \draw[->, thick,shorten <=.2cm, shorten >=.2cm] (D2) to [bend right=15,swap]
(D4); \draw[->, thick,shorten <=.2cm, shorten >=.2cm] (D3) to [bend right=20,swap] (D5);
\draw[->, thick,shorten <=.2cm, shorten >=.2cm] (D4) to [bend right=15,swap] (D6); \draw[->,
thick,shorten <=.2cm, shorten >=.2cm] (D5) to [bend right=15,swap] (D1); \draw[->, thick,shorten <=.2cm, shorten >=.2cm] (D6)
to [bend right=20,swap] (D2);
\draw[<->, thick] (D1) to[bend right=0] (D4);
\draw[<->, thick] (D2) to[bend right=0] (D3);
\draw[<->, thick] (D5) to[bend right=0] (D6);
\end{scope}
\end{tikzpicture}
\end{equation}
\end{minipage}

\caption{The $\Gal(\bk/\k)$-actions on the hexagon of a del
Pezzo surface of degree $6$.}\label{fig:action-dP6n}
\end{figure}

\begin{lem}\label{lem:Hexagon}
	The action of $\rho(\Gal(\bar\k/\k))$ on the hexagon of a del Pezzo surface of degree $6$ is as in Figure~\ref{fig:action-dP6n}.
\end{lem}

\begin{proof}
The dihedral group $D_6=\langle r,s \mid r^6=s^2=\id, srs=r^{-1}\rangle$ equals $\sym_3\times\Z/2=\langle r^2,s\rangle\times \langle r^3\rangle$, where $r$ is a rotation of order $6$, as in Figure~\ref{fig:action-dP6n}(\ref{dp6:7n}) and $s$ is a reflection as in Figure~\ref{fig:action-dP6n}(\ref{dp6:2n}).
	Writing $t=r^3s$ for a reflection as in Figure~\ref{fig:action-dP6n}(\ref{dp6:3n}), the action of $\rho(\Gal(\bar\k/\k))$ is one of the following, up to conjugation. (For convenience, we give references to the lemmas that deal with the respective cases.)
	\begin{enumerate}
		\item The trivial subgroup as in Figure~\ref{fig:action-dP6n}(\ref{dp6:1n}), see Lemma~\ref{prop:DP1} 
		\item $\langle s \rangle = \Z/2$ as in Figure~\ref{fig:action-dP6n}(\ref{dp6:2n}), see Lemma~\ref{prop:DP2} 
		\item  $\langle t \rangle = \Z/2$ as in Figure~\ref{fig:action-dP6n}(\ref{dp6:3n}), see Lemma~\ref{prop:DP3}
		\item $\langle r^3 \rangle = \Z/2$ as in Figure~\ref{fig:action-dP6n}(\ref{dp6:4n}), see Lemma~\ref{prop:DP4}  
		\item $\langle r^3, s \rangle = \Z/2\times\Z/2$ as in Figure~\ref{fig:action-dP6n}(\ref{dp6:5n}), see Lemma~\ref{prop:DP9} 
		\item $\langle r^2 \rangle = \Z/3$ as in Figure~\ref{fig:action-dP6n}(\ref{dp6:6n}), see Lemma~\ref{prop:DP5}
		\item $\langle r \rangle = \Z/6$ as in Figure~\ref{fig:action-dP6n}(\ref{dp6:7n}), see Lemma~\ref{prop:DP6}  and Remark~\ref{rmk:correction Lemma 4.6 and Lemma 4.7}
		\item $\langle r^2,t \rangle = \sym_3$ as in Figure~\ref{fig:action-dP6n}(\ref{dp6:8n}), see Lemma~\ref{prop:DP7} 
		\item $\langle r,s \rangle = \sym_3\times \Z/2$ as in Figure~\ref{fig:action-dP6}(\ref{dp6:9}), see  Lemma~\ref{prop:DP8}  and Remark~\ref{rmk:correction Lemma 4.6 and Lemma 4.7}
		\item $\langle r^2,s \rangle = \sym_3$ as in Figure~\ref{fig:action-dP6n}(\ref{dp6:10n}), see Lemma~\ref{lem:DP10} below.
	\end{enumerate}
\end{proof}

We now discuss the cases where $\rho(\Gal(\bar\k/\k))$ acts transitively on the six edges of the hexagon, that is, (\ref{dp6:7n}), (\ref{dp6:9n}) and the missing case (\ref{dp6:10n}) in Figure~\ref{fig:action-dP6n}.

\begin{rmk}\label{rmk:correction Lemma 4.6 and Lemma 4.7}
In Lemma~\ref{prop:DP6}(1) and Lemma~\ref{prop:DP8}(1), which describe the situation of (\ref{dp6:7n}) and (\ref{dp6:9n}) in Figure~\ref{fig:action-dP6n}, the correct statement is that there is a quadratic extension $L/\k$ and a Galois extension $F/L$ with respectively
\begin{enumerate}
	\item $\Gal(F/L)\simeq\Z/3$ and $\Gal(F/\k)\simeq\Z/6$ for Lemma~\ref{prop:DP6}, and
	\item $\Gal(F/L)\simeq\sym_3$ and $\Gal(F/\k)\simeq\sym_3\times\Z/2$ for Lemma~\ref{prop:DP8}.
\end{enumerate}
This comes out of the analysis of the del Pezzo surfaces in Figure~\ref{fig:action-dP6n}(\ref{dp6:7n}) and (\ref{dp6:9n}) in Lemmas~\ref{prop:DP6} and~\ref{prop:DP8}  and is stated incorrectly in the statements of the two lemmas.

Similarly, the correct assumptions in Example~\ref{ex:4.8}  should be $F/L$ as above. In this case it holds that $\Gal(FL/L)\simeq \Gal(F/\k)$, as claimed in said example. In fact, the missing case is $\Gal(F/L)\simeq\Z/3$ and $\Gal(F/\k)=\sym_3$.
\end{rmk}

\begin{lem}\label{lem:DP10}
Let $X$ be a del Pezzo surface of degree $6$ such that $\rho(\Gal(\bk/\k))\simeq\sym_3$ as indicated in Figure~\ref{fig:action-dP6n}(\ref{dp6:10n}). Then $X\to\ast$ is a Mori fibre space, and the following hold:
\begin{enumerate}
\item\label{DP10:1}
There exists a quadratic field extension $L/\k$ and a point $p=\{p_1,p_2,p_3\}$ in $\p^2_L$ of degree $3$ such that $X_L$ is isomorphic to the blow-up of\, $\p^2_L$ in $p$. Moreover, there is a cyclic extension $F/L$ of degree $3$ such that each $(-1)$-curve in the hexagon of $X$ is defined over $F$.
\item\label{DP10:2} Two such surfaces $X$ and $X'$ are isomorphic if and only if there exists an isomorphism $F\to F'$ over $\k$ that sends $L$ onto $L'$.
\item\label{DP10:3} If\, $X_L$ is $L$-rational, then the action of $\Aut_\k(X)$ on the hexagon of $X$ induces a split exact sequence
\[1\rightarrow(\pi^{-1}\Aut_L(\p^2_L,p_1,p_2,p_3)\pi)^{g}\rightarrow\Aut_\k(X)\stackrel{\hat{\rho}}\rightarrow \Z/2\rightarrow1,\]
where $\Z/2$ is generated by a rotation and $g$ is the generator of\, $\Gal(L/\k)$.
\end{enumerate}
\end{lem}
\begin{proof}
Every $(-1)$-curve in the hexagon of $X$ is contained in the same $\Gal(\bk/\k)$-orbit; hence $X$ is a Mori fibre space.

(\ref{DP10:1})\&(\ref{DP10:3}) The group $\Aut_\k(X)$ acts by symmetries on the hexagon of $X$. The only element of $\sym_3\times\Z/2$ that commutes with the action of $\Gal(\bk/\k)$ is the rotation of order $2$, so $\hat\rho(\Aut_\k(X))\subset\Z/2$. Let us show that $\hat\rho(\Aut_\k(X))=\Z/2$ and that $\hat\rho$ has a section.
Let $F/\k$ be the splitting field of a $(-1)$-curve in the hexagon of $X$, \textit{i.e.}\ the smallest normal field extension of $\k$ over which the curve is defined. Then $\Gal(F/\k)\simeq\sym_3$ with action on the hexagon as in Figure~\ref{fig:action-dP6n}(\ref{dp6:10n}).
Let $r\in\Gal(F/\k)$ be the rotation of order $3$ indicated in Figure~\ref{fig:action-dP6n}(\ref{dp6:10n}) and  $L:=F^{r} \subset F$ the subfield of $F$ fixed by $r$.
Then $\Gal(F/L)=\langle r\rangle\simeq\Z/3$ and $\Gal(L/\k)\simeq\Gal(F/\k)/\Gal(F/L)\simeq\Z/2$ is generated by the reflection in Figure~\ref{fig:action-dP6n}(\ref{dp6:10n}). The action of $\Gal(\bk/L)$ on the hexagon of $X_L$ factors through $\Gal(F/L)$,
forming two $\Gal(\bk/L)$-orbits of pairwise disjoint $(-1)$-curves of size $3$. This gives (\ref{DP10:1}).

If $X_L$ is $L$-rational, the hexagon of $X_L$ is as in Lemma~\ref{prop:DP5}  and the contraction of one of the $\Gal(\bk/L)$-orbits in the hexagon of $X_L$ onto a point $p$ of degree $3$ is a birational morphism $\pi\colon X_L\to \p^2_L$.
The kernel of $\widehat\rho$ is isomorphic to $(\pi^{-1}\Aut_L(\p^2_L,p_1,p_2,p_3)\pi)^{g}$, where $p_1,p_2,p_3\in\p^2(L)$ are the geometric components of $p$ and $g$ is the generator of $\Gal(L/\k)\simeq\Z/2$.
The only non-trivial element of $\sym_3\times\Z/2$ commuting with $\rho(\Gal(\bk/\k))$ is the rotation of order $2$; hence $\widehat\rho(\Aut_\k(X))\subset\Z/2$.
There exists a quadratic involution $\varphi_p\in\Bir_L(\p^2)$ with base-point $p$ that induces a
rotation of order $2$
on the hexagon of $X_L$ (see Lemma~\ref{prop:DP5}\eqref{DP5:4}).
Since $X_L$ is $L$-rational, it has an $L$-rational point $q$. It is not contained in the hexagon of $X_L$, and we assume, moreover, that $\varphi_p$ fixes $\pi(q)$.
It remains to check that $\psi:=\pi^{-1}\varphi_p\pi\in\Aut(X_L)$ is defined over $\k$. The automorphism $\psi^{-1}g\psi g$ of $X_L$ is conjugate by $\pi$ to an automorphism of $\p^2_L$ fixing $\pi(q)$ and each geometric component of $p$, and is hence the identity map.
Thus the involution $\varphi_p$ lifts to a $\k$-automorphism of $X$. It acts like a rotation of order $2$ on the hexagon of $X$, and thus $\hat\rho(\Aut_\k(X))=\Z/2$ and the sequence splits.

(\ref{DP10:2})
If there is an automorphism $\tau\colon F\to F'$ that sends $L$ onto $L'$ and fixes $\k$, we can identify $L$ and $L'$. Then the surfaces $X_L$ and $X'_{L'}$ are $L$-isomorphic by Lemma~\ref{prop:DP5}. Since $\tau$ fixes $\k$, it induces a $\k$-isomorphism of $X$ and $X'$ by the above construction.
On the other hand, suppose that $X$ and $X'$ are $\k$-isomorphic. Then the smallest normal field extensions $F$ and $F'$ over which all $(-1)$-curves of $X_{\bk}$ and $X'_{\bk}$ are defined are $\k$-isomorphic. This isomorphism sends the fixed field $L=F^r$ onto the fixed field $L'=F'^r$.
\end{proof}

\begin{ex}[Construction of rational del Pezzo surfaces as in Figure~\ref{fig:action-dP6n}(\ref{dp6:7n}),(\ref{dp6:9n}),(\ref{dp6:10n})]
(See also Example~\ref{ex:4.8}.)
Let $q=\{q_1,q_2\}$ in $\p^2$ be a point of degree $2$, with splitting field $L/\k$ being a quadratic extension, and let $p=\{p_1,p_2,p_3\}$ in $\p^2$ be a point of degree $3$, with splitting field $F/\k$ (which is an extension of degree $3$ or~$6$). Assume that the components of $p$ and $q$ are in general position, that is, no three of the five geometric components are collinear. Denote by $D\subset\p^2$ the conic passing through the five geometric components.

Blowing up $q$ and contracting the line passing through $q$ gives a $k$-birational map $\p^2\dashrightarrow \Ql^L$, where $\Ql^L$ is a del Pezzo surface of degree $8$ as in Definition~\ref{def:Q}. The image of $p$ in $\Ql^L$, again denoted by $p$, is in general position, and the strict transform of $D$ is the diagonal passing through $p$. Blowing up $p$ and then contracting the strict transform of $D$ gives a $k$-birational map $\Ql^L\dashrightarrow X$, where $X$ is a del Pezzo surface of degree $6$, and the action on the hexagon of $X$ is one of the following:
\begin{enumerate}
	\item If $\Gal(F/\k)\simeq\Z/3\Z$ and so $\Gal(LF/\k)\simeq \Z/6\Z$, then the hexagon of $X$ is as in Figure~\ref{fig:action-dP6n}(\ref{dp6:7n}). (For example, take $\k=\F_q$, $L=\F_{q^2}$, $F=\F_{q^3}$, or $\k=\Q(\zeta)$, where $\zeta$ is a primitive third root of $1$, and $L=\k[\sqrt[2]{3}]$ and $F=\k[\sqrt[3]{2}]$.)
	\item If $\Gal(F/\k)\simeq \sym_3$ and $L\subset F$, hence $\Gal(LF/\k)\simeq \sym_3$, then the hexagon of $X$ is as in Figure~\ref{fig:action-dP6n}(\ref{dp6:10n}). (For example, take $\k=\Q$, $L=\Q(\zeta)$, and $F=\Q[\sqrt[3]{2},\zeta]$.)
	\item If $\Gal(F/\k)\simeq \sym_3$ and $\Gal(LF/\k)\simeq \sym_3\times\Z/2\Z$, then the hexagon of $X$ is as in Figure~\ref{fig:action-dP6n}(\ref{dp6:9n}). (For example, take $\k=\Q$, $L=\Q[i]$, $F=\Q[\sqrt[3]{2},\zeta]$.)
\end{enumerate}
\end{ex}


\begin{proof}[Proof of Theorem~\ref{thm:1n}]

The proof remains the same as that of Theorem~\ref{thm:1}, but we have to add the missing case (\ref{1:5c3n}) from the previous section: if the $\Gal(\bk/\k)$-action is as in Figure~\ref{fig:action-dP6n}(\ref{dp6:10n}), Lemma~\ref{lem:DP10} implies that $(X,\Aut(X))$ is as in Theorem~\ref{thm:1n}(\ref{1:5c3n}).
\end{proof}


We now check whether the missing del Pezzo surface case is $\Aut_\k(X)$-birationally rigid or superrigid or neither of them.
Let $X$ be as in Figure~\ref{fig:action-dP6n}(\ref{dp6:10n}).
Then the field $\k$ is infinite because any finite field extension of a finite field is cyclic.
Suppose that $X$ is $\k$-rational and pick $q\in X(\k)$. Let $\pi\colon X_L\to \p^2_L$ be the contraction from Lemma~\ref{lem:DP10}.
As in Lemma~\ref{lem:fixedpoints}, one shows that the map
\[
H:=(\pi^{-1}\Aut_L(\p^2_L,p_1,p_2,p_3)\pi)^{g}\to X(\k),\quad \alpha\longmapsto \pi^{-1}\alpha\pi(q)
\]
is a bijection.
Since $\k$ is infinite, the group $H$ is infinite, and it acts faithfully on the $\bk$-points of $X$ outside the hexagon.

\begin{lem}\label{lem:rigidity-DP10}
Let $X$ be a del Pezzo surface of degree $6$ as in Figure~\ref{fig:action-dP6n}(\ref{dp6:10n}). Then there are no $\Aut_\k(X)$-orbits with at most $5$ geometric components.
In particular, there are no $\Aut_\k(X)$-equivariant Sarkisov links starting from $X$.
\end{lem}
\begin{proof}
The group $\Aut_\k(X)$ acts transitively on the $(-1)$-curves of the hexagon of $X$. By the above remark, $H$ has no orbits with at most $5$ geometric components outside the hexagon of $X$. Therefore, there are no $\Aut_\k(X)$-equivariant Sarkisov links starting from $X$.
\end{proof}

While the statement of Theorem~\ref{thm:2} remains correct, we need to complete its proof with the additional case in the classification of del Pezzo surfaces of degree $6$.

\begin{proof}[Completion of the proof of Theorem~\ref{thm:2}]
The proof remains the same; however, we have to add the case of the del Pezzo surface $X$ of degree $6$ in Figure~\ref{fig:action-dP6n}. By the above Lemma~\ref{lem:rigidity-DP10}, $X$ is $\Aut_\k(X)$-birationally superrigid.
\end{proof}

The proof of Corollary~\ref{cor:1} remains correct because the statement of Theorem~\ref{thm:2} is correct.

\begin{proof}[Proof of Theorem~\ref{thm:3n}]
The proof remains the same as that of Theorem~\ref{thm:3} except for case (\ref{1:5c}). The correction for the parameter in class (\ref{1:5c1n}) and (\ref{1:5c2n}) follows from Remark~\ref{rmk:correction Lemma 4.6 and Lemma 4.7}.
The parameter for the additional class (\ref{1:5c3n}) follows from Lemma~\ref{lem:DP10}.
\end{proof}


Finally, the new class (\ref{1:5c3n}) of del Pezzo surfaces $X$ of degree $6$ are Mori fibre spaces by Lemma~\ref{lem:DP10}. As Remark~\ref{rmk:homomorphism} explains, this means that the homomorphism
\[
\Psi\colon \mathrm{BirMori}(\p^2)\to (\bigoplus_{\chi\in M(\F_1)}\Z/2)\bigast_{C\in J_5}(\bigoplus_{\chi\in M(C)}\Z/2)\ast(\bigast_{C\in J_6}\bigoplus_{\chi\in M(C)}\Z/2)
\]
sends the conjugacy class of $\Aut_\k(X)$ onto zero. Therefore, the statements and proofs of Theorem~\ref{thm:julia} and Proposition~\ref{thm:4}  remain correct.

\subsection*{Acknowledgments}

We'd like to thank Aurore Boitrel for pointing out the missing case in the classification.

\renewcommand{\addcontentsline}[3]{}

\end{document}